\documentstyle[11pt]{article}

\begin{document}
\noindent
\begin{center}
  {\LARGE A Homotopy Theory of Orbispaces}
  \end{center}

  \noindent
  \begin{center}

                    Weimin Chen\footnote{partially 
supported by the National Science Foundation}\\[5pt]
      Department of Mathematics, SUNY-Stony Brook,\\ 
                 Stony Brook, NY 11794 \\[5pt]
              e-mail: wechen@math.sunysb.edu\\[5pt]
              \end{center}

              \def \U{{\cal U}}
              \def \V{{\cal V}}
              \def \E{{\cal E}}
              \def \F{{\cal F}}
              \def \L{{\cal L}}
              \def \O{{\cal O}}
              \def \C{{\bf C}}
              \def \R{{\bf R}}
              \def \Z{{\bf Z}}
              \def \D{{\bf D}}
              \def \N{{\cal N}}

\abstract An orbifold is a singular space which is locally modeled on
the quotient of a smooth manifold by a smooth action of a finite
group.  It appears naturally in geometry and topology when group
actions on manifolds are involved and the stabilizer of each fixed
point is finite. The concept of an orbifold was first introduced by
Satake under the name ``$V$-manifold'' in a paper where he also
extended the basic differential geometry to his newly defined singular
spaces (cf. \cite{Sa}).  The local structure of an orbifold -- being the
quotient of a smooth manifold by a finite group action -- was merely
used as some ``generalized smooth structure''. A different aspect of
the local structure was later recognized by Thurston, who gave the
name ``orbifold'' and introduced an important concept -- the
fundamental group of an orbifold (cf. \cite{Th}).

In 1985, physicists Dixon, Harvey, Vafa and Witten studied string
theories on Calabi-Yau orbifolds (cf. \cite{DHVW}).  An interesting
discovery in their paper was the prediction that a certain physicist's
Euler number of the orbifold must be equal to the Euler number of any
of its crepant resolutions. This was soon related to the so called
McKay correspondence in mathematics (cf. \cite{McK}). Later developments
include orbifold or stringy Hodge numbers (cf. \cite{V, Z, BD}), 
mirror symmetry of
Calabi-Yau orbifolds (cf. \cite{Ro}), and most recently the Gromov-Witten
invariants of symplectic orbifolds (cf. \cite{CR1, CR2}). 
One common feature
of these studies is that certain contributions from singularities,
which are called ``twisted sectors'' in physics, have to be properly
incorporated. This is called the ``stringy aspect'' of an orbifold
(cf. \cite{R}).

This paper makes an effort to understand the stringy aspect of
orbifolds in the realm of ``traditional mathematics''. Surprisingly, we
were led to a refinement of Thurston's discovery!

\tableofcontents

\section{Introduction}

\hspace{5mm} The purpose of this paper is to formally introduce the
category of orbispaces, and to construct a homotopy theory, to which
the basic aspects of the (ordinary) homotopy theory of topological
spaces such as relative homotopy groups and the associated exact
sequence, the covering theory, the fibration theory, etc. are
extended. The category of orbispaces may be regarded as a
natural  extension of
the equivariant category (i.e. the category of $G$-spaces where 
morphisms are equivariant maps), and so
does the corresponding homotopy theory. More precisely, each (locally
connected) $G$-space canonically determines an orbispace whose
underlying topological space is the orbit space of the given
$G$-space, and the constructed homotopy groups of the orbispace are
naturally isomorphic to the corresponding homotopy groups of the Borel
space (cf. \cite{Bo}) of the $G$-space. In general, an orbispace is
just a (locally connected) topological space equipped with a
compatible system of local $G$-space structures, and a morphism
between two orbispaces is just an equivalence class of compatible
systems of local equivariant maps. The notion of an orbispace is not
new. For example, there were relevant definitions by Haefliger
(cf. \cite{Ha2}) and Kontsevich (cf. \cite{Kon}), where the orbispace
structure was described by an (\'{e}tale) topological groupoid. The
definition of orbispaces in this paper is formulated along the lines
of the definition of $V$-manifolds (or orbifolds in Thurston's
terminology) given by Satake (cf. \cite{Sa}). In particular, an
orbifold is naturally an orbispace in the sense of this
paper. Our definition of orbispaces is slightly more general than
those by Haefliger or Kontsevich in the sense that the local group
actions in our definition need not to be discrete. We shall call an
orbispace \'{e}tale if all of the local group actions are discrete. 

The homotopic-homological invariants of a topological groupoid were
usually defined as the ordinary invariants of the corresponding
classifying space of the topological groupoid. For example, the
equivariant (co)homology of a $G$-space is defined as the ordinary
(co)homology of the associated Borel space, and the homotopy groups
and (co)homology groups of an orbifold were defined by Haefliger as
the corresponding ordinary groups of the classifying space of the
topological groupoid associated to the orbifold (cf. \cite{Ha1}).  (It
is shown that our definition of homotopy groups of orbifolds is
equivalent to Haefliger's.) The novelty of this paper is the
introduction of the notion of morphisms between orbispaces with which
they form a category, and to propose to construct
algebraic-topological invariants of the category of orbispaces from
the spaces of morphisms, rather than from the classifying spaces of
the topological groupoids. The homotopy theory of orbispaces
constructed in this paper is obtained by extending the usual
construction on the based loop spaces of a topological space to the
based loop spaces of an orbispace (i.e. the space of based morphisms
from $S^1$ into the orbispace). Another interesting point of this
paper is the introduction of the notion of free loop space of an
orbifold, which is the space of all smooth morphisms (in the sense of
orbispaces) from $S^1$ into the orbifold. The free loop space of an
orbifold has a natural pre-Hilbert orbifold structure, and generalizes
the notion of twisted loop spaces introduced by physicists in the
study of string theories on Calabi-Yau orbifolds, where the orbifolds
are quotients of smooth manifolds by finite (or more general discrete)
groups.

\vspace{1.5mm}

String theories on orbifolds were first introduced by Dixon, Harvey,
Vafa and Witten in \cite{DHVW}, in which the physicists analyzed
string propagation on a compact Calabi-Yau manifold $Y$ equipped with
a finite group action of $G$ preserving the Calabi-Yau structure of
$Y$ (more general $Y$ can be non-compact and $G$ be discrete with the
quotient $X=Y/G$ being a compact orbifold). For the purpose of
symmetry breaking, it was necessary to not only consider strings
$y(t)$ satisfying periodic boundary conditions but also boundary
conditions periodic up to the action of $G$:
$$
y(t+2\pi)=g\cdot y(t), \hspace{3mm} \mbox{ for some } g\in G. \leqno (1.1)
$$
These more general boundary conditions are the so-called ``twisted
boundary conditions''. One consistent framework proposed in
\cite{DHVW} for analyzing propagation of such strings on $Y$ is to
consider the (closed) string theory on the quotient orbifold $X=Y/G$. 

In string theories on a smooth manifold $M$, the Euler characteristic
$\chi(M)$ of $M$ is interpreted as twice of the ``number of
generations'' in the physical theory. In calculating the number of
generations in the orbifold case, the modular invariance of closed
string theories also requires the inclusion of contributions from the
string sectors with twisted boundary conditions. Using path integral
method, the physicists obtained the following formula for the ``Euler
characteristic'' of the orbifold $X=Y/G$ as twice of the number of
generations in the physical theory on $X$: 
$$
\chi_{orb}(X)=\frac{1}{|G|}\sum_{gh=hg}\chi(g,h), \leqno (1.2)
$$
where $\chi(g,h)$ is the Euler characteristic of the common
fixed-point set of $g$ and $h$ in $Y$. This was later reformulated by
Hirzebruch and H\"{o}fer (cf. \cite{HH}) as 
$$
\chi_{orb}(X)=\chi(X)+\sum_{(g),g\neq 1}\chi(Y^g/C(g)), \leqno (1.3)
$$
where $Y^g$ is the fixed-point set of $g$ in $Y$, $(g)$ stands for the
conjugacy class of $g$ in $G$ and $C(g)$ is the centralizer of $g$ in $G$.

One of the advantages of introducing string theories on orbifolds is
that string propagation on an orbifold may be regarded as an arbitrarily good
approximation to string propagation on any of the smooth resolutions
of the orbifold. This relation has the following remarkable
consequence: the Euler characteristic of a Calabi-Yau orbifold $X$ as
defined in $(1.2)$ or equivalently in $(1.3)$ must be equal to the Euler
characteristic of any smooth crepant resolution $\widehat{X}$ of $X$
$$
\chi_{orb}(X)=\chi(\widehat{X}). \leqno (1.4)
$$

The physicists' prediction $(1.4)$ was soon related to the so-called
``McKay correspondence'' in mathematics (cf. \cite{McK}). A weak
version of McKay correspondence may be stated as follows:\\
Let $G\subset SL(n,\C)$ be a finite subgroup, and
$\pi:Y\rightarrow X=\C^n/G$ be a crepant resolution. There exist
``natural'' bijections between conjugacy classes of $G$ and basis of
$H_\ast(Y,\Z)$ (cf. \cite{Re}).

There have been subsequent developments since \cite{DHVW}, both in
mathematics and physics. Very recently, the Gromov-Witten invariants
of symplectic or projective orbifolds were constructed in 
\cite{CR1, CR2}, where certain cohomologies from the singular set of the
orbifold were involved in a very essential way. More concretely, one
can associate each orbifold $X$ a space
$$
\widetilde{X}:=\{(p,(g)_{G_p})|p\in X, g\in G_p \}, \leqno (1.5)
$$
where $G_p$ is the ``stabilizer'' group at $p$ and $(g)_{G_p}$ is the
conjugacy class of an element $g$ of $G_p$. The space $\widetilde{X}$
was first introduced by Kawasaki in \cite{Ka} in the study of index
theory over $V$-manifolds. In the case when $X=Y/G$,
$$
\widetilde{X}=\sqcup_{(g)} Y^g/C(g)
=X\sqcup\sqcup_{\{(g),g\neq 1\}} Y^g/C(g), \leqno (1.6)
$$ 
which appears in the formulation of the physicist's Euler
characteristic of the orbifold $X$ in $(1.3)$. The main results of
\cite{CR1, CR2} can be summarized as follows: with suitable degree
shifting on the cohomology groups of the components of $\widetilde{X}$
coming from singularity of the orbifold $X$, there is a new cup
product on the total cohomology of the space $\widetilde{X}$, and the
quantum cohomology of $X$ is a (quantum) deformation of the new
cohomology ring (which is the total cohomology of $\widetilde{X}$ with
suitable degree shifting, rather than the total ordinary cohomology
ring of $X$). The components of $\widetilde{X}$ coming from
singularity of $X$ are called ``twisted sectors'' in \cite{CR1, CR2}.

\vspace{1.5mm}

The new mathematics of orbifolds, in which the space $\widetilde{X}$
instead of the orbifold $X$ itself seems to be a more natural object
to study, is given a name as the ``stringy aspect'' of orbifolds,
reflecting its origin from string theory (cf. \cite{R}). Despite its
rich structures and great influence on mathematics, string theory
still remains mysterious and uncertain to most of mathematicians.  The
question of how much of the ``stringy aspect'' of orbifolds can be
understood in the realm of classical mathematics without going into
the ``quantum level'' is the main motivation of this article. With the
formulation of the category of orbispaces and by regarding an orbifold
as an orbispace, we introduced the free loop space of an orbifold, 
which is defined as the space of all morphisms (in the sense of
orbispaces) from $S^1$ into the orbifold. The free loop space of an
orbifold has the following features:\\
1). It generalizes the space of strings satisfying twisted boundary
conditions $(1.1)$ when the orbifold is a quotient $Y/G$.\\
2). The space $\widetilde{X}$ defined in $(1.5)$ embeds in the free
loop space of the orbifold $X$ as the fixed-point set of the canonical
$S^1$-action on the free loop space defined by rotating the domain of
each loop. \\
3). It has a nice geometric structure, being a pre-Hilbert orbifold. 

The first two features of the free loop space seem to provide some
evidence that the category of orbispaces might serve as an appropriate
mathematical framework for the string theories on orbifolds in
physics. On the other hand, if one regards the category of orbispaces
as a natural extension of the equivariant category, it is very
tempting to see the appearance of the ``twisted sector'' components of
the space $\widetilde{X}$ defined in $(1.5)$ as a result of some kind
of generalizations of the localization theorems in the equivariant
category (cf. \cite{Se1, Qu}). This is consistent with 
the result of Atiyah and Segal
(cf. \cite{AS, BC1}) which interpreted the Euler characteristic
$\chi_{orb}(X)$ defined by $(1.2)$ as the Euler characteristic of the
equivariant K-theory of the $G$-space $Y$ when the orbifold $X$ is a
global quotient $Y/G$. Therefore, it seems to be a very interesting
problem to develop a general framework for cohomology theories of
orbispaces, of which certain localization theorems still hold. In
particular, it would be interesting to work out the localization
theorem (if there is any) for the K-theory of orbispaces. This
viewpoint would provide an appropriate interpretation of the so-called
``contributions from singularities'', and the physicist's prediction
$(1.4)$ for a general Calabi-Yau orbifold would still be consistent
with the general philosophy of McKay correspondence --- relating
``representation theory'' to ``resolution of singularities''.

The pre-Hilbert orbifold structure on the free loop space of an
orbifold makes it possible to extend some geometrical-analytical
constructions on the free loop space of a smooth manifold to the free
loop space of an orbifold. Especially what we have in mind is Witten's 
interpretation of elliptic genera and the proof of rigidity
(cf. \cite{W, BT, T}). If an analogous construction can
be carried out in the orbifold case, we would have a general
definition of ``orbifold elliptic genera'' and some kind of rigidity
property of orbifold elliptic genera\footnote{relevant work has been
done recently, cf. \cite{Liu}}. One may even ask if there is a 
corresponding elliptic cohomology theory for the category of
orbispaces (cf. \cite{Se2, De}). Of course, this would first
require to establish a general framework for (co)homology theories of
orbispaces. 

\vspace{1.5mm}

Next we turn to a particular viewpoint of the paper \cite{DHVW} of
Dixon, Harvey, Vafa and Witten, which seems to deserve more
attention. Recall that the physicists were trying to analyze string
propagation in an equivariant setting, i.e., the propagation of
strings on a Calabi-Yau manifold equipped with a compatible action of
a finite group. The strings are required to satisfy the twisted
boundary conditions $(1.1)$. As a general framework for analyzing such
string propagation, the physicists proposed to consider it as a closed
string theory on the quotient orbifold. The remarkable prediction
$(1.4)$ on the Euler characteristics is a natural by-product of this
setup. Similar ideas have also appeared in combinatorial group
theory. The classical Bass-Serre theory recovers group actions on
graphs from the so-called ``graphs of groups''. For group actions on
simplicial complexes, the corresponding notion is ``complex of
groups''. Haefliger associated to each complex of groups a more
geometric but equivalent object, an orbihedron (cf. \cite{Ha2}). We
will show that an orbihedron is naturally an orbispace in the sense of
this paper, and a homomorphism between two complexes of groups
corresponds to a morphism between the corresponding orbihedra
(regarded as orbispaces in the sense of this paper). The formulation
of the category of orbispaces would provide a general setting for
adopting similar ideas to other problems involving group actions. For
example, it seems to be an interesting problem to examine the
Baum-Connes theory of discrete groups (cf. \cite{BC1,BC2,BC3}) 
from this point of view.

Complexes of groups not only provide a natural class of examples of
orbispaces in the sense of this paper, but may also play a
fundamental role in the theory of orbispaces. In order to explain this
aspect, let us preview some notations in this paper. The orbispace
structure of an orbispace $X$ is given by a collection $\U$ of
connected open subsets satisfying a set of axioms (cf. Definition
2.1.2). Each open subset in $\U$ is called a ``basic open set'' of
$X$. Given any cover $\{U_i\}$ of $X$ by basic open sets, there is an
associated complex of groups, denoted by $G(\{U_i\})$. (Here we assume
$X$ is \'{e}tale, in general we need to slightly modify the definition
of complex of groups given in \cite{Ha2}.) The cover $\{U_i\}$
determines a small category $C(\{U_i\})$ as follows. The objects of
$C(\{U_i\})$ are the set of connected components of intersections of
non-ordered $n$-tuples $(U_{i_1},U_{i_2},\cdots U_{i_n})$ of distinct
elements of $\{U_i\}$, $n=1,2,\cdots$, such that $U_{i_1}\cap
U_{i_2}\cap\cdots\cap U_{i_n}\neq\emptyset$. The morphism between two
objects is the natural inclusion. It is easily seen that $C(\{U_i\})$
is indeed a small category without loop. The barycentric subdivision
of the simplicial complex on which $G(\{U_i\})$ is defined is the
simplicial complex canonically associated to the small category
$C(\{U_i\})$. Such a complex of groups seems to serve as the
``structure group'' of the orbispace $X$. For example, the cohomology
of the complex of groups $G(\{U_i\})$ (cf. \cite{Ha3}) may serve as the
coefficient ring in the cohomology theory of $X$ (cf. Example 2.1.3 c).

\vspace{2mm}

Now we summarize the main results in this paper, followed by a brief
description of the organization of this paper. The reader is referred
to the corresponding sections for the more precise statement and
detailed proof of each result.  

\vspace{2mm}

\noindent{\bf Theorem A:}{\em

\begin{enumerate}
\item For each integer $k\geq 1$, there is a functor $\pi_k$ of $k$-th 
homotopy group from the category of orbispaces with a base-point structure to 
the category of groups (abelian groups when $k\geq 2$), which is homotopy 
invariant (in a sense to be defined in this paper). When the orbispace
is \'{e}tale, an element of the $k$-th homotopy group may be
represented by a morphism from $S^k$ into the orbispace. 
Moreover, for each orbispace $X$, there is a natural homomorphism 
$\Pi_X: \pi_k(X,\ast)\rightarrow \pi_k(X_{top},\ast)$, where $X_{top}$ stands 
for the underlying topological space of the orbispace $X$.
\item Any based pseudo-embedding of orbispaces 
$\tilde{i}:(Y,\underline{q})\rightarrow (X,\underline{p})$ is associated with
a set of relative homotopy groups $\pi_k(X,Y,\tilde{i})$, $k\geq 1$, 
and an exact homotopy sequence:
$$
\cdots\rightarrow\pi_{k+1}(X,Y,\tilde{i})\stackrel{\partial}{\rightarrow}
\pi_k(Y,\underline{q})
\stackrel {i_{\#}}{\rightarrow}\pi_k(X,\underline{p})
\stackrel{j_\#}{\rightarrow}\pi_k(X,Y,\tilde{i})\rightarrow\cdots.
\leqno (1.7)
$$
\item  A notion of orbispace covering is formulated, generalizing the
usual notion in the topological category, and all the basic results in
the topological covering theory are extended to the orbispace
category.  Specializing at the case of orbifolds, our notion of
orbispace covering is equivalent to Thurston's version of orbifold
covering. In particular, the fundamental group of an orbifold defined
here coincides with Thurston's orbifold fundamental group
(cf. \cite{Th}). 
\item A notion of orbispace fibration is formulated,
generalizing the usual notion of fibration, and an analogous Serre's
long exact sequence of homotopy groups is derived.
\item The Seifert-Van Kampen Theorem is extended to the category of
orbispaces.
\item When the orbispace $X$ is canonically defined from a
$G$-space $(Y,G)$, there is a natural isomorphism
$$
\theta_k:\pi_k(X,\ast)\rightarrow \pi_k(Y_G,\ast) \leqno (1.8)
$$
for each $k\geq 1$, where $Y_G=EG\times_G Y$ is the Borel space of
$(Y,G)$.
\item The definition of homotopy groups of orbifolds given in this paper
is equivalent to Haefliger's definition in \cite{Ha1}.
\item For any complex of groups, the fundamental group of its associated
orbihedron (viewed as an orbispace in the sense of this paper) is 
isomorphic to the fundamental group of the complex of groups as 
defined by Haefliger in \cite{Ha2}.
\end{enumerate}
}

\vspace{2mm}

\noindent{\bf Theorem B:}\hspace{2mm}{\em
The free loop space of an orbifold $X$, defined as the space of smooth
morphisms from $S^1$ into $X$, is naturally a pre-Hilbert orbifold. It
admits a canonical $S^1$ action defined by rotating the domain of each
loop, and the space $\widetilde{X}$ as defined in $(1.5)$ is naturally
embedded in the free loop space as the fixed-point set of the
canonical $S^1$ action.
}

\vspace{2mm}

The paper is organized as follows. 

The definition of orbispaces is
given in section 2.1, followed by a discussion on three primary
classes of orbispaces, i.e., the orbispaces canonically defined from
$G$-spaces, orbifolds, and complexes of groups and the associated
orbihedra. A collection of remarks is given after these examples to
further explain various aspects of the definition  
of orbispaces. Section 2.1 is concluded with a pathological example.

The definition of morphisms between orbispaces is introduced in
section 2.2. This section is concluded with a theorem stating that
orbispaces together with the so-defined morphisms form a category.

Section 2.3 introduces the notion of base-point structure of an
orbispace. The based version of the category of orbispaces is
formulated, which will be used in the construction of homotopy theory
of orbispaces. Section 2.3 also introduces the notions of
sub-orbispaces, orbispace embedding, pseudo-embedding, Cartesian
product, etc. The notion of pseudo-embedding will be needed in the
formulation of relative homotopy groups in section 3.3. 

The based loop space of an orbispace with a base-point structure is
introduced in section 3.1. A natural ``compact-open'' topology was
given to the based loop space, and a canonical neighborhood of a based
loop was described. This description is technically useful in many
proofs. A homotopy associative multiplication and a homotopy inverse
were defined using composition and inversion of based loops and the
based loop space was shown to be an H-group under these
operations. Section 3.1 is concluded by observing that the based loop
space defines a functor $\Omega$ from the category of based orbispaces
to the category of H-groups.

Section 3.2 starts with the definition of homotopy groups of a based
orbispace. (The $k$-homotopy group of a based orbispace is defined to
be the $(k-1)$-th homotopy group of the based loop space.) Then a 
homotopy equivalence was introduced amongst based morphisms as being
path-connected in the space of based morphisms equipped with the
``compact-open'' topology. The homotopy groups were shown to be 
homotopy invariant. Path-connectedness of orbispaces was defined
and homotopy groups with path-connected base points were shown to be
canonically isomorphic. The section is concluded with a theorem
stating that the fundamental group of a complex of groups $G(X)$ defined in
\cite{Ha2} is isomorphic to $\pi_1(X)$ where $X$ is the orbihedron
associated to $G(X)$, viewed canonically as an orbispace in the sense
of this paper.

The relative homotopy theory was developed in section 3.3. The notion of
based relative loop space was introduced, associated to any given
pseudo-embedding between two based orbispaces. The $k$-th relative
homotopy group was defined to be the $(k-1)$-th homotopy group of the
based relative loop space. The bulk of section 3.3 was devoted to the
proof of the corresponding homotopy exact sequence $(1.7)$ relating
the relative homotopy groups associated to a pseudo-embedding with the
homotopy groups of each orbispace.

In section 3.4 we show that the homotopy groups of an orbispace which
is canonically defined from a $G$-space are isomorphic to the
(ordinary) homotopy groups of the Borel space of the $G$-space, and
the isomorphism is natural with respect to equivariant maps between
the $G$-spaces.

In section 3.5 we introduce the free loop space of an orbispace, in
particular the free loop space of an orbifold. We identify the free
loop space of an orbispace defined from a $G$-space with the so-called
twisted loop spaces in physics literature. There is a canonical $S^1$ action
on the free loop space defined by rotating the domain of each loop. In
the case of orbifolds, it is shown that the free loop space has a
natural pre-Hilbert orbifold structure, and the space $\widetilde{X}$
defined by $(1.5)$ is embedded in the free loop space as the fixed
point set of the canonical $S^1$ action. The section is ended with
some remarks on further exploiting the pre-Hilbert orbifold structure
on the free loop space of an orbifold. 

The orbispace covering theory was developed in section 4.1. After
introducing the notion of orbispace covering morphism, its based-path
lifting property was established. As an immediate consequence, we
showed that a covering morphism induces isomorphisms on $\pi_k$ with
$k\geq 2$ and a monomorphism on $\pi_1$. We also showed that a based
morphism can be lifted uniquely to a covering space if and only if the
image of $\pi_1$ under the said based morphism is contained in the
image of the $\pi_1$ of the covering space (under the so-called locally
strongly path-connected condition). We then introduce the notion of
universal orbispace covering, which is any orbispace covering space
with trivial $\pi_1$. Under an additional condition of semilocally
1-connectedness, we established the existence of orbispace covering
associated to any given subgroup of the $\pi_1$ of a connected
orbispace, in particular, the universal covering of the orbispace. 
Finally, the notion of deck transformations was given, and a short
exact sequence was established, relating the group of deck
transformations with the normalizer of the image of the $\pi_1$ of the
covering space. As an application, we showed that the definition of
$\pi_1$ of an orbifold in this paper is equivalent to Thurston's
definition in \cite{Th}. The section was ended with a criterion for 
an orbispace to have a universal covering with trivial orbispace
structure. This criterion recovers Haefliger's criterion for the
developability of a complex of groups in \cite{Ha2}.

Section 4.2 concerns the theory of orbispace fibrations. An orbispace
fibration is a morphism of which each local equivariant map is a 
fibration (in the ordinary sense) and the homomorphism between the
corresponding topological groups is a surjective fibration. For each
given orbispace fibration, we constructed the ``fiber'' of the
fibration over a given base-point structure of the base space. The ``fiber'' 
will not be a sub-orbispace of the total space in general, but an
orbispace with a canonical pseudo-embedding into the total
space. Associated to each orbispace fibration, there is an analogous
Serre's exact sequence of homotopy groups, which is derived from the
long exact sequence $(1.7)$ associated to the pseudo-embedding of the
``fiber'' into the total space. As an application, we show that the
homotopy groups of orbifolds defined in this paper are isomorphic to
those defined by Haefliger in \cite{Ha1}. The section ends with some
examples of orbispace fibrations such as orbispace fiber bundles, 
Seifert fibrations, and a canonical fibration from a normal orbispace 
onto its canonical reduction. 

The classical Seifert-Van Kampen Theorem was generalized to the
orbispace category in section 4.3.

\vspace{3mm} 

This is a substantial revision of an early version of the paper where
Haefliger's work, which is very relevant to the discussion here,  
was not properly mentioned. This revision owes its existence to Eugene
Lerman who had pointed out Haefliger's paper \cite{Ha1} to the author 
shortly after the early version was posted in the e-print
archive. After further searching on Haefliger's work, the author found
\cite{Ha2, Ha3}. Thanks also go to Charles Boyer for explaining
the results of Haefliger in \cite{Ha1} which is written in French. Finally,
the author is indebted to Dennis Sullivan for an enlightening
conversation after a preliminary version of this paper was completed.

\section{Category of Orbispaces}

\hspace{5mm}This section is devoted to the foundation of category of 
orbispaces. An 
orbispace is locally modeled on spaces with group actions. Recall that given 
a topological group $G$, a space $Y$ is called a $G$-space if there is a
continuous action $G\times Y\rightarrow Y$, written $(g,y)\mapsto g\cdot y$,
satisfying $g\cdot(g^\prime\cdot y)=(gg^\prime)\cdot y$ and $1_G\cdot y=y$.
In this paper, we will not confine ourselves to actions by a fixed group. 
Hence in order to simplify the exposition, we shall call any space with a 
group action a $G$-space by abusing the notation.

\subsection{Orbispaces}

\hspace{5mm}
Throughout this paper, we shall generally assume that the topological spaces 
under consideration are {\it locally connected}. Recall that a space $X$ is 
locally connected if for any point $p\in X$, and any neighborhood $V$ of $p$ 
in $X$, there is a connected neighborhood $U$ of $p$ such that $U\subset V$. 
It is easily seen that any open subspace of a locally connected space is 
locally connected. The practical reason for which we assume 
locally-connectedness is that a locally connected space can be decomposed into 
a disjoint union of open connected subspaces, each of which is called a 
{\it connected component} of the space. 

Locally-connectedness is preserved by the orbit space of a
$G$-space. More precisely, let $Y$ be a locally connected space with a
continuous action of a topological group $G$, then the space of orbits
$X:=Y/G$, given with the quotient topology, is also locally connected.
In fact, for any open subset $U\subset X$, consider the inverse image
$\pi^{-1}(U)\subset Y$, where $\pi:Y\rightarrow X$ is the natural
projection. As an open subset of a locally connected space,
$\pi^{-1}(U)$ is decomposed into a disjoint union of connected
components. The action of $G$ on $Y$ gives rise to an action on the
set of these connected components. The image of each component under
$\pi$ is the same within an orbit, hence $U$ is decomposed into a
disjoint union of open connected subsets. This particularly implies
that the orbit space $X$ is locally connected.

Let $U$ be a connected, locally connected topological space. By a 
{\it $G$-structure} on $U$ we mean a triple $(\widehat{U},G_U,\pi_U)$ where 
$(\widehat{U},G_U)$ is a connected, locally connected 
$G$-space and $\pi_U:\widehat{U}\rightarrow U$ is a continuous 
map inducing a homeomorphism between the orbit space $\widehat{U}/G_U$ 
and $U$. An {\it isomorphism} between two $G$-structures on $U$, 
$(\widehat{U}_i,G_{i,U},\pi_{i,U})$ for $i=1,2$, is a pair $(\phi,\lambda)$ 
where $\lambda:G_{1,U}\rightarrow G_{2,U}$ is an isomorphism and $\phi:
\widehat{U}_1\rightarrow\widehat{U}_2$ is a $\lambda$-equivariant homeomorphism
such that $\pi_{2,U}\circ \phi=\pi_{1,U}$. By the {\it domain
(resp. range)} of an isomorphism $(\phi,\lambda)$ of $G$-structures we mean
the domain (resp. range) of $\phi$ and $\lambda$.
Note that each $g\in G_U$ induces
an automorphism $(\phi_g,\lambda_g)$ on $(\widehat{U},G_U,\pi_U)$, defined 
by setting $\phi_g(x)= g\cdot x,\forall x\in\widehat{U}$ and 
$\lambda_g(h)=ghg^{-1},\forall h\in G_U$.
However, it might not be true that every automorphism arises in this way. 
We shall only take into consideration the automorphisms $(\phi_g,\lambda_g)$, 
$g\in G_U$. More precisely, we {\it define} the automorphism group of the 
$G$-structure $(\widehat{U},G_U,\pi_U)$ to be $G_U$ via the induced 
isomorphisms $(\phi_g,\lambda_g)$ on it. We would like to  
point out that since the action of $G_U$ on $\widehat{U}$ is not required to
be effective, two {\it different} automorphisms of the $G$-structure could 
have the {\it same} induced map.

\vspace{1.5mm}

Given a $G$-structure $(\widehat{U},G_U,\pi_U)$ on $U$, we consider the 
inverse image $\pi_U^{-1}(W)$ in $\widehat{U}$, where $W$ is a 
connected open subset of $U$. Denote by $\widehat{W}$ one of the connected
components of $\pi_U^{-1}(W)$, by $G_W$ the subgroup of $G_U$ consisting of 
elements $g\in G_U$ such that $g\cdot\widehat{W}=\widehat{W}$, and let 
$\pi_W=(\pi_U)|_{\widehat{W}}$. 

\vspace{1.5mm}

\noindent{\bf Lemma 2.1.1:}
{\em The triple $(\widehat{W},G_W,\pi_W)$ defines a $G$-structure on $W$. 
Moreover, the action of $G_U$ on $\widehat{U}$ induces a transitive action on 
the set of all such $G$-structures on $W$ for which the following holds: 
Let $(\widehat{W_i}, G_{W,i}, \pi_{W,i})$, $i=1,2$, be two such 
$G$-structures, and $\widehat{W_2}=g\cdot\widehat{W_1}$ for some $g\in G_U$, 
then $G_{W,2}=gG_{W,1}g^{-1}$ in $G_U$. The stabilizer of the action at 
$G$-structure $(\widehat{W},G_W,\pi_W)$ is precisely the subgroup $G_W$ in 
$G_U$.
}

\vspace{2mm}

\noindent{\bf Proof:}
The action of $G_U$ on the set of connected components of $\pi^{-1}_U(W)$ is
transitive because $W$ is connected. This implies that for any connected 
component $\widehat{W}$, the map $\pi_W:=\pi_U|_{\widehat{W}}:\widehat{W}
\rightarrow W$ is surjective. On the 
other hand, if $\pi_W(x)=\pi_W(y)$ for some $x,y\in \widehat{W}$, then there 
is a $g\in G_U$ such that $g\cdot x=y$. Clearly $g$ preserves $\widehat{W}$, 
hence $g$ lies in $G_W$. It follows easily now that $\pi_W:\widehat{W}
\rightarrow W$ induces a homeomorphism between the orbit space 
$\widehat{W}/G_W$ and $W$, hence $(\widehat{W},G_W,\pi_W)$ defines a 
$G$-structure on $W$. 

The rest of the lemma is straightforward, and we leave the details to the 
readers. 

\hfill $\Box$

We will say that $(\widehat{W},G_W,\pi_W)$ is {\it induced} from the 
$G$-structure $(\widehat{U},G_U,\pi_U)$. Note that the subgroup $G_W$ is 
both closed and open in $G_U$. In fact, the space of cosets
$G_U/G_W=\{gG_W|g\in G_U\}$ inherits a discrete topology from $G_U$. 

\vspace{1.5mm}

With these preparations, we are ready to define orbispace as in the following

\vspace{1.5mm}

\noindent{\bf Definition 2.1.2:}
{\em Let $X$ be a locally connected topological space. An {\it orbispace 
structure} on $X$ is a collection $\U$ of open subsets of $X$ satisfying the 
following conditions:
\begin{enumerate}
\item Each element $U$ of $\U$ is connected and $\U$ is a base of $X$.
\item Each element $U$ of $\U$ is assigned with a $G$-structure $(\widehat{U},
G_U,\pi_U)$ satisfying the following conditions:
\begin{itemize}
\item [{a)}] For any $U_\alpha,U_\beta\in\U$ such that $U_\alpha\cap U_\beta
\neq\emptyset$, there is a non-empty set $Tran(U_\alpha,U_\beta)
=\{(\phi,\lambda)\}$, where each $(\phi,\lambda)$ is an isomorphism
from a $G$-structure of a connected component of $U_\alpha\cap U_\beta$ induced
from the $G$-structure of $U_\alpha$ to a $G$-structure induced from the
$G$-structure of $U_\beta$. Each $(\phi,\lambda)$ will be called a transition 
map.
\item [{b)}] For any transition maps $(\phi_i,\lambda_i)\in 
Tran(U_\alpha,U_\beta)$, $i=1,2$,
which are isomorphisms between induced $G$-structures of the same 
connected component of $U_\alpha\cap U_\beta$, there are 
$g_\alpha\in G_{U_\alpha}$, $g_\beta\in G_{U_\beta}$ such that 
$(\phi_2,\lambda_2)=g_\beta^{-1}\circ (\phi_1,\lambda_1)\circ g_\alpha$.
Here $g_\alpha$, $g_\beta$ are regarded as elements of the automorphism group
of the $G$-structures of $U_\alpha$ and $U_\beta$ respectively. Moreover, for 
any $(\phi,\lambda)\in Tran(U_\alpha,U_\beta)$, $g_\alpha\in G_{U_\alpha}$, 
$g_\beta\in G_{U_\beta}$, the composition 
$g_\beta^{-1}\circ (\phi,\lambda)\circ g_\alpha$ is
in $Tran(U_\alpha,U_\beta)$, and 
$g_\beta^{-1}\circ (\phi,\lambda)\circ g_\alpha$, $(\phi,\lambda)$ are 
regarded as the same element in $Tran(U_\alpha,U_\beta)$ if and only if 
$\lambda(g_\alpha)=g_\beta$ (in particular, $g_\alpha$ is in the domain of
$\lambda$).
\item [{c)}] The identity map $Id:(\widehat{U},G_U,\pi_U)\rightarrow 
(\widehat{U},G_U,\pi_U)$ is contained in $Tran(U,U)$ for any $U\in \U$. Hence
by virture of $b)$, $Tran(U,U)$ is naturally identified with $G_U$ as the
automorphism group of the $G$-structure $(\widehat{U},G_U,\pi_U)$.
\item [{d)}] Let $U_\alpha,U_\beta,U_\gamma\in \U$ be any three elements
such that $U_\alpha\cap U_\beta\cap U_\beta\neq\emptyset$. 
For any $(\phi_{\gamma\beta},\lambda_{\gamma\beta})\in
Tran(U_\beta,U_\gamma)$ and $(\phi_{\beta\alpha},\lambda_{\beta\alpha})\in
Tran(U_\alpha,U_\beta)$, if the domain of $(\phi_{\beta\alpha},
\lambda_{\beta\alpha})$ intersects the range of 
$(\phi_{\gamma\beta},\lambda_{\gamma\beta})$, then there is a transition map in
$Tran(U_\alpha,U_\gamma)$, which will be called the {\it composition} of 
$(\phi_{\beta\alpha},\lambda_{\beta\alpha})$ with
$(\phi_{\gamma\beta},\lambda_{\gamma\beta})$ and denoted by 
$(\phi_{\gamma\beta},\lambda_{\gamma\beta})\circ 
(\phi_{\beta\alpha},\lambda_{\beta\alpha})$, such that 
the restriction of $(\phi_{\gamma\beta},\lambda_{\gamma\beta})\circ 
(\phi_{\beta\alpha},\lambda_{\beta\alpha})$ to $\pi_{U_\alpha}^{-1}
(U_\alpha\cap U_\beta\cap U_\gamma)$ is equal to the composition of 
$(\phi_{\beta\alpha},\lambda_{\beta\alpha})$ with 
$(\phi_{\gamma\beta},\lambda_{\gamma\beta})$ as maps. Furthermore, the 
operation of composition of transition maps is associative. When 
$\alpha=\beta$ (resp. $\beta=\gamma$), the composition 
$(\phi_{\gamma\beta},\lambda_{\gamma\beta})\circ 
(\phi_{\beta\alpha},\lambda_{\beta\alpha})$ coincides with the usual 
composition as understood in $b)$ where $(\phi_{\beta\alpha},
\lambda_{\beta\alpha})$ (resp. $(\phi_{\gamma\beta},
\lambda_{\gamma\beta})$) is understood as an element of $G_{U_\alpha}$ 
(resp. $G_{U_\gamma}$) by $c)$. 
\end{itemize}
\end{enumerate}

The topological space $X$ equipped with the orbispace structure $\U$
is called an {\it orbispace}, and will be denoted by $(X,\U)$ in
general.  
} 

\hfill $\Box$ 

Each element $U$ of $\U$ is called a {\it basic open set} of $X$.  For
any point $p\in X$, which is contained in a basic open set $U$, pick a
$x\in \widehat{U}$ in the inverse image $\pi_U^{-1}(p)$. We define the
{\it isotropy group} of $p$ to be the stabilizer $G_x$ of $x$ in $G_U$
(i.e.  $G_x=\{g\in G_U|g\cdot x=x\}$), and denote it by $G_p$. Clearly
different choices of $x$ result in the same conjugacy class in $G_U$,
and different choices of $U$ give rise to isomorphic groups because of
the existence of transition maps. An orbispace structure $\U$ is
called {\it trivial} if $G_U$ is trivial for each $U\in\U$. (Every
locally connected topological space is canonically an orbispace with a
trivial orbispace structure.) An orbispace $(X,\U)$ is called {\it
\'{e}tale} if $G_U$ is discrete for each $U\in\U$. As a notational
convention, we very often only write $X$ for an orbispace $(X,\U)$,
and write $X_{top}$ for the underlying topological space $X$ for
simplicity.

\vspace{3mm}

A collection of remarks will be given to further explain the various
aspects of Definition 2.1.2. But we shall first look at the following three
primary classes of examples of orbispaces.

\vspace{1.5mm}

\noindent{\bf Example 2.1.3 a:}\hspace{2mm}
For any locally connected $G$-space $(Y,G)$, the orbit space $Y/G$ 
canonically inherits an
orbispace structure, of which $\U$ is taken to be the set of all connected open
subsets of the orbit space $Y/G$. The $G$-structure assigned to each element 
of $\U$ is a fixed choice of the $G$-structures induced from the God-given 
$G$-structure $(Y,G)$ on the orbit space $Y/G$. Each set of transition maps is
obtained by restricting $G$ to the corresponding induced $G$-structures. The
verification of Definition 2.1.2 for this case is straightforward. An orbispace
will be called {\it global} if it arises as the orbit space of a $G$-space 
equipped with the canonical orbispace structure as discussed in this example.

\hfill $\Box$

\noindent{\bf Example 2.1.3 b:}\hspace{2mm} 
Any orbifold is naturally
an orbispace. The set $\U$ is taken to be the set of all convex
geodesic neighborhoods (assuming a Riemannian metric is given).  A
$G$-structure is just a uniformizing system, and transition maps are
to be derived from {injections} between uniformizing systems. Recall 
that for any inclusion $U\subset W$ of geodesic
neighborhoods, an injection from the uniformizing system
$(\widehat{U},G_U,\pi_U)$ of $U$ into the uniformizing system
$(\widehat{W},G_W,\pi_W)$ of $W$ is just an isomorphism onto one of
the uniformizing systems of $U$ induced from
$(\widehat{W},G_W,\pi_W)$. The set of injections for each inclusion
$U\subset W$ is a $G_W$-homogeneous space over a point.

For any two elements $U_\alpha$, $U_\beta$ in $\U$ such that 
$U_\alpha\cap U_\beta\neq\emptyset$, the set of transition maps 
$Tran(U_\alpha,U_\beta)$ is defined as follows: Set $U_{\alpha\beta}
=U_\alpha\cap U_\beta$, we define $Tran(U_\alpha,U_\beta)$
to be the set of all $(b,i_b)\circ (a,i_a)^{-1}$ where $(a,i_a)$ 
(resp. $(b,i_b)$) is an injection from $(\widehat{U}_{\alpha\beta},
G_{U_{\alpha\beta}},\pi_{U_{\alpha\beta}})$
into $(\widehat{U}_\alpha,G_{U_\alpha},\pi_{U_\alpha})$ (resp.
$(\widehat{U}_\beta,G_{U_\beta},\pi_{U_\beta})$). One is ready to verify
that two $(b,i_b)\circ (a,i_a)^{-1}$ and $(b,i_b)^\prime\circ 
((a,i_a)^\prime)^{-1}$ are the same if and only if $(a,i_a)^\prime=(a,i_a)\circ
g$ and $(b,i_b)^\prime=(b,i_b)\circ g^{-1}$ for some 
$g\in G_{U_{\alpha\beta}}$. The axioms {\em 2-a),b),c)} in Definition 2.1.2 
are obvious. As for {\em 2-d)} of Definition 2.1.2, suppose $U_\alpha$, 
$U_\beta$, $U_\gamma$ are geodesic neighborhoods satisfying
$U_\alpha\cap U_\beta\cap U_\gamma\neq\emptyset$. Set 
$U_{\alpha\beta\gamma}=U_\alpha\cap U_\beta\cap U_\gamma$. Let $(b,i_b)\circ
(a,i_a)^{-1}$ be a transition map in $Tran(U_\alpha,U_\beta)$ and 
$(d,i_d)\circ (c,i_c)^{-1}$ in $Tran(U_\beta,U_\gamma)$. We need to define 
the composition $((d,i_d)\circ (c,i_c)^{-1})\circ ((b,i_b)\circ (a,i_a)^{-1})$.
It is done as follows: Pick an injection $(\xi,i_\xi)$ from 
$(\widehat{U}_{\alpha\beta\gamma},G_{U_{\alpha\beta\gamma}},
\pi_{U_{\alpha\beta\gamma}})$ into $(\widehat{U}_{\alpha\beta},
G_{U_{\alpha\beta}},\pi_{U_{\alpha\beta}})$, then there is a unique injection
$(\eta,i_\eta)$ such that $(b,i_b)\circ (\xi,i_\xi)=(c,i_c)\circ 
(\eta,i_\eta)$. It is easily seen that the restriction of $(d,i_d)\circ 
(c,i_c)^{-1}\circ (b,i_b)\circ (a,i_a)^{-1}$ to $\pi_{U_\alpha}^{-1}
(U_{\alpha\beta\gamma})$ is equal to 
$((d,i_d)\circ (\eta,i_\eta))\circ ((a,i_a)\circ (\xi,i_\xi))^{-1}$. 
We pick an injection $(\theta,i_\theta)$ from 
$(\widehat{U}_{\alpha\beta\gamma},G_{U_{\alpha\beta\gamma}},
\pi_{U_{\alpha\beta\gamma}})$ into $(\widehat{U}_{\alpha\gamma},
G_{U_{\alpha\gamma}},\pi_{U_{\alpha\gamma}})$, then there are unique 
injections $(e,i_e)$ and $(f,i_f)$ such that $(e,i_e)\circ (\theta,i_\theta)=
(a,i_a)\circ (\xi,i_\xi)$ and $(f,i_f)\circ (\theta,i_\theta)=
(d,i_d)\circ (\eta,i_\eta)$. We simply define the composition by
$$
((d,i_d)\circ (c,i_c)^{-1})\circ ((b,i_b)\circ (a,i_a)^{-1}):=(f,i_f)\circ
(e,i_e)^{-1} \leqno (2.1.1)
$$ 
in $Tran(U_\alpha,U_\gamma)$. One can verify that it is 
well-defined and associative.

\hfill $\Box$

\noindent{\bf Example 2.1.3 c:}\hspace{2mm}
In this example we shall discuss the notions of complex of groups and
orbihedron as given by Haefliger in \cite{Ha2}. 

Let $X$ be a simplicial cell complex. We set $V(X)$ for the set of
barycenters of cells of $X$, and $E(X)$ for the set of edges of the
barycentric subdivision of $X$. Each edge $a\in E(X)$ has a natural
orientation: if the initial point $i(a)$ is the barycenter of a
cell $\sigma$ and the terminal point $t(a)$ is the barycenter of a
cell $\tau$, then $\dim\tau < \dim\sigma$. Two edges $a,b\in E(X)$ are
said to be {composable} if $i(a)=t(b)$ and the composition $c=ab$
is the edge $c$ with $i(c)=i(b)$ and $t(c)=t(a)$ such that $a,b$ and
$c$ form the boundary of a 2-simplex of the barycentric subdivision.

A complex of groups $G(X)=(X,G_\sigma,\psi_a,g_{a,b})$ on $X$ is given
by \\
1). a group $G_\sigma$ (with discrete topology) for each cell 
$\sigma\in V(X)$;\\
2). an injective homomorphism $\psi_a: G_{i(a)}\rightarrow G_{t(a)}$
for each edge $a\in E(X)$;\\
3). for any composable edges $a,b$ ($i(a)=t(b)$), an element
$g_{a,b}\in G_{t(a)}$ is given such that 
$$
g_{a,b}\psi_{ab}(h)g_{a,b}^{-1}=\psi_a(\psi_b(h)), \hspace{2mm}
\forall h\in G_{i(b)}, \leqno (2.1.2a)
$$
and the set of elements $\{g_{a,b}\}$ satisfies the following cocycle
condition for any triple $a,b,c$ of composable edges
$$
\psi_a(g_{b,c})g_{a,bc}=g_{a,b}g_{ab,c}. \leqno (2.1.2b)
$$

An orbihedron structure on a simplicial cell complex $X$ is given by
the following data:\\
1). For each cell $\sigma$ of $X$, a simplicial cell complex
$Lk\tilde{\sigma}$ is given, as well as a simplicial action without
inversion of a group $G_\sigma$ on $Lk\tilde{\sigma}$ and a simplicial
projection $Lk(p_\sigma):Lk\tilde{\sigma}\rightarrow Lk\sigma$
inducing an isomorphism of $Lk\tilde{\sigma}/G_\sigma$ with
$Lk\sigma$. On the join of the closure $\bar{\sigma}$ of $\sigma$ with
$Lk\tilde{\sigma}$, there is a simplicial action without inversion of
$G_\sigma$ which is given by the join of the identity action on
$\bar{\sigma}$ with the given action on $Lk\tilde{\sigma}$. Denote by
$\tilde{\sigma}$ the simplex $\sigma$ in the join, by
$St\tilde{\sigma}$ the star of $\sigma$ in the join. Then the action
of $G_\sigma$ restricts to $St\tilde{\sigma}$ and there is a
simplicial projection $p_\sigma:St\tilde{\sigma}\rightarrow St\sigma$
inducing an isomorphism of $St\tilde{\sigma}/G_\sigma$ with $St\sigma$.\\
2). For any edge $a\in E(X)$ with $i(a)=\tau$ and $t(a)=\sigma$, there
is an injective homomorphism $\psi_a:G_\tau\rightarrow G_\sigma$ and a
$\psi_a$-equivariant simplicial map 
$f_a:p_\tau^{-1}(Sta)\rightarrow p_\sigma^{-1}(Sta)$ which is a homeomorphism
onto an open subset projecting onto $Sta\subset St\sigma$. Moreover,
for any point $x\in St\tilde{\tau}$, the restriction of $\psi_a$ to the
stabilizer of $x$ in $G_\tau$ is an isomorphism onto the stabilizer of
$f_a(x)$ in $G_\sigma$. \\
3). For any composable edges $a,b\in E(X)$, an element $g_{a,b}\in
G_{t(a)}$ is given such that $f_a\circ f_b=g_{a,b}\circ f_{ab}$ and
$g_{a,b}\psi_{ab}(h)g_{a,b}^{-1}=\psi_a(\psi_b(h)), \forall h\in
G_{i(b)}$,  and that $\psi_a(g_{b,c})g_{a,bc}=g_{a,b}g_{ab,c}$ for any
triple $a,b,c$ of composable edges in $E(X)$. 

Given the above data of an orbihedron structure on $X$, a canonical
\'{e}tale topological groupoid $\Gamma=\Gamma(X)$ can be constructed 
as follows. The space of units is the disjoint union of
the $St\tilde{\sigma}$'s. The restriction $\Gamma_\sigma$ of $\Gamma$
is the topological groupoid $G_\sigma\times St\tilde{\sigma}$
associated to the action of $G_\sigma$ on $St\tilde{\sigma}$. For any
edge $a\in E(X)$ with $i(a)=\sigma$ and $t(a)=\tau$, the space
$\Gamma_a$ of elements $\gamma\in \Gamma$ with $i(\gamma)\in
p_\sigma^{-1}(Sta)$ and $t(\gamma)\in St\tilde{\tau}$ is the product
$G_\tau \times p_\sigma^{-1}(Sta)$. Such an element may be represented
by $\gamma=(g,f_a,x)$, where $x\in p_\sigma^{-1}(Sta)$ and $g\in
G_\tau$, with $i(\gamma)=x$ and $t(\gamma)=gf_a(x)$. The space of
elements $\gamma\in \Gamma$ with $i(\gamma)\in St\tilde{\sigma}$ and
$t(\gamma) \in St\tilde{\tau}$ is the union of the $\Gamma_a$'s where
$a\in E(X)$ satisfying $i(a)=\sigma$ and $t(a)=\tau$. The composition
$(h,y)\circ (g,f_a,x)$, where $y=gf_a(x)$, is $(hg,f_a,x)$, and the
composition $(g,f_a,x)\circ (k,z)$, where $k\in G_{i(a)}$ and
$z=k^{-1}x$, is $(g\psi_a(k),f_a,z)$. For a pair of composable edges
$a,b\in E(X)$, the composition $(g,f_a,x)\circ (h,f_b,y)$, where
$x=hf_b(y)$, is defined to be $(gg_{a,b}\psi_a(h),f_{ab},y)$. The
associativity of such compositions is ensured by the conditions in
3) above.

Observe that in an orbihedron structure on $X$, a complex of
groups on $X$ is naturally contained, which is called the complex of
groups associated to the orbihedron structure. On the other hand, 
given a complex of groups $G(X)$ on $X$, Haefliger canonically
constructed an orbihedron structure on $X$ whose associated complex of
groups is $G(X)$. The reader is referred to \cite{Ha2} for the details
of the construction.

Given an orbihedron structure on a simplicial cell complex $X$, we can
canonically define an orbispace structure on $X$ as follows. We let
$\U=\{St\sigma,\sigma\in V(X)\}$. The $G$-structure of $St\sigma$ is
defined to be $(St\tilde{\sigma},G_\sigma,p_\sigma)$. Observe that for
any two different cells $\sigma$ and $\tau$, $St\sigma\cap
St\tau\neq\emptyset$ if and only if $\dim\sigma\neq \dim\tau$ and
there is an edge $a\in E(X)$ such that $i(a)=\sigma$ and $t(a)=\tau$
(assuming $\dim\sigma > \dim\tau$), and in the latter case,
$St\sigma\cap St\tau$ is the disjoint union of all the $Sta$'s where
$a\in E(X)$ satisfy $i(a)=\sigma$ and $t(a)=\tau$. For each such an
edge $a$, we define the set of transition maps to be $\{(g\circ
f_a,g\psi_ag^{-1}), g\in G_{t(a)}\}$. The axioms 2-a),b),c) in
Definition 2.1.2 are naturally satisfied, and the axiom 2-d) is
ensured by the conditions in 3) of the definition of orbihedron
structure. Thus $\U$ determines an orbispace structure on $X$.

Let $G(X)=(X,G_{\sigma},\psi_a,g_{a,b})$ and
$G(X^\prime)=(X^\prime,G_{\sigma^\prime}^\prime,\psi_{a^\prime}^\prime,
g_{a^\prime,b^\prime}^\prime)$ be complexes of groups on simplicial
cell complexes $X$ and $X^\prime$ respectively. Let $\theta:X\rightarrow
X^\prime$ be a simplicial map. A homomorphism $\Phi$ from $G(X)$ to
$G(X^\prime)$ over $\theta$ is given by\\
1). a homomorphism $\phi_\sigma:G_\sigma\rightarrow G_{\theta(\sigma)}$ for
each cell $\sigma$ of $X$;\\
2). an element $g_a^\prime\in G_{t(\theta(a))}^\prime$ for each edge $a\in
E(X)$ such that 
$$
g_a^\prime\psi_{\theta(a)}^\prime(\phi_{i(a)}(h))(g_a^\prime)^{-1}=
\phi_{t(a)}(\psi_a(h)),\forall h\in G_{i(a)}, \leqno (2.1.3a)
$$
and such that for composable edges $a,b$ of $E(X)$
$$
\phi_{t(a)}(g_{a,b})g_{ab}^\prime=g_a^\prime\psi_{\theta(a)}^\prime(g_b^\prime)
g^\prime_{\theta(a),\theta(b)}. \leqno (2.1.3b)
$$

Given a homomorphism $\Phi=(\theta,\phi_\sigma,g_a^\prime): 
G(X)\rightarrow G(X^\prime)$ over $\theta:
X\rightarrow X^\prime$, there is a canonically defined morphism
between $X$ and $X^\prime$, where $X$ and $X^\prime$ are regarded as
orbispaces canonically defined from the orbihedron structures
constructed by Haefliger. (Morphisms between orbispaces are defined in
the next section.) This amounts to give a collection of maps $\{u_\sigma: 
St\tilde{\sigma}\rightarrow St\widetilde{\theta(\sigma)},\sigma\in V(X)\}$ 
where $u_\sigma$ is equivariant with respect to the homomorphism
$\phi_\sigma:G_\sigma\rightarrow G^\prime_{\theta(\sigma)}$, and a
collection of maps $\{\rho_a:\{(g\circ f_a,g\psi_ag^{-1}),g\in
G_{t(a)}\}\rightarrow \{(g^\prime\circ f_{\theta(a)}^\prime,
g^\prime\psi_{\theta(a)}^\prime (g^\prime)^{-1}), g^\prime\in 
G^\prime_{t(\theta(a))}\}, a\in E(X)\}$ satisfying (cf. $(2.2.1a)$, 
$(2.2.1b)$)
$$
(u_\sigma,\phi_\sigma)\circ (g\circ f_a,g\psi_ag^{-1})=
\rho_a((g\circ f_a,g\psi_ag^{-1}))\circ (u_\tau,\phi_\tau),
\forall g\in G_\sigma, \leqno (2.1.4a)
$$
and 
$$
(\phi_\sigma(g), Ad(\phi_\sigma(g)))\circ \rho_a((f_a,\psi_a))=
\rho_a((g\circ f_a, g\psi_ag^{-1})), \leqno (2.1.4b)
$$
where $a\in E(X)$ satisfies $i(a)=\tau$ and $t(a)=\sigma$, 
$(\phi_\sigma(g), Ad(\phi_\sigma(g)))$ is the self-equivariant map of
$(St\widetilde{\theta(\sigma)}, G_{\theta(\sigma)}^\prime)$ which is
given by the action of $\phi_\sigma(g)$ on the first factor and
conjugation by $\phi_\sigma(g)$ on the second factor, and 
$$
\rho_{ab}((g_{a,b}\circ f_{ab},g_{a,b}\psi_{ab}g_{a,b}^{-1}))=
\rho_a((f_a,\psi_a))\circ \rho_b((f_b,\psi_b)) \leqno (2.1.4c)
$$
for any composable edges $a,b$ of $E(X)$. In the construction of
orbihedron structure on $X$ associated to a complex of groups $G(X)=
(X,G_\sigma,\psi_a,g_{a,b})$, the simplicial cell complex
$Lk\tilde{\sigma}$ for each cell $\sigma$ consists of cells labelled
by $(g\psi_a(G_{i(a)}),a)$ where $g\in G_\sigma$ and $a\in E(X)$ with 
$t(a)=\sigma$. The maps $u_\sigma:
St\tilde{\sigma}\rightarrow St\widetilde{\theta(\sigma)}$ generated
by correspondence 
$$
(g\psi_a(G_{i(a)}),a)\mapsto
(\phi_\sigma(g)g_a^\prime\psi_{\theta(\sigma)}^\prime(G^\prime_{i(\theta(a))}),
\theta(a)), \forall a\in E(X) \mbox{ s.t.} \hspace{2mm}
t(a)=\sigma, \leqno (2.1.5a)
$$
which are easily checked to be $\phi_\sigma$ equivariant, and the maps 
$$
\rho_a: (g\circ f_a,g\psi_ag^{-1})\mapsto (\phi_{t(a)}(g)\circ
g_a^\prime\circ f_{\theta(a)}^\prime,
\phi_{t(a)}(g)g_a^\prime\psi_{\theta(a)}^\prime
(\phi_{t(a)}(g)g_a^\prime)^{-1}) \leqno (2.1.5b)
$$
can be checked to satisfy the conditions $(2.1.4a-c)$, hence define
the desired morphism between the orbispaces $X$ and $X^\prime$. 

Our formulation of category of orbispaces suggests a 
slightly more general definition of complex of groups, which we call a
{\it complex of topological groups}, and a generalized form of
homomorphisms between complexes of (topological) groups. We define the
notion of complex of topological groups by further requiring in the
definition of complex of groups that each $G_\sigma$ is a topological
group and the image of each injective homomorphism $\psi_a:G_{i(a)}\rightarrow
G_{t(a)}$ is both open and closed. We define the generalized form of
homomorphisms by allowing the simplicial map $\theta$ not being
homeomorphic over each cell but just onto, and requiring that in the
case $i(\theta(a))=t(\theta(a))$, $\phi_{\theta(a)}^\prime$ be the
identity homomorphism of $G^\prime_{t(\theta(a))}$. 

Let $(X,\U)$ be an orbispace. Given any cover $\{U_i\}$ of $X$ where
each $U_i\in \U$, we can canonically construct a complex of topological groups
$G(\{U_i\})$, unique up to a coboundary deduction
(cf. \cite{Ha2}). The simplicial cell complex $X(\{U_i\})$ over which 
$G(\{U_i\})$ is defined is constructed as follows. The n-cells of $X(\{U_i\})$
are labelled by $(\{U_{i_0},\cdots, U_{i_n}\},j)$, where $U_{i_0},
\cdots, U_{i_n}$ are distinct elements of $\{U_i\}$ such that 
$U_{i_0}\cap\cdots\cap U_{i_n}\neq \emptyset$, $j$ stands
for a connected component of $U_{i_0}\cap\cdots\cap U_{i_n}$. The
faces of $(\{U_{i_0},\cdots, U_{i_n}\},j)$ are obtained by deleting
one or several $U_{i_k}$'s. The topological group associated to
the cell $(\{U_{i_0},\cdots, U_{i_n}\},j)$ is the group in the 
$G$-structure of $j$. Each edge of the barycentric subdivision of
$X(\{U_i\})$ is given by the inclusion $j\hookrightarrow j_{\{i_k\}}$
where $j_{\{i_k\}}$ is the corresponding connected component of the
intersection of $\{U_{i_0},\cdots, U_{i_n}\}$ with one or several 
$U_{i_k}$'s deleted. For such an edge $a$, the injective homomorphism $\psi_a$
is taken to be the group monomorphism in a fixed choice of transition maps
associated to the inclusion $j\hookrightarrow j_{\{i_k\}}$. Then the
axioms 2-b),d) in Definition 2.1.2 ensure the existence of $g_{a,b}$'s
satisfying $(2.1.2a-b)$ (the cocycle condition $(2.1.2b)$ follows from
the associativity of compositions of transition maps). Different
choices of transition maps  in defining $\psi_a$'s result in a complex of
topological groups differed by a coboundary. 

Suppose $\{V_k\}$ is a refinement of $\{U_i\}$ where each
$V_k\in\U$. Then there is a simplicial map
$\theta:X(\{V_k\})\rightarrow X(\{U_i\})$ canonically induced from the
refinement relation between $\{V_k\}$ and $\{U_i\}$, which is onto over each
cell. Fixing a choice of transition maps associated to each inclusion
$\sigma\hookrightarrow \theta(\sigma)$ inherited from the refinement
relation between $\{V_k\}$ and $\{U_i\}$, we take the group
homomorphism in the chosen transition map for the homomorphism
$\phi_\sigma$. Then again by the axioms 2-b),d) in Definition 2.1.2,
there exist $\{g_a^\prime\}$ satisfying $(2.1.3a-b)$ so that 
$(\theta,\phi_\sigma,g_a^\prime)$ defines a
generalized homomorphism from $G(\{V_k\})$ to $G(\{U_i\})$. We
leave the details to the readers. 

\hfill $\Box$ 

The following remarks are intended to further explain certain aspects of 
Definition 2.1.2.

\vspace{1.5mm}

\noindent{\bf Remark 2.1.4 a:}\hspace{2mm}
When a basic open set $U_\alpha$ is included in another one $U_\beta$,
a transition map in $Tran(U_\alpha,U_\beta)$ is just an isomorphism from 
the $G$-structure of $U_\alpha$ onto one of the $G$-structures of $U_\alpha$ 
induced from the $G$-structure of $U_\beta$, something analogous to 
injections in the case of orbifolds. As a matter of fact, one can formulate 
the definition of orbispaces only using such transition maps. We choose to 
introduce the more general notion of transition maps because it turns out to
be much more convenient in defining and studying morphisms
between orbispaces. We will continue to call a transition map in 
$Tran(U_\alpha,U_\beta)$ an {\it injection} when $U_\alpha\subset U_\beta$ 
holds.

\hfill $\Box$

\noindent{\bf Remark 2.1.4 b:}\hspace{2mm}
Suppose $U_\alpha$, $U_\beta$, and $U_\gamma$ are basic open sets such that
$U_\alpha\cap U_\beta\cap U_\gamma\neq \emptyset$. In the equation 
$$
(\phi_{\gamma\alpha},\lambda_{\gamma\alpha})=
(\phi_{\gamma\beta},\lambda_{\gamma\beta})\circ 
(\phi_{\beta\alpha},\lambda_{\beta\alpha}),\leqno (2.1.6)
$$
where $(\phi_{\beta\alpha},\lambda_{\beta\alpha})\in Tran(U_\alpha,U_\beta)$,
etc., any one of the three transition maps is uniquely determined by the other
two. As for the solvability of $(2.1.6)$, $(\phi_{\gamma\alpha},
\lambda_{\gamma\alpha})$ exists as postulated by {\em 2-d)} of Definition 
2.1.2, as long as the domain of $(\phi_{\beta\alpha},\lambda_{\beta\alpha})$ 
intersects the range of $(\phi_{\gamma\beta},\lambda_{\gamma\beta})$. 
Furthermore, $(\phi_{\beta\alpha},\lambda_{\beta\alpha})$ exists as long as
the range of $(\phi_{\gamma\beta},\lambda_{\gamma\beta})$ intersects the range
of $(\phi_{\gamma\alpha},\lambda_{\gamma\alpha})$, and 
$(\phi_{\gamma\beta},\lambda_{\gamma\beta})$ exists as long as the domain of 
$(\phi_{\beta\alpha},\lambda_{\beta\alpha})$ intersects the domain of 
$(\phi_{\gamma\alpha},\lambda_{\gamma\alpha})$. In particular, it follows that
$Tran(U_\alpha,U_\beta)=\{(\phi,\lambda)|(\phi,\lambda)^{-1}:=(\phi^{-1},
\lambda^{-1})\in Tran(U_\beta,U_\alpha)\}$. 

\hfill $\Box$

\noindent{\bf Remark 2.1.4 c:}\hspace{2mm}
In order to minimize the dependence of orbispace structure on the choice of the
set $\U$ of basic open sets so that it becomes more intrinsic, it is appropiate
to introduce an equivalence relation as follows: Two orbispace structures 
are said to be {\it directly equivalent} if one is contained in the other. 
The equivalence relation to be introduced is just a finite chain of direct
equivalence.

On the other hand, for any connected open subset $W$ of a basic open set $U$, 
we can assign a $G$-structure on $W$ by fixing a choice of the $G$-structures 
induced from the $G$-structure of $U$, and define a set of transition maps 
$Tran(W,V)$ for any basic open set $V$ by restricting each element of
$Tran(U,V)$ to the $G$-structure on $W$. It is easily seen that the given set
$\U$ of basic open sets can be enlarged by adding all the connected open
subsets in this manner, which will be assumed throughout this paper.

With these understood, let us consider the following situation: Suppose 
$(X,\U)$ is an orbispace which has a cover $\{U_\alpha,\alpha\in\Lambda\}$, 
$U_\alpha\in\U$, such that for each $\alpha\in\Lambda$, $\pi_{U_\alpha}:
\widehat{U_\alpha}\rightarrow U_\alpha$ is a covering map. 
Then by adding all the elementary neighborhoods of each $U_\alpha$ 
(w.r.t. the covering map $\pi_{U_\alpha}$) to $\U$, it is easily seen 
that the orbispace structure on $(X,\U)$ is equivalent to the trivial one.

\hfill $\Box$

\noindent{\bf Remark 2.1.4 d:}\hspace{2mm}
Each set of transition maps $Tran(U_\alpha,U_\beta)$ carries a natural
topology as follows. By the axiom {\em 2-b)} in Definition 2.1.2, for any 
transition maps $(\phi_i,\lambda_i)\in Tran(U_\alpha,U_\beta)$, $i=1,2$,
which are isomorphisms between induced $G$-structures of the same 
connected component of $U_\alpha\cap U_\beta$, there are 
$g_\alpha\in G_{U_\alpha}$, $g_\beta\in G_{U_\beta}$ such that 
$(\phi_2,\lambda_2)=g_\beta^{-1}\circ (\phi_1,\lambda_1)\circ g_\alpha$,
where $g_\alpha$, $g_\beta$ are regarded as elements of the automorphism groups
of the $G$-structures of $U_\alpha$ and $U_\beta$ respectively. 
This structure allows $Tran(U_\alpha,U_\beta)$ to inherit a topology from the
topological groups $G_{U_\alpha}$ and $G_{U_\beta}$. 
Recall that for any basic open sets $W,U$ such that $W\subset U$, the
space of cosets $G_U/G_W$ has a discrete topology. Hence each space of
transition maps $Tran(U_\alpha,U_\beta)$ is partitioned into a
disjoint union of copies of $G_{W_i}$, where 
$\{W_i|i\in\Lambda_{\alpha\beta}\}$ is the set of connected components
of $U_\alpha\cap U_\beta$. We will call this partition the {\it
canonical partition} of $Tran(U_\alpha, U_\beta)$. Under this topology, 
the composition of transition maps postulated in {\em 2-d)} of Definition 
2.1.2 is continuous.

\hfill $\Box$

\noindent{\bf Remark 2.1.4 e:}\hspace{2mm}
Let $(X,\U)$ be an orbispace. For each basic open set $U\in\U$, 
let $K_U$ be the subgroup of $G_U$ defined by 
$K_U=\{g|g\in G_U, \hspace{2mm}g\cdot x=x, \hspace{2mm}
\forall x\in\widehat{U}\}$. It is easy to see that $K_U$ is 
a normal subgroup and for any basic open set $W\subset U$, 
$K_U$ is contained in the image of $K_W$ under any injection 
in $Tran(W,U)$. We will call $(X,\U)$
{\it normal} if $K_W$ is sent isomorphically to $K_U$ under any injection in 
$Tran(W,U)$. For a normal orbispace $(X,\U)$, the isomorphism class of 
$K_U$, $U\in\U$, depends only on the connected component of $X$ that
contains $U$. For a connected normal orbispace $(X,\U)$, 
we will call the abstract group which is isomorphic to all $K_U$ the 
{\it kernel} of $(X,\U)$, and denote it by $K_{(X,\U)}$, or simply
$K_X$. A normal orbispace $(X,\U)$ is called {\it reduced} if $K_U$ 
is trivial for all $U\in\U$. For any normal orbispace $X$, 
there is a reduced orbispace canonically associated to it, which 
is obtained by replacing the $G$-structure $(\widehat{U},G_U,\pi_U)$ 
of each basic open set $U$ with a new $G$-structure
$(\widehat{U},G_U/K_U,\pi_U)$. We will call the resulting orbispace the
{\it canonical reduction} of $X$, and denote it by $X^{red}$.

Consider the following example. Let $X=S^1=\R/\Z$. Set $I_1=(-\frac{1}{8},
\frac{5}{8})$, $I_2=(\frac{3}{8},\frac{9}{8})$, then $I_1,I_2$ give rise
to two open sets in $S^1$, denoted by $U_1$ and $U_2$ respectively, such that
$S^1=U_1\cup U_2$ and $U_1\cap U_2$ has two connected components, $W_1$, 
$W_2$, with $W_1=(\frac{3}{8},\frac{5}{8})$. On the other hand, suppose $G$ 
is a group and there is an automorphism
$\tau:G\rightarrow G$ which is not an inner automorphism, i.e., for any 
$g\in G$, $\tau(h)\neq ghg^{-1}$ for some $h\in G$. Now we give $U_i$, $i=1,2$,
a $G$-structure $(U_i,G,\pi_i)$ where $G$ acts on $U_i$ trivially. We specify
the set of transition maps by 
$$
Tran(U_1,U_2):=\{h^{-1}\circ (Id_{W_1},Id_{G})\circ g,\hspace{1.5mm}
h^{-1}\circ (Id_{W_2},\tau^{-1})\circ g|h,g\in G\}. \leqno (2.1.7)
$$
It is easily seen that there is an orbispace structure on $S^1$ generated by
$U_1,U_2$ with the given $Tran(U_1,U_2)$. We denote it by $\U_\tau$. 
The orbispace $(S^1,\U_\tau)$ is clearly normal, with kernel 
$K_{(S^1,\U_\tau)}=G$. However, we wish to point out that $K_{(S^1,\U_\tau)}$ 
only lives as an abstract group. For example, there is no compatible action
of $K_{(S^1,\U_\tau)}$ on $S^1$, because $(S^1,\U_\tau)$ is not equivalent to 
the orbispace $(S^1,\U_0)$ derived from the $G$-space $(S^1,G)$ where
$G$ acts trivially. In fact, we will show that the fundamental group
of $(S^1,\U_0)$ is $\pi_0(G)\times \Z$ while the fundamental group of 
$(S^1,\U_\tau)$ is the semi-direct product of $\pi_0(G)$ by $\Z$ with respect
to the homomorphism $\Z\rightarrow Aut(\pi_0(G))$ given by 
$1\mapsto \tau_\ast$, where $\tau_\ast:\pi_0(G)\rightarrow \pi_0(G)$ is 
induced by $\tau:G\rightarrow G$ (cf. Theorem 3.2.6, Example 4.2.9).

\hfill $\Box$

In a certain sense, the definition of orbispace structure in this paper is 
of a less intrinsic style than the definition of orbifold structure 
originally adopted in \cite{CR1}, where the language of germs was
used. The reason for which we 
choose such a formulation lies in the following facts: First of all, 
an automorphism of a $G$-structure $(\widehat{U},G_U,\pi_U)$, i.e., a pair 
$(\phi,\lambda)$ where $\lambda\in Aut(G_U)$ and $\phi$ is a
$\lambda$-equivariant homeomorphism of $\widehat{U}$ inducing identity map
on $U$, might not be induced by an action of $G_U$ on $\widehat{U}$ (at least
we have a burden of proving so). But we want to only consider those arising 
from the action of $G_U$. Secondly, suppose $W$ is a connected open subset 
of $U$, both of which are given with $G$-structures $(\widehat{W},G_W,\pi_W)$ 
and $(\widehat{U},G_U,\pi_U)$ respectively. Then an open embedding 
$\phi:\widehat{W}\rightarrow \widehat{U}$, which is $\lambda$-equivariant for 
some monomorphism $\lambda:G_W\rightarrow G_U$ and induces the inclusion 
$W\hookrightarrow U$, might not be an isomorphism onto one of the 
$G$-structures of $W$ induced from $(\widehat{U},G_U,\pi_U)$. However, this 
kind of pathology can be ruled out if each $\widehat{U}$ is Hausdorff,
$G_U$ acts on $\widehat{U}$ effectively, and for any convergent sequence 
$\{g_n\cdot x_n\}$ in $\widehat{U}$ where $g_n\in G_U$ and 
$x_n\in\widehat{U}$, there is a convergent subsequence of $\{g_n\}$ in $G_U$.
These conditions are clearly satisfied in the case of orbifolds.

\vspace{1.5mm}

We end this subsection by including such a pathological example.

\vspace{2mm}

\noindent{\bf Example 2.1.5:}\hspace{2mm}
Let $D_1, D_2$ be two copies of the unit disc in the plane of complex 
numbers, and $\rho:D_1\setminus\{0\}\rightarrow D_2\setminus\{0\}$ be the 
identity map. We define $Y=(D_1\sqcup D_2)/\{z=\rho(z)\}$ to be the space 
obtained by gluing $D_1$ and $D_2$ via $\rho$. Then the identity map 
$\tau:D_1\rightarrow D_2$ generates an involution on $Y$. Set 
$G=<\tau>\cong \Z_2$. We will analyze the global orbispace $X:=Y/G$ 
canonically constructed from the $G$-space $(Y,G)$. 

We identify the underlying topological space of $X$ with the unit 
disc $D$ in the plane of complex numbers. For any connected open subset $U$ 
of $D$, if $U$ does not contain $0$, then in its $G$-structure 
$(\widehat{U},G_U,\pi_U)$, we have $\widehat{U}=U$ with $G_U=\Z_2$ acting 
trivially. If $U$ contains $0$, then $\widehat{U}$ is obtained by gluing two 
copies of $U$ along the non-zero numbers in $U$, with $G_U=\Z_2$ acting 
non-trivially, which is generated by the identity map between the two
copies of $U$. 

On the other hand, each connected open subset $U$ of $D$ has
a trivial $G$-structure, and there exist equivariant open embeddings into the
$G$-structures induced from $(Y,G)$. Clearly none of them is an isomorphism
onto. As a matter of fact, if these trivial $G$-structures were allowed, the
orbispace $X$ would have lost the non-trivial local $G$-structures 
inherited from $(Y,G)$ as its orbit space.

Finally, we point out that the orbispace $X$ is clearly not normal.

\subsection{Morphisms between orbispaces}

\hspace{5mm}Let $(X,\U)$, $(X^\prime,\U^\prime)$ be two orbispaces. A 
{\it morphism}
from $(X,\U)$ to $(X^\prime,\U^\prime)$ will be defined as an equivalence 
class of a system of local equivariant continuous maps satisfying certain 
compatibility conditions. Such a system consists of following ingredients: 
(1) a collection of basic open sets of $X$, $\{U_\alpha|\alpha\in\Lambda\}$, 
which forms a cover of $X$, (2) a collection of basic open sets of $X^\prime$, 
$\{U_\alpha^\prime|\alpha\in\Lambda\}$, which is labelled by the same indexing
set $\Lambda$ (here $U_\alpha^\prime$ may coincide with
$U_\beta^\prime$ with $\alpha\neq\beta$), (3) a collection of 
continuous maps $\{f_\alpha:\widehat{U_\alpha}\rightarrow 
\widehat{U_\alpha^\prime}|\alpha\in\Lambda\}$ defined between the 
$G$-structures of $U_\alpha$ and $U_\alpha^\prime$ for each 
$\alpha\in\Lambda$, and (4) a collection 
$\{\rho_{\beta\alpha}|\alpha,\beta\in\Lambda\}$ where each 
$\rho_{\beta\alpha}$ is a continuous map from $Tran(U_\alpha,U_\beta)$ to 
$Tran(U_\alpha^\prime,U_\beta^\prime)$. The maps $\{f_\alpha\}$ and 
$\{\rho_{\beta\alpha}\}$ are required to satisfy the following compatibility 
conditions:
$$
f_\beta\circ\imath=\rho_{\beta\alpha}(\imath)\circ f_\alpha, \forall 
\alpha,\beta\in\Lambda,\imath\in Tran(U_\alpha,U_\beta)\leqno (2.2.1 a)
$$
$$
\rho_{\gamma\alpha}(\jmath\circ\imath)=\rho_{\gamma\beta}(\jmath)\circ
\rho_{\beta\alpha}(\imath), \forall \alpha, \beta, \gamma\in\Lambda, \imath\in 
Tran(U_\alpha,U_\beta), \jmath\in Tran(U_\beta,U_\gamma). \leqno (2.2.1 b)
$$
In particular, for each $\alpha\in\Lambda$, 
$\rho_{\alpha\alpha}:Tran(U_\alpha,U_\alpha)\rightarrow 
Tran(U_\alpha^\prime,U_\alpha^\prime)$ is actually a homomorphism 
$G_{U_\alpha}\rightarrow G_{U_\alpha^\prime}$ by $(2.2.1 b)$, which will be 
abbreviated by $\rho_\alpha$. Equation $(2.2.1 a)$ then implies that each 
$f_\alpha$ is $\rho_\alpha$-equivariant. Such a system will be denoted by 
$(\{U_\alpha\},\{U_\alpha^\prime\},\{f_\alpha\},\{\rho_{\beta\alpha}\})$. We 
will say that the maps $\{f_\alpha\}$, $\{\rho_{\beta\alpha}\}$ are defined 
over $(\{U_\alpha\},\{U_\alpha^\prime\})$. In a context where no confusion 
arises, we will abbreviate the notation of such a system by 
$(\{f_\alpha\},\{\rho_{\beta\alpha}\})$. Clearly each of such systems induces 
a continuous map $f:X_{top}\rightarrow X_{top}^\prime$ between the underlying 
topological spaces. 

\vspace{1.5mm}

\noindent{\bf Definition 2.2.1:}
{\em Two systems of maps $(\{f_{\alpha,1}\},\{\rho_{\beta\alpha,1}\})$
and $(\{f_{\alpha,2}\},\{\rho_{\beta\alpha,2}\})$ defined over 
$(\{U_\alpha\},\{U_\alpha^\prime\})$ are said to be isomorphic if there is a 
collection $\{\delta_\alpha\}$ where each $\delta_\alpha$ is an automorphism 
of the $G$-structure 
$(\widehat{U_\alpha^\prime},G_{U_\alpha^\prime},\pi_{U_\alpha^\prime})$ on 
$U_\alpha^\prime$, such that the following conditions hold:
$$
f_{\alpha,2}=\delta_\alpha\circ f_{\alpha,1}, \forall \alpha\in\Lambda \leqno 
(2.2.2 a)
$$
$$
\rho_{\beta\alpha,2}(\imath)=\delta_\beta\circ\rho_{\beta\alpha,1}(\imath)
\circ (\delta_\alpha)^{-1}, \forall \alpha, \beta\in\Lambda, \imath\in 
Tran(U_\alpha,U_\beta). \leqno (2.2.2 b)
$$
}

We remark that isomorphic systems induce the same continuous map between the 
underlying topological spaces.

\vspace{1.5mm}

The equivalence relation to be introduced on these systems of maps, which will
eventually define a morphism between orbispaces, involves a procedure of taking
refinements of a given system of maps. We shall describe this procedure first.

Let $\sigma=(\{U_\alpha\},\{U_\alpha^\prime\},\{f_\alpha\},
\{\rho_{\beta\alpha}\})$ be a system satisfying $(2.2.1a)$ and $(2,2.1b)$, 
which induces a continuous map $f:X_{top}\rightarrow X_{top}^\prime$. 
Suppose we are given a cover $\{U_i|U_i\in\U, i\in I\}$ of $X$ which is 
a refinement of $\{U_\alpha|\alpha\in\Lambda\}$, i.e., there is a mapping
$\theta:I\rightarrow \Lambda$ between the indexing sets such that $U_i\subset 
U_{\theta(i)}$ for each $i\in I$, and a collection of basic open sets of 
$X^\prime$, $\{U_i^\prime|i\in I\}$, which is labelled by the same
indexing set $I$ (here $U_i^\prime$ may coincide with $U_j^\prime$
with $i\neq j$) and satisfies $f(U_i)\subset U_i^\prime$ for 
each $i\in I$.

\vspace{1.5mm}

\noindent{\bf Lemma 2.2.2:}
{\em The system 
$\sigma=(\{U_\alpha\},\{U_\alpha^\prime\},
\{f_\alpha\},\{\rho_{\beta\alpha}\})$ 
canonically induces a family of systems defined over 
$(\{U_i\},\{U_i^\prime\})$ which are isomorphic to each other in the 
sense of Definition 2.2.1.      
}
\vspace{2mm}

\noindent{\bf Proof:}
We first construct the continuous map $f_i:\widehat{U}_i\rightarrow
\widehat{U}_i^\prime$ for each index $i\in I$. By assumption $U_i\subset
U_{\theta(i)}$, there is an injection $\imath\in Tran(U_i,U_{\theta(i)})$.
We take the composition $f_{\theta(i)}\circ \imath: \widehat{U}_i\rightarrow
\widehat{U}_{\theta(i)}^\prime$. Since $f(U_i)\subset U_i^\prime$, 
it follows that $f(U_i)\subset U_i^\prime\cap U_{\theta(i)}^\prime\neq 
\emptyset$, and hence there exists a transition map
$\imath^\prime\in Tran(U_i^\prime,U_{\theta(i)}^\prime)$ such that 
$f_{\theta(i)}\circ \imath (\widehat{U}_i)\subset Im\;\imath^\prime$. 
We define $f_i:\widehat{U}_i\rightarrow \widehat{U}_i^\prime$ by setting
$$
f_i:=(\imath^\prime)^{-1}\circ f_{\theta(i)}\circ\imath. \leqno (2.2.3a)
$$ 
The map $f_i$ depends on the choices of transition maps $\imath$ and 
$\imath^\prime$. However, one is ready to verify 
that any two such maps factor through an automorphism of 
$(\widehat{U}_i^\prime, G_{U_i^\prime},\pi_{U_i^\prime})$.

Suppose we have fixed a choice of transition maps $\imath$,
$\imath^\prime$ and the map $f_i$ for each $i\in I$. We will now 
define a continuous map $\xi_{ji}:Tran(U_i,U_j)\rightarrow Tran(U_i^\prime,
U_j^\prime)$ for each pair of indices $(i,j)$. To fix the notation, let us
assume that $f_i:=(\imath^\prime)^{-1}\circ f_{\theta(i)}\circ\imath$,
$f_j:=(\jmath^\prime)^{-1}\circ f_{\theta(j)}\circ\jmath$ for some fixed 
choices of transition maps $\imath,\imath^\prime,\jmath,\jmath^\prime$.
For any $a\in Tran(U_i,U_j)$, set $\theta(a)=\jmath\circ a \circ\imath^{-1}$ 
which is in $Tran(U_{\theta(i)},U_{\theta(j)})$. We then define 
$\xi_{ji}:Tran(U_i,U_j)\rightarrow Tran(U_i^\prime,U_j^\prime)$ by setting
$$
\xi_{ji}(a):=(\jmath^\prime)^{-1}\circ\rho_{\theta(j)\theta(i)}(\theta(a))
\circ\imath^\prime \leqno (2.2.3b)
$$
in $Tran(U_i^\prime,U_j^\prime)$. It is not hard to see that $(2.2.1a-b)$ are
satisfied for $(\{U_i\},\{U_i^\prime\},\{f_i\},\{\xi_{ji}\})$, and two 
different such induced systems are related by $(2.2.2a-b)$. 

Finally, we observe that different choices for the index mapping $\theta:
I\rightarrow\Lambda$ result in isomorphic induced systems.

\hfill $\Box$

\noindent{\bf Definition 2.2.3:}
{\em We call each induced system $(\{f_i\},\{\xi_{ji}\})$ a refinement
of $(\{f_\alpha\},\{\rho_{\beta\alpha}\})$ associated to 
$(\{U_i\},\{U_i^\prime\})$. Systems $\sigma=(\{U_\alpha\},
\{U_\alpha^\prime\},\{f_\alpha\},\{\rho_{\beta\alpha}\})$ and 
$\tau=(\{U_a\},\{U_a^\prime\},\{g_a\},\{\xi_{ba}\})$ are said to be directly 
equivalent if there exists a $(\{U_i\},\{U_i^\prime\})$, of which the 
associated refinements of the systems $\sigma$ and $\tau$ are isomorphic.
}

\hfill $\Box$

We remark that any refinement of a given system induces the same continuous
map between the underlying topological spaces, so do any two directly
equivalent systems.

\vspace{1.5mm}

\noindent{\bf Lemma 2.2.4:}
{\em Direct equivalence is an equivalence relation.}

\vspace{2mm}

\noindent{\bf Proof:}
We shall prove that any two refinements of a given system have a common 
refinement. The lemma then follows from the fact that refinement relation is 
transitive.

Suppose two systems $\sigma=(\{U_\alpha\},\{U_\alpha^\prime\},\{f_\alpha\},
\{\rho_{\beta\alpha}\})$ and $\tau=(\{U_a\},\{U_a^\prime\},\{g_a\},
\{\xi_{ba}\})$ are refinements of a given system $\kappa$. In order to 
construct a common refinement of $\sigma$ and $\tau$, we first take a common
refinement $\{U_i|i\in I\}$ of $\{U_\alpha\}$ and $\{U_a\}$. We can arrange so
that for each $i\in I$, there is a basic open set $U_i^\prime$ of $X^\prime$
such that $f(U_i)\subset U^\prime_i$, where $f:X_{top}\rightarrow 
X_{top}^\prime$ is the continuous map induced by $\kappa$. By Lemma 2.2.2,
there are refinements of $\sigma$ and $\tau$ associated to 
$(\{U_i\},\{U_i^\prime\})$, which must also be refinements of $\kappa$ because
refinement relation is transitive. Hence they are isomorphic, and therefore
$\sigma$ and $\tau$ have a common refinement.

\hfill $\Box$

\noindent{\bf Definition 2.2.5:}
{\em A morphism from $(X,\U)$ to $(X^\prime,\U^\prime)$ is an equivalence 
class of systems 
$$
(\{U_\alpha\},\{U_\alpha^\prime\},\{f_\alpha\},\{\rho_{\beta\alpha}\})
$$
satisfying $(2.2.1a-b)$. 
}

\hfill $\Box$

As a notational convention, we will write the equivalence class of 
$(\{U_\alpha\},\{U_\alpha^\prime\},\{f_\alpha\},\{\rho_{\beta\alpha}\})$
as $(\tilde{f},\tilde{\rho})$ or simply $\tilde{f}$, and the induced continuous
map between underlying topological spaces as $f$.

\vspace{1.5mm}

Suppose $X$ and $X^\prime$ are global orbispaces which are defined canonically
from $G$-spaces $(Y,G)$ and $(Y^\prime,G^\prime)$ respectively. Then any pair 
$(f,\rho)$, where $\rho:G\rightarrow G^\prime$ is a homomorphism and 
$f:Y\rightarrow Y^\prime$ is a $\rho$-equivariant continuous map, defines a 
morphism $X\rightarrow X^\prime$. However, simple examples show that not 
every morphism $X\rightarrow X^\prime$ arises in this way. For instance, taking
$S^1$ as an orbispace with a trivial orbispace structure, a morphism from $S^1$
into an orbispace $X$ will be called a {\it free loop} in $X$. When $X$ is
global given by a $G$-space $(Y,G)$, one can show that a free loop in $X$ is
just a conjugate class of a pair $(\gamma,g)$ where $g\in G$ and 
$\gamma:[0,1]\rightarrow Y$ is a path satisfying
$\gamma(1)=g\cdot\gamma(0)$ (cf. Lemma 3.5.1).
This is precisely the so-called {\it twisted loop} in physics literature.
Certainly, one can regard $S^1$ as a trivial $G$-space, and easily see that
a free loop in $X$ does not necessarily arise as a global equivariant map.

On the other hand, in the case of orbifolds, a class of $C^\infty$ maps 
was singled out in \cite{CR2}, which was called {\it good maps}, 
by the property that pull-back orbibundles are well-defined. A natural 
isomorphism relation was introduced on the pull-back orbibundles. 
One can show that a ($C^\infty$) morphism between orbifolds is
equivalent to giving a good map together with an isomorphism class of 
pull-back orbibundles \footnote{an equivalent notion under the 
name ``strong maps'' had been introduced earlier, cf. \cite{Pr}}. 

Finally, as we pointed out in Example 2.1.3c, the notion of
homomorphisms between complexes of groups is equivalent to the notion
of morphisms between the associated orbihedra, viewed canonically as
orbispaces in the sense of this paper.  

\vspace{2mm}

\noindent{\bf Lemma 2.2.6:}
{\em Any two morphisms $\tilde{f}:X\rightarrow Y$ and 
$\tilde{g}:Y\rightarrow Z$ canonically determine a morphism, 
denoted by $\tilde{g}\circ\tilde{f}:X\rightarrow Z$,
which satisfies the associativity law: 
$\tilde{h}\circ (\tilde{g}\circ\tilde{f})=
(\tilde{h}\circ \tilde{g})\circ\tilde{f}$. 
}

\vspace{2mm}

\noindent{\bf Proof:}
Suppose the morphisms $\tilde{f}:X\rightarrow Y$ and 
$\tilde{g}:Y\rightarrow Z$ are represented by systems 
$\sigma=(\{U_\alpha\},\{V_\alpha\},\{f_\alpha\},\{\rho_{\beta\alpha}\})$ and 
$\tau=(\{V_a\},\{W_a\},\{g_a\},\{\eta_{ba}\})$ respectively. Consider a 
collection of basic open sets of $X$, $\{U_{\alpha a}\}$, which consists of 
all the connected components of $U_\alpha\cap f^{-1}(V_a)$ for all $\alpha,a$. 
Obviously $\{U_{\alpha a}\}$ is a refinement of $\{U_\alpha\}$. We take 
$V_{\alpha a}$ to be an appropriated component of $V_\alpha\cap V_a$ 
satisfying $f(U_{\alpha a})\subset V_{\alpha a}$. By Lemma 2.2.2, there is a 
refinement $\sigma^\prime$ of $\sigma$ defined over 
$(\{U_{\alpha a}\},\{V_{\alpha a}\})$. 
On the other hand, consider the collection of basic open sets of $Y$, 
$\{V_{\alpha a}\}$, which consists of all the connected components of 
$V_\alpha\cap V_a$ for all $\alpha,a$. Obviously $\{V_{\alpha a}\}$ is a 
refinement of $\{V_a\}$. We take $W_{\alpha a}$ to be $W_a$ which satisfies 
$g(V_{\alpha a})\subset W_a$. Again by Lemma 2.2.2, we have a refinement 
$\tau^\prime$ of $\tau$ which is defined over 
$(\{V_{\alpha a}\},\{W_{\alpha a}\})$. We define the morphism 
$\tilde{g}\circ\tilde{f}:X\rightarrow Z$ to be the equivalence class of the 
system $\tau^\prime\circ\sigma^\prime$ obtained by composing $\sigma^\prime$ 
with $\tau^\prime$.

In order to see that $\tilde{g}\circ\tilde{f}$ is well-defined, we simply 
observe that if we change $\sigma^\prime$ or $\tau^\prime$ by isomorphic 
systems, the resulting composition $\tau^\prime\circ\sigma^\prime$ will be 
changed to an isomorphic one, and if we replace $\sigma$ or $\tau$ by a 
refinement, we will obtain a refinement of $\tau^\prime\circ\sigma^\prime$. 

Finally, the associativity of composition follows easily from the nature of 
the construction. We leave the details for readers.

\hfill $\Box$

The morphism $\tilde{g}\circ\tilde{f}$ will be called the {\it composition}
of $\tilde{f}$ with $\tilde{g}$.

\vspace{2mm}

\noindent{\bf Theorem 2.2.7:}
{\em The set of all orbispaces forms a category.}

\subsection{Based orbispaces and pseudo-embedding}

\hspace{5mm}The construction of homotopy theory requires introducing a 
{\it base-point 
structure} to each orbispace. An orbispace equipped with a base-point 
structure will be called a {\it based orbispace}. In fact, there is a based 
version of the category of orbispaces, in which each object is a based 
orbispace and each morphism is an equivalence class of systems of maps which 
are required to preserve the base-point structures. We shall only sketch the 
construction since it is parallel to the case without base-point structures.

Let $(X,\U)$ be an orbispace. A base-point structure of $(X,\U)$ is a triple 
$(o,U_o,\hat{o})$ where $o$ is a point in $X$, $U_o\in\U$ is a basic open set 
containing $o$, and $\hat{o}$ is a point in the inverse image 
$\pi_{U_o}^{-1}(o)$. We will write $\underline{o}=(o,U_o,\hat{o})$ for the 
base-point structure, and $(X,\underline{o})$ for the based orbispace.

We sketch the definition of a morphism between based orbispaces next. Let 
$(X,\underline{o})$ and $(X^\prime,\underline{o^\prime})$ be based orbispaces, 
where $\underline{o}=(o,U_o,\hat{o})$ and 
$\underline{o^\prime}=(o^\prime,U_{o}^\prime,\hat{o}^\prime)$. A {\it based 
system} from $(X,\underline{o})$ to $(X^\prime,\underline{o^\prime})$ is just 
an ordinary system 
$(\{U_\alpha\},\{U_\alpha^\prime\},\{f_\alpha\},\{\rho_{\beta\alpha}\})$ 
satisfying $(2.2.1a-b)$ together with the following additional requirements: 
$U_o\in\{U_\alpha\}$, $U_o^\prime\in\{U_\alpha^\prime\}$, $U_o\mapsto 
U_o^\prime$ under the correspondence $U_\alpha\mapsto U_\alpha^\prime$, and 
the map $f_o:\widehat{U_o}\rightarrow \widehat{U_o^\prime}$ satisfies the 
equation $f_o(\hat{o})=\hat{o}^\prime$. An isomorphism relation between based 
systems is defined as in the ordinary case in Definition 2.2.1, but with an 
additional condition that the automorphism $\delta_o$ of the $G$-structure 
$(\widehat{U_o^\prime},G_{U_o^\prime},\pi_{U_o^\prime})$ of $U_o^\prime$ must 
be $1_{G_{U_o^\prime}}\in G_{U_o^\prime}$. As for the construction of 
refinement of a based system 
$(\{U_\alpha\},\{U_\alpha^\prime\},\{f_\alpha\},\{\rho_{\beta\alpha}\})$, we 
assume that we are given a refinement $\{U_i\}$ of $\{U_\alpha\}$ satisfying 
the condition that $U_o\in \{U_i\}$ and the index mapping $\alpha=\theta(i)$ 
which defines the refinement relation obeys $\theta(o)=o$, and moreover, 
$U_o^\prime\in \{U_i^\prime\}$ and the map $f_o:\widehat{U_o}\rightarrow 
\widehat{U_o^\prime}$ is unchanged throughout. Parallel to the case without 
base-point structures, a direct equivalence relation on based systems can be 
defined, which is also an equivalence relation. An equivalence class of based 
systems will be called a {\it based morphism}. One can similarly show that the 
composition of two based morphisms is well-defined. 

\vspace{1.5mm}

The following theorem is parallel to Theorem 2.2.7.

\vspace{1.5mm}

\noindent{\bf Theorem 2.3.1:}
{\em The set of based orbispaces with based morphisms forms a category.}

\vspace{2mm}

In what follows, we will introduce several definitions regarding {\it 
sub-orbispaces, orbispace embedding, pseudo-embedding}, and {\it cartesian
product}. The notion of pseudo-embedding is involved in defining relative 
homotopy groups in section 3.3.

Let $(X,\U)$ be an orbispace, $Y$ be a locally connected subspace of the 
underlying topological space $X$ with induced topology, such that for each 
$\U_\alpha\in\U$, the subspace $\pi_{U_\alpha}^{-1}(Y\cap U_\alpha)$ in 
$\widehat{U_\alpha}$ is locally connected. We will canonically construct an 
orbispace structure on $Y$ from the given orbispace structure on $X$ as 
follows. The set of basic open sets of $Y$, which will be denoted by 
$\V=\{V_i|i\in I\}$, is defined to be the collection of all the connected 
components of $Y\cap U_\alpha$ for all $U_\alpha\in \U$. The $G$-structure 
$(\widehat{V}_i,G_{V_i},\pi_{V_i})$ on $V_i$ is given in the following 
way. Suppose $V_i$ is a connected component of $Y\cap U_\alpha$ for some 
$\alpha$. Then $\widehat{V}_i$ is taken to be a connected component 
of $\pi_{U_\alpha}^{-1}(V_i)$ in $\widehat{U}_\alpha$, $G_{V_i}$ to be the 
subgroup of $G_{U_\alpha}$ consisting of elements $g$ such that  
$g\cdot \widehat{V}_i=\widehat{V}_i$, and 
$\pi_{V_i}=\pi_{U_\alpha}|_{\widehat{V}_i}$. 
As for transition maps, let $V_i$, $V_j$ be two basic open sets in $\V$ 
such that $V_i\cap V_j\neq\emptyset$, and suppose $V_i$, $V_j$ are connected 
components of $Y\cap U_\alpha$ and $Y\cap U_\beta$ respectively. 
We define $Tran(V_i,V_j)$ to be the set of all restrictions of transition maps 
in $Tran(U_\alpha,U_\beta)$ on $\widehat{V}_i$. One can verify that this 
indeed defines an orbispace structure on $Y$. The orbispace $(Y,\V)$ is 
called a {sub-orbispace} of $(X,\U)$. Clearly different choices for the 
$G$-structure $(\widehat{V}_i, G_{V_i},\pi_{V_i})$ on $V_i$ will result in 
isomorphic sub-orbispace structures. As examples, any open subset of an
orbispace is canonically a sub-orbispace.

There is a canonical way to assign a base-point structure to a sub-orbispace
of a based orbispace. If $X$ is given a base-point structure 
$\underline{o}=(o,U_o,\hat{o})$ such that $o\in Y$, we can give a base-point 
structure $(o,V_o,\hat{o})$ to $Y$, where $V_o$ is the connected component 
of $Y\cap U_o$ containing $o$. The $G$-structure of $V_o$ is taken to be 
the induced $G$-structure that contains $\hat{o}$. The sub-orbispace 
$Y$ with the canonical base-point structure $(o,V_o,\hat{o})$ is called a 
{\it based sub-orbispace} of $(X,\underline{o})$.

A morphism $\tilde{i}:Y\rightarrow X$ is called an {\it 
orbispace embedding} if $\tilde{i}$ is an isomorphism onto the sub-orbispace 
$i(Y)$ of $X$, where $i:Y_{top}\rightarrow X_{top}$ is the continuous map 
induced by $\tilde{i}$. There is a notion of {\it orbifold embedding} defined 
as follows: a $C^\infty$ map betweem 
orbifolds, $f:X\rightarrow X^\prime$, is called an orbifold embedding if 
there are compatible local smooth liftings $f_p:V_p\rightarrow V_{f(p)}^\prime$
between uniformizing systems and isomorphisms $\rho_p:G_p\rightarrow 
G_{f(p)}^\prime$ such that each $f_p$ is a $\rho_p$-equivariant embedding.
It is easy to see that orbispace embedding is equivalent to orbifold embedding 
in the case of orbifolds.

There is a more general notion of orbispace embedding. A morphism 
$\tilde{i}:Y\rightarrow X$ is called a {\it pseudo-embedding}, 
if there is a representing system 
$\sigma=(\{V_\alpha\},\{U_\alpha\},\{i_\alpha\},\{\rho_{\beta\alpha}\})$, 
in which each $i_\alpha:\widehat{V}_\alpha\rightarrow\widehat{U}_\alpha$ is 
an embedding, and each $\rho_\alpha:G_{V_\alpha}\rightarrow G_{U_\alpha}$ is 
a monomorphism such that the natural projection 
$G_{U_\alpha}\rightarrow G_{U_\alpha}/\rho_\alpha(G_{V_\alpha})$ is a weak 
fibration. The injectivity of each $\rho_\alpha$ implies that the
restriction of $\rho_{\beta\alpha}:Tran(V_\alpha,V_\beta)\rightarrow
Tran(U_\alpha, U_\beta)$ on each component of the canonical partition
of $Tran(V_\alpha,V_\beta)$ (cf. Remark 2.1.4 d) is an embedding for
any $\alpha,\beta$. 
Note that in the construction of sub-orbispace structure, each 
map $G_{U_\alpha}\rightarrow G_{U_\alpha}/G_{V_i}$ is a weak fibration 
because $G_{U_\alpha}/G_{V_i}$ has a discrete topology. Hence
orbispace embedding is a special case of pseudo-embedding. 
The based version of pseudo-embedding is similarly defined.

For any two orbispaces $(X,\U)$ and $(X^\prime,\U^\prime)$, the {\it
cartesian product}, which will be denoted by $(X\times
X^\prime,\U\times\U^\prime)$, is defined as follows. The orbispace
structure $\U\times\U^\prime$ consists of open sets $U\times U^\prime$
where $U\in\U$ and $U^\prime\in\U^\prime$, and the $G$-structure of
$U\times U^\prime$ is
$(\widehat{U}\times\widehat{U^\prime},G_{U}\times
G_{U^\prime},\pi_U\times\pi_{U^\prime})$ (connectedness and
locally-connectedness are preserved under cartesian product). As for
transition maps, we define
$$
Tran(U_\alpha\times U^\prime_{\alpha^\prime}, 
U_\beta\times U^\prime_{\beta^\prime}):=
Tran(U_\alpha,U_\beta)\times Tran(U^\prime_{\alpha^\prime},
U^\prime_{\beta^\prime}), \leqno (2.3.1)
$$
which clearly satisfy the axioms in Definition 2.1.2. The cartesian product 
is associative, and is defined for any finitely many orbispaces. If 
$\underline{o}=(o,U_o,\hat{o})$ and 
$\underline{o}^\prime=(o^\prime,U^\prime_o,\hat{o}^\prime)$ are base-point 
structures of $(X,\U)$ and $(X^\prime,\U^\prime)$  respectively, there is a 
canonical base-point structure $\underline{o}\times\underline{o}^\prime:
=((o,o^\prime),U_o\times U_o^\prime,(\hat{o},\hat{o}^\prime))$ of the 
cartesian product $(X\times X^\prime,\U\times\U^\prime)$. 

\vspace{1.5mm}

For an example of pseudo-embedding which is not an orbispace embedding, we
consider the cartesian product $X\times X^\prime$. Pick a point $o\in X$, and
a basic open set $U$ containing $o$. We further fix a point $\hat{o}\in 
\widehat{U}$ such that $\pi_U(\hat{o})=o$. We define a morphism 
$\tilde{i}:X^\prime\rightarrow X\times X^\prime$ by a system 
$(\{U_\alpha^\prime\},\{U\times U_\alpha^\prime\},\{i_\alpha\},
\{\rho_{\beta\alpha}\})$ where $\{U_\alpha^\prime\}$ is a cover of $X^\prime$
by basic open sets, $i_\alpha:\widehat{U_\alpha^\prime}\rightarrow\widehat{U}
\times\widehat{U_\alpha^\prime}$ is given by $x\mapsto (\hat{o},x)$, and
$\rho_{\beta\alpha}$ is defined by $\xi\mapsto (1_{G_U},\xi)$ for any
$\xi\in Tran(U_\alpha^\prime,U_\beta^\prime)$. The morphism $\tilde{i}:
X^\prime\rightarrow X\times X^\prime$ is clearly a pseudo-embedding, but 
certainly not an orbispace embedding when the isotropy group $G_o$ of $o$ is 
non-trivial, or $G_U$ does not have a discrete topology.

\section{Loop Spaces and Homotopy Groups}

\subsection{Based loop space}

\hspace{5mm}
Let $(X,\underline{o})$ be a based orbispace where $\underline{o}=
(o,U_o,\hat{o})$. The based loop space of $(X,\underline{o})$, which will be
denoted by $\Omega(X,\underline{o})$, is by definition the space of all 
``based morphisms'' from $S^1$ into $(X,\underline{o})$. However, the
``based morphism'' referred here is slightly different from the version
described in Section 2.3, hence we shall give it a detailed account next.

To begin with,  let $S^1=[0,1]/\{0\sim 1\}$. We will fix $0\in S^1$ as a base 
point, and write $(S^1,\ast)$ for the corresponding based topological space. 
In defining the based loop space, $(S^1,\ast)$ will be regarded as a based 
orbispace with a trivial orbispace structure, with the understanding that in 
its base-point structure $(\ast,I_0,\hat{\ast})$, the interval $I_0$ is 
allowed to be changed. 

Let $\{I_i| i=0, 1, \cdots, n\}$ be a cover of $(S^1,\ast)$ by open intervals 
such that $\ast\in I_0$ and $I_i\cap I_j\neq\emptyset$ only if $j=i+1$ or 
$j=n$ and $i=0$, and $\{U_i| i=0, 1, \cdots, n\}$ be a collection of basic 
open sets of $X$ where $U_0=U_o$. We shall consider systems such as  
$(\{I_i\},\{U_i\},\{f_i\},\{\rho_{ji}\})$, where each $f_i:I_i\rightarrow 
\widehat{U_i}$ is a continuous map and each $\rho_{ji}$ is in 
$Tran(U_i, U_j)$ where either $j=i+1$ or $j=n$, $i=0$. Furthermore, we 
require that on each $I_i\cap I_j$, equation $f_j=\rho_{ji}\circ f_i$ holds, 
and $f_0(\ast)=\hat{o}$ in $\widehat{U_o}$. An isomorphism relation on these 
systems is defined as in Definition 2.2.1, with an additional
requirement that the automorphism $\delta_o$ of the $G$-structure 
$(\widehat{U_o},G_{U_o},\pi_{U_o})$ be $1_{G_{U_o}}$ in $G_{U_o}$.  A system 
$(\{J_k\},\{V_k\},\{g_k\},\{\xi_{lk}\})$, where $k$ is running from $0$ to 
$m$, is said to be a refinement of $(\{I_i\},\{U_i\},\{f_i\},\{\rho_{ji}\})$ 
if (1) $\{J_k\}$ is a refinement of $\{I_i\}$ given by a mapping 
$\theta:\{0,\cdots,m\}\rightarrow \{0,\cdots,n\}$ such that $\theta(0)=0$,  
and (2) there are transition maps $\jmath_k\in Tran(V_k, U_{\theta(k)})$ such 
that $g_k=(\jmath_k)^{-1}\circ f_{\theta(k)}|_{J_k}$ and 
$\xi_{lk}=(\jmath_l)^{-1}\circ\rho_{\theta(l)\theta(k)}\circ\jmath_k$, where 
the transition map $\jmath_0\in Tran(V_0,U_0)$ is required to be 
$1_{G_{U_o}}$ in $G_{U_o}$ (note that $V_0=U_0=U_o$ as required). With 
these understood,  an equivalence relation on these systems can be similarly 
defined. 

\vspace{1.5mm}

\noindent{\bf Definition 3.1.1:}
{\em The based loop space of $(X,\underline{o})$, denoted by 
$\Omega(X,\underline{o})$, is the set of equivalence classes of systems 
from $(S^1,\ast)$ into $(X,\underline{o})$, as described in
the previous paragraph.
}

\hfill $\Box$

The based loop space $\Omega(X,\underline{o})$ carries a natural 
``compact-open'' topology which we shall describe next.

A subbase of the topology on $\Omega(X,\underline{o})$ is constructed as
follows. Let $\sigma=(\{I_i\},\{U_i\},\{f_{i,0}\},\{\rho_{ji,0}\})$ be any
based system from $(S^1,\ast)$ to $(X,\underline{o})$. Given $K=\{K_i\}$ where
each $K_i\subset I_i$ is compact, $V=\{V_i\}$ where each $V_i\subset
\widehat{U_i}$ is open such that $f_{i,0}(K_i)\subset V_i$, and $A=\{A_{ji}\}$
where each $A_{ji}$ is an open neighborhood of $\rho_{ji,0}$ in 
$Tran(U_i,U_j)$, we define a subset $\O(\sigma,K,V,A)$ of 
$\Omega(X,\underline{o})$ by setting
$$
\O(\sigma,K,V,A):=\{\tilde{f}|\tilde{f}=[(\{I_i\},\{U_i\},\{f_{i}\},
\{\rho_{ji}\})], f_i(K_i)\subset V_i, \rho_{ji}\in A_{ji}\}.
\leqno (3.1.1)
$$
The topology on $\Omega(X,\underline{o})$ is the one generated by the set of
all $\O(\sigma,K,V,A)$.

\vspace{1.5mm}

We shall next give a canonical description of a neighborhood of a based loop. 
Let $\tilde{f}_0$ be a based loop in $\Omega(X,\underline{o})$ and 
$\sigma_0=(\{I_i\},\{U_i\},\{f_{i,0}\},\{\rho_{ji,0}\})$, $i=0,1,\dots,n$,  
be a system representing it. For any $\tilde{f}\in\O(\sigma_0,K,V,A)$ with 
some $K=\{K_i\}$, $V=\{V_i\}$, and $A=\{A_{ji}\}$, there is a representing 
system $\sigma=(\{I_i\},\{U_i\},\{f_{i}\},\{\rho_{ji}\})$ satisfying 
$f_i(K_i)\subset V_i$ and $\rho_{ji}\in A_{ji}$. Observe that for each 
$\rho_{ji}$, there are $g_j\in G_{U_j}$, $h_i\in G_{U_i}$ such that 
$\rho_{ji}=g^{-1}_j\circ\rho_{ji,o}\circ h_i$. Set $W$ to be the connected 
component of $U_i\cap U_j$ over whose $G$-structure the transition map 
$\rho_{ji,0}$ is defined. Then both $g_j$ and $h_i$ are in $G_W$ because 
$G_W$ is open in both $G_{U_j}$ and $G_{U_i}$. Hence there are neighborhoods 
$O_j$ of $1_{G_{U_j}}$ in $G_{U_j}$ and $\delta_j\in O_j$, $j=1,\dots, n$, 
such that $\delta_1\circ\rho_{10}=\rho_{10,0}$ and 
$\delta_j\circ\rho_{ji}\circ\delta_i^{-1}=\rho_{ji,0}$ for $1<j\leq n$. 
Set $f_j^\prime=\delta_j\circ f_j$. Then $\tilde{f}$ has a representing 
system $(\{I_i\},\{U_i\},\{f_{i}^\prime\},\{\rho_{ji}^\prime\})$ where 
$f_j^\prime(K_j)\subset O_j\cdot V_j$ and $\rho_{ji}^\prime=\rho_{ji,0}$ for 
$j=1,\cdots,n$, and $f_0^\prime(K_0)\subset V_0$ and $\rho_{0n}^\prime
=\delta\circ\rho_{0n,0}$ for some $\delta$ in a neighborhood $A_o$ of 
$1_{G_{U_o}}$. Since such a representing system is clearly unique, and each 
open set in $\widehat{U_i}$ can be written as a union of open sets of form 
$O\cdot V$ where $O$ is an open neighborhood of $1_{G_{U_i}}$ and 
$V$ is an open set in $\widehat{U_i}$. We actually have proved the following 

\vspace{1.5mm}

\noindent{\bf Lemma 3.1.2:}
{\em Let $\tilde{f}_0$ be a based loop which is represented by a system 
$\sigma_0=(\{I_i\},\{U_i\},\{f_{i,0}\},\{\rho_{ji,0}\})$, $i=0,1,\dots,n$. 
Then a neighborhood of $\tilde{f}_0$ in $\Omega(X,\underline{o})$ can be 
identified with a set of systems $(\{f_{i}\},\{\rho_{ji}\})$ defined over 
$(\{I_i\},\{U_i\})$ such that each $f_i$ is in a neighborhood of $f_{i,0}$ 
in the space of continuous maps $C^0(I_i,\widehat{U_i})$, which is given 
with the compact-open topology, and each $\rho_{ji}=\rho_{ji,0}$ for 
$j=1,\cdots,n$, and $\rho_{0n}=\delta\circ\rho_{0n,0}$ 
for some $\delta$ in a neighborhood of $1_{G_{U_o}}$ in $G_{U_o}$. 
}

\hfill $\Box$

We shall call a neighborhood of a based loop as described in the above lemma 
{\it a canonical neighborhood} associated to a given representing system. 

\vspace{1.5mm}

There is a special element in $\Omega(X,\underline{o})$, which is the 
equivalence class of $(\{I_i\},\{U_i\},\{f_i\},\{\rho_{ji}\})$ where for each
index $i$, $U_i=U_o$ and $f_i$ is the constant map into $\hat{o}$, and 
$\rho_{ji}=1_{G_{U_o}}\in G_{U_o}=Tran(U_o,U_o)$. We denote this element by 
$\tilde{o}$ and fix it as the base point of $\Omega(X,\underline{o})$. 

\vspace{1.5mm}

We shall next construct continuous maps $\#:\Omega(X,\underline{o})
\times\Omega(X,\underline{o})\rightarrow \Omega(X,\underline{o})$ and
$\nu:\Omega(X,\underline{o})\rightarrow \Omega(X,\underline{o})$ so that
$((\Omega(X,\underline{o}),\tilde{o}),\#,\nu)$ becomes an $H$-group (cf.
\cite{Sw}).

The map $\#:\Omega(X,\underline{o})\times\Omega(X,\underline{o})\rightarrow 
\Omega(X,\underline{o})$, written $(\tilde{\gamma}_1,\tilde{\gamma}_2)\mapsto
\tilde{\gamma}_1\#\tilde{\gamma}_2$, is defined as follows. 
Let $\tilde{\gamma}_1$ and $\tilde{\gamma}_2$ be represented by systems
$(\{I_\alpha,I_0\},\{U_\alpha,U_o\},\{\gamma_{\alpha,1},\gamma_{0,1}\},
\{\xi_{\beta\alpha},\xi_{0\alpha}\})$ and $(\{J_a,J_0\},\{V_a,U_o\},
\{\gamma_{a,2},\gamma_{0,2}\},\{\xi_{ba},\xi_{0a}\})$ respectively. We 
decompose $I_0$ into $I_{0,-}\cup I_{0,+}$, and $J_0$ into $J_{0,-}\cup
J_{0,+}$ along the base point $\ast\in S^1$, and set $H=I_{0,-}\cup J_{0,+}$
and $H_0=J_{0,-}\cup I_{0,+}$. Then the system
$(\{I_\alpha,J_a,H,H_0\},\{U_\alpha,V_a,U_o,U_o\},\{\gamma_{\alpha,1},
\gamma_{a,2},\gamma,\gamma_0\}, \{\xi_{\beta\alpha},\xi_{ba},\xi_{0\alpha},
\xi_{0a}\})$ defines a based morphism from $([0,2]/\{0\sim 2\},\ast)$ into 
$(X,\underline{o})$, where $\gamma:H\rightarrow \widehat{U}_o$ and 
$\gamma_0: H_0\rightarrow\widehat{U}_o$ are defined piecewisely from
$\gamma_{0,1}$ and $\gamma_{0,2}$ in the obvious way. We pre-compose this based
morphism by the homeomorphism $\times_2:(S^1,\ast)\rightarrow 
([0,2]/\{0\sim 2\},\ast)$ where $\times_2$ is obtained from multiplying 
by $2$, and define $\tilde{\gamma}_1\#\tilde{\gamma}_2$ to be the resulting
based morphism in $\Omega(X,\underline{o})$. We leave the verification that 
$\#$ is well-defined to the readers. (The key point here is that in
the definition of isomorphism relation between based systems, the
automorphism $\delta_o$ of $(\widehat{U}_o,G_{U_o},\pi_{U_o})$ is
required to be $1_{G_{U_o}}\in G_{U_o}$ in Definition 2.2.1.) 
The map $\#:\Omega(X,\underline{o})\times
\Omega(X,\underline{o})\rightarrow \Omega(X,\underline{o})$ is clearly 
continuous. As for the map $\nu:\Omega(X,\underline{o})\rightarrow
\Omega(X,\underline{o})$, we define for any $\tilde{\gamma}\in
\Omega(X,\underline{o})$ the inverse of $\tilde{\gamma}$, written as 
$\nu(\tilde{\gamma})$, to be the equivalence class of the systems obtained 
by pre-composing the representatives of $\tilde{\gamma}$ with the 
self-homeomorphism $t\mapsto 1-t$ of $(S^1,\ast)$. The such defined map 
$\nu:\Omega(X,\underline{o})\rightarrow\Omega(X,\underline{o})$ is clearly 
continuous.

\vspace{1.5mm}

\noindent{\bf Lemma 3.1.3:}
{\em The based topological space $(\Omega(X,\underline{o}),\tilde{o})$ is an
$H$-group with homotopy associative multiplication $\#$ and homopoty inverse
$\nu$.
}

\vspace{2mm}

\noindent{\bf Proof:}
We refer to \cite{Sw} for the definition of $H$-group. Here we only sketch
a proof that $\nu$ is a homotopy inverse, i.e., both maps $\#\circ (\nu,1),
\#\circ (1,\nu):\Omega(X,\underline{o})\rightarrow\Omega(X,\underline{o})$
are homopotic to the constant map onto $\tilde{o}\in\Omega(X,\underline{o})$.
The proof of homotopy associativity of $\#$ is completely parallel.
 
Recall that we identify $S^1$ with $[0,1]/\{0\sim 1\}$. Now in this 
identification we decompose $[0,1]$ into $[0,\frac{1}{2}]\cup 
[\frac{1}{2},1]$. For any $s\in [0,1]$, set $I_s=([0,\frac{1}{2}(1-s)]\cup 
[\frac{1}{2}(1+s),1])/\{\frac{1}{2}(1-s)\sim\frac{1}{2}(1+s)\}$, and 
$S^1_s=I_s/\{0\sim 1\}$ with the base-point $\ast=0$.
It is easy to see that for any $\tilde{\gamma}\in\Omega(X,\underline{o})$ and
$s\in [0,1]$, the restriction of any of $\tilde{\gamma}\#\nu(\tilde{\gamma})$ 
or $\nu(\tilde{\gamma})\#\tilde{\gamma}$ to $I_s$ defines a based morphism 
from $(S^1_s,\ast)$ into $(X,\underline{o})$. Let $\eta_s:(S^1,\ast)
\rightarrow (S^1_s,\ast)$ be the family of homeomorphisms defined by
$$
\eta_s(t):=\left\{\begin{array}{cc}
(1-s)t & 0\leq t\leq\frac{1}{2}\\
(1-s)(t-1)+1 & \frac{1}{2}\leq t\leq 1.
\end{array} \right. \leqno (3.1.2)
$$
Then the family of compositions $(\tilde{\gamma}\#\nu(\tilde{\gamma}))
\circ\eta_s$ or $(\nu(\tilde{\gamma})\#\tilde{\gamma})\circ\eta_s$, 
$s\in [0,1]$, provides the required homotopy from $\#\circ (1,\nu)$ or 
$\#\circ (\nu,1)$ to the constant map onto 
$\tilde{o}\in\Omega(X,\underline{o})$. 

\hfill $\Box$

Thus we have constructed, for each based orbispace $(X,\underline{o})$,
a natural based topological space $(\Omega(X,\underline{o}),\tilde{o})$,
which is in fact an H-group with homotopy
associative multiplication $\#$ and homotopy inverse $\nu$.  
Moreover, for each based morphism 
$\tilde{f}:(X,\underline{o})\rightarrow (X^\prime,\underline{o^\prime})$, 
there is a corresponding map 
$\Omega(\tilde{f}):(\Omega(X,\underline{o}),\tilde{o})\rightarrow 
(\Omega(X^\prime,\underline{o^\prime}),\tilde{o}^\prime)$ defined by
$\Omega(\tilde{f})(\tilde{\gamma})=\tilde{f}\circ\tilde{\gamma}$. The map
$\Omega(\tilde{f})$ clearly satisfies the equations 
$\Omega(\tilde{f})(\tilde{\gamma}_1\#\tilde{\gamma}_2)=
\Omega(\tilde{f})(\tilde{\gamma}_1)\#\Omega(\tilde{f})(\tilde{\gamma}_2)$ and
$\Omega(\tilde{f})(\nu(\tilde{\gamma}))=
\nu(\Omega(\tilde{f})(\tilde{\gamma}))$. It is also continuous, and we have
actually obtained a functor $\Omega$ from the category of based
orbispaces to the category of H-groups. 

The continuity of
$\Omega(\tilde{f}):\Omega(X,\underline{o})\rightarrow 
\Omega(X^\prime,\underline{o^\prime})$ can be seen
as follows. Let $\tilde{\gamma}^\prime=\tilde{f}\circ\tilde{\gamma}$
be a based loop in $(X^\prime,\underline{o^\prime})$ which is the
image of a based loop $\tilde{\gamma}\in \Omega(X,\underline{o})$
under the map $\Omega(\tilde{f})$. Suppose we are given a neighborhood of 
$\tilde{\gamma}^\prime$, which is assumed without loss of generality
to be a canonical neighborhood associated to a representing system
$$
\sigma^\prime=(\{I_i\},\{U_i^\prime\},\{\gamma_i^\prime\},
\{\xi_{ji}^\prime\}), \hspace{2mm} i=0,1,\cdots, n, \leqno (3.1.3)
$$
as described in Lemma 3.1.2. 
We may further assume that the system $\sigma^\prime$ given in
$(3.1.3)$ is the composition of a representing system of 
$\tilde{\gamma}$ with a representing system of $\tilde{f}$. More
precisely, there are systems $\tau=(\{U_\alpha\},
\{U_\alpha^\prime\},\{f_\alpha\},\{\rho_{\beta\alpha}\})$ and
$\sigma=(\{I_i\},\{U_i\},\{\gamma_i\},\{\xi_{ji}\})$, and there is a
mapping of indexes $i\mapsto \alpha=\theta(i)$ with
$U_i=U_{\theta(i)}$ such that
$$
U_i^\prime=U^\prime_{\theta(i)}, \hspace{2mm} 
\gamma_i^\prime=f_{\theta(i)}\circ \gamma_i, \hspace{2mm}
\xi_{ji}^\prime=\rho_{\theta(j)\theta(i)}(\xi_{ji}). 
\leqno (3.1.4)
$$
Suppose the canonical neighborhood of $\tilde{\gamma}^\prime$
associated to the representing system $\sigma^\prime$ of $(3.1.3)$ is
given by the set of systems
$(\{\bar{\gamma}_i^\prime\},\{\bar{\xi}_{ji}^\prime\})$ 
defined over $(\{I_i\},\{U_i^\prime\})$, $i=0,1,\cdots,n$, satisfying 
$$ 
\bar{\gamma}_i^\prime(K_i)\subset V_i^\prime, \hspace{2mm}
\bar{\xi}_{ji}^\prime=\xi_{ji}^\prime \; \forall j=1,\cdots,n,
\hspace{2mm}
\bar{\xi}_{0n}^\prime=\delta^\prime\circ \xi_{0n}^\prime \hspace{2mm}
\mbox{ with } \delta^\prime\in O^\prime,
\leqno (3.1.5)
$$
where $K_i$ is a compact subset of $I_i$, $V_i^\prime$ is an open
subset of $\widehat{U^\prime_{\theta(i)}}$, and $O^\prime$ is a
neighborhood of $1_{G_{U_o^\prime}}$ in 
$G_{U_o^\prime}$. We set $V_i=f^{-1}_{\theta(i)}(V_i^\prime)\subset 
\widehat{U_{\theta(i)}}$ and $O=\rho_o^{-1}(O^\prime)$, which is a
neighborhood of $1_{G_{U_o}}$ in $G_{U_o}$. Then the
canonical neighborhood of $\tilde{\gamma}$ associated to the
representing system $\sigma$, which is given by the set of systems 
$(\{\bar{\gamma}_i\},\{\bar{\xi}_{ji}\})$ 
defined over $(\{I_i\},\{U_i\})$, $i=0,1,\cdots,n$, satisfying 
$$ 
\bar{\gamma}_i(K_i)\subset V_i, \hspace{2mm}
\bar{\xi}_{ji}=\xi_{ji} \; \forall j=1,\cdots,n,
\hspace{2mm}
\bar{\xi}_{0n}=\delta \circ \xi_{0n} \hspace{2mm}
\mbox{ with } \delta\in O, \leqno (3.1.6)
$$
is mapped into the canonical neighborhood of $\tilde{\gamma}^\prime$
described in $(3.1.5)$ under the map $\Omega(\tilde{f})$. Hence 
$\Omega(\tilde{f}):\Omega(X,\underline{o})\rightarrow 
\Omega(X^\prime,\underline{o^\prime})$ is continuous.

Finally, observe also that there is a natural continuous map $\Pi_X:
(\Omega(X,\underline{o}),\tilde{o})\rightarrow(\Omega(X_{top},o),o)$ defined 
by forgetting the orbispace structure on $X$, i.e., $\Pi_X(\tilde{\gamma})=
\gamma$ for each $\tilde{\gamma}\in\Omega(X,\underline{o})$, which clearly
preserves the H-group structure on the corresponding spaces.

\subsection{Homotopy groups}

\noindent{\bf Definition 3.2.1:}
{\em Let $(X,\underline{o})$ be a based orbispace. For any $k\geq 1$, 
the k-th homotopy group of $(X,\underline{o})$, denoted by 
$\pi_k(X,\underline{o})$, is defined to be the (k-1)-th homotopy group of the 
based loop space $(\Omega(X,\underline{o}),\tilde{o})$.
}

\hfill $\Box$

Because $(\Omega(X,\underline{o}),\tilde{o})$ has a structure of H-group, 
$\pi_1(X,\underline{o})$ is a group whose multiplication is induced by the 
homotopy associative multiplication $\#: \Omega(X,\underline{o})\times
\Omega(X,\underline{o})\rightarrow \Omega(X,\underline{o})$. Moreover, 
for $k\geq 2$, $\pi_k(X,\underline{o})$ is an abelian group. There is an action
of $\pi_1(X,\underline{o})$ on each $\pi_k(X,\underline{o})$, written as 
$C:\pi_1(X,\underline{o})\rightarrow Aut(\pi_k(X,\underline{o}))$, which is 
defined as follows: for any $a\in\pi_1(X,\underline{o})$ and 
$\phi\in\pi_k(X,\underline{o})$, where $a$ is represented by $\tilde{\gamma}\in
\Omega(X,\underline{o})$ and $\phi$ is represented by $\Phi:(S^{k-1},\ast)
\rightarrow (\Omega(X,\underline{o}),\tilde{o})$, we define $C(a)(\phi)$ to be
the class represented by the continuous map $(S^{k-1},\ast)\rightarrow 
(\Omega(X,\underline{o}),\tilde{o})$ given by $x\mapsto (\tilde{\gamma})^{-1}
\#\Phi(x)\#\tilde{\gamma}$. The action 
$C:\pi_1(X,\underline{o})\rightarrow Aut(\pi_k(X,\underline{o}))$ will be 
called the {\it canonical action} of $\pi_1(X,\underline{o})$ on 
$\pi_k(X,\underline{o})$.

In summary, we have defined for each $k\geq 1$ a functor $\pi_k$ from
the category of based orbispaces to the category of groups (abelian
groups for $k\geq 2$). When restricted to the sub-category of based
locally connected topological spaces, the functor $\pi_k$ coincides 
with the functor of $k$-th homotopy group of a topological space, 
defined as the set of all homotopy classes of based continuous maps 
from $(S^k,\ast)$ into the based topological space. As a notational 
convention, for any based morphism $\tilde{f}:
(X,\underline{o})\rightarrow (X^\prime,\underline{o^\prime})$, we write
$f_{\#}:\pi_k(X,\underline{o})\rightarrow
\pi_k(X^\prime,\underline{o^\prime})$ for the corresponding
homomorphism between homotopy groups. 

Each based morphism $\tilde{f}:(S^k,\ast)\rightarrow
(X,\underline{o})$ gives rise to an element $f_{\#}(1)\in
\pi_k(X,\underline{o})$ where $1\in \pi_k(S^k,\ast)$ is the class of
the identity map $S^k\rightarrow S^k$. We claim that the converse is
true when the orbispace $X$ is \'{e}tale. More concretely, for any
element $u$ of $\pi_k(X,\underline{o})$ with $X$ being \'{e}tale,
there is a based morphism $\tilde{f}:(S^k,\ast)\rightarrow
(X,\underline{o})$ such that $f_{\#}(1)=u$. (An orbispace $(X,\U)$ is
called \'{e}tale if $G_U$ is discrete for every $U\in\U$.) The claim
can be seen as follows. Let $u:(S^{k-1},\ast)\rightarrow
(\Omega(X,\underline{o}),\tilde{o})$ be a continuous map. We subdivide
$S^{k-1}$ into a union of finitely many closed simplexes $\cup_{a\in
A}K_a$, which is fine enough so that the image of each $K_a$ under $u$
lies in a canonical neighborhood of some based loop as described in
Lemma 3.1.2. Hence there are parametrized representing systems of
based loops $\sigma_a$, $a\in A$,
$$
\sigma_a=(\{I_i\times K_a\},\{U_i^a\},\{\gamma_i^a\},\{\xi_{ji}^a\}),
\hspace{2mm} i=0,1,\cdots,n, \leqno (3.2.1)
$$ 
where each $\gamma_i^a:I_i\times K_a \rightarrow \widehat{U_i^a}$ is
continuous, each $\xi_{ji}^a$ is a transition map in
$Tran(U_i^a,U_j^a)$ which is constant in $x\in K_a$ (here we have used the
assumption that the orbispace $X$ is \'{e}tale), and the equations 
$\gamma_j^a=\xi_{ji}^a\circ\gamma_i^a$ are satisfied. Here the number
of intervals $I_i$ can be chosen to be independent of $K_a$ by passing
to refinement. When restricted to the common face of $K_a$ and $K_b$,
the parametrized systems $\sigma_a$ and $\sigma_b$ represent the same
map $u|_{K_a\cap K_b}$, hence there are transition maps
$\eta_{ba}^i\in Tran(U_i^a,U_i^b)$ such that 
$$
\xi_{ji}^b=\eta_{ba}^j\circ \xi_{ji}^a\circ (\eta_{ba}^i)^{-1},
\hspace{2mm} \gamma_i^b=\eta_{ba}^i\circ \gamma_i^a \hspace{2mm}
\mbox{ on } I_i\times (K_a\cap K_b).
\leqno (3.2.2)
$$
Here again each $\eta_{ba}^i$ is constant in $x\in K_a\cap K_b$ by the
assumption that the orbispace $X$ is \'{e}tale. 
Now we regard $S^k$ as the suspension $SS^{k-1}$ of $S^{k-1}$, and
define a system of maps 
$$
\tau=(\{V_i^a\},\{U_i^a\},\{f_i^a\},\{\rho_{(j,b)(i,a)}\}),
\hspace{2mm} i=0,1,\cdots,n, \; a\in A \leqno (3.2.3)
$$
where each $V_i^a=I_i\times O_a$ for some small regular neighborhood $O_a$
of the closed simplex $K_a$, $f_i^a:V_i^a\rightarrow \widehat{U_i^a}$
is defined from a suitable extension of $\gamma_i^a$, and
$\rho_{(j,b)(i,a)}=\eta_{ba}^j\circ\xi_{ji}^a=\xi_{ji}^b\circ\eta_{ba}^i$ 
is a transition map in $Tran(U_i^a,U_j^b)$. The extension $f_i^a$ 
of $\gamma_i^a$ is
chosen in such a way that $f_i^b=\eta_{ba}^i\circ \gamma_i^a$. It is
easy to verify that the system $\tau$ as defined in $(3.2.3)$ 
satisfies the compatibility
conditions $(2.2.1 a-b)$ so that it determines a based morphism
$\tilde{f}: (S^k,\ast)\rightarrow (X,\underline{o})$. By the nature of
construction, we have $f_{\#}(1)=[u]$ where $1\in \pi_k(S^k,\ast)$ is
the class of the identity map $S^k\rightarrow S^k$. This concludes the
proof of the claim.

\vspace{1.5mm}

In what follows, we shall formulate and establish the homotopy invariance 
of groups $\pi_k(X,\underline{o})$, and investigate how they behave under 
the change of base-point structures.

\vspace{1.5mm}

Let $(X,\underline{o})$, $(X^\prime,\underline{o^\prime})$ be two based 
orbispaces. The space of based morphisms between $(X,\underline{o})$ and 
$(X^\prime,\underline{o^\prime})$, denoted by $Mor\{(X,\underline{o}),
(X^\prime,\underline{o^\prime})\}$, carries a natural ``compact-open''
topology.  Let $\sigma=(\{U_\alpha\},\{U_\alpha^\prime\},\{f_{\alpha,0}\},
\{\rho_{\beta\alpha,0}\})$ be a based system from 
$(X,\underline{o})$ to $(X^\prime,\underline{o^\prime})$. 
Let $K=\{K_\alpha\}$ be a finite set where each 
$K_\alpha\subset\widehat{U_\alpha}$ is compact, $V=\{V_\alpha\}$ be a finite
set where each $V_\alpha\subset\widehat{U_\alpha^\prime}$ is open such that
$f_{\alpha,0}(K_\alpha)\subset V_\alpha$, and $L=\{L_{\beta\alpha}\}$ and 
$O=\{O_{\beta\alpha}\}$ be finite sets where each $L_{\beta\alpha}$ is a 
compact subset in $Tran(U_\alpha, U_\beta)$ and each $O_{\beta\alpha}$ is an 
open subset of $Tran(U_\alpha^\prime,U_\beta^\prime)$ such that 
$\rho_{\beta\alpha,0}(L_{\beta\alpha})\subset O_{\beta\alpha}$. We introduce 
a subset $\O(\sigma,K,V,L,O)$ of 
$Mor\{(X,\underline{o}),(X^\prime,\underline{o^\prime})\}$, which is the set 
of all based morphisms represented by systems 
$(\{U_\alpha\},\{U_\alpha^\prime\},\{f_{\alpha}\},\{\rho_{\beta\alpha}\})$ 
satisfying $f(K_\alpha)\subset V_\alpha$ and 
$\rho_{\beta\alpha}(L_{\beta\alpha})\subset O_{\beta\alpha}$ for all 
$K_\alpha\in K$, $V_\alpha\in V$, $L_{\beta\alpha}\in L$ and 
$O_{\beta\alpha}\in O$. The ``compact-open'' topology on 
$Mor\{(X,\underline{o}),(X^\prime,\underline{o^\prime})\}$ is the one 
generated by all these subsets $\O(\sigma,K,V,L,O)$.

\vspace{1.5mm}

\noindent{\bf Definition 3.2.2:}
{\em Two based morphisms $\tilde{f}_1,\tilde{f}_2: (X,\underline{o})
\rightarrow (X^\prime,\underline{o^\prime})$ are said to be 
homotopic if $\tilde{f}_1,\tilde{f}_2$ are path-connected in the space of all based 
morphisms $Mor\{(X,\underline{o}),(X^\prime,\underline{o^\prime})\}$ 
equipped with the ``compact-open'' topology. A based morphism 
$\tilde{f}:(X,\underline{o})\rightarrow (X^\prime,\underline{o^\prime})$ is 
said to be a homotopy equivalence if there is a based morphism 
$\tilde{g}:(X^\prime,\underline{o^\prime})\rightarrow (X,\underline{o})$ 
such that both of $\tilde{g}\circ\tilde{f}$ and $\tilde{f}\circ\tilde{g}$ 
are homotopic to the identity morphism. 
}

\hfill $\Box$

\noindent{\bf Theorem 3.2.3:}
{\em Let $\tilde{f}_1,\tilde{f}_2: (X,\underline{o})\rightarrow 
(X^\prime,\underline{o^\prime})$ be two based morphisms. Then $(f_1)_{\#}=
(f_2)_{\#}:\pi_k(X,\underline{o})\rightarrow
\pi_k(X^\prime,\underline{o^\prime})$ for all $k\geq 1$ if 
$\tilde{f}_1,\tilde{f}_2$ are homotopic. As a consequence, 
$f_{\#}:\pi_k(X,\underline{o})\rightarrow\pi_k(X^\prime,\underline{o^\prime})$
is an isomorphism for any $k\geq 1$ if $\tilde{f}:(X,\underline{o})\rightarrow 
(X^\prime,\underline{o^\prime})$ is a homotopy equivalence.
}

\vspace{2mm}

\noindent{\bf Proof:}
Let $\tilde{f}_s$, $s\in [1,2]$, be a path connecting $\tilde{f}_1$ and
$\tilde{f}_2$ in $Mor\{(X,\underline{o}),(X^\prime,\underline{o^\prime})\}$. 
We consider the map $F:(\Omega(X,\underline{o}),\tilde{o})\times [1,2]
\rightarrow (\Omega(X^\prime,\underline{o^\prime}),\tilde{o}^\prime)$ defined 
by $F(\tilde{\gamma},s)=\Omega(\tilde{f}_s)(\tilde{\gamma})$, where 
$\Omega(\tilde{f}_s):(\Omega(X,\underline{o}),\tilde{o})\rightarrow 
(\Omega(X^\prime,\underline{o^\prime}),\tilde{o}^\prime)$ is given by
$\Omega(\tilde{f_s})(\tilde{\gamma})=\tilde{f_s}\circ\tilde{\gamma}$.
We need to show that for any continuous map $u:S^k\rightarrow
\Omega(X,\underline{o})$, the composition $\widehat{F}:S^k\times [1,2]
\rightarrow\Omega(X^\prime,\underline{o^\prime})$ defined by 
$\widehat{F}(x,s):=F(u(x),s)$ is continuous.

For any fixed $(x_0,s_0)\in S^k\times [1,2]$, set $u(x_0)=\tilde{\gamma}_0
\in\Omega(X,\underline{o})$. We will show that for any neighborhood of 
$F(\tilde{\gamma}_0,s_0)$ in $\Omega(X^\prime,\underline{o^\prime})$ of form 
$(3.1.1)$, $\O(\sigma^\prime,K^\prime,V^\prime,A^\prime)$, there is a 
neighborhood $\O(\sigma,K^\prime,V,A)$ of $\tilde{\gamma}_0$, such that for
a given compact subset ${\cal{C}}\subset\O(\sigma,K^\prime,V,A)$, there is
an interval $I$ containing $s_0$, such that for any $(\tilde{\gamma},s)\in 
{\cal{C}}\times I$, we have 
$F(\tilde{\gamma},s)\in\O(\sigma^\prime,K^\prime,V^\prime,A^\prime)$. 
This is seen as follows. 
Let $\sigma=(\{I_i\},\{U_i\},\{\gamma_{i,0}\},\{\xi_{ji,0}\})$ be a
system representing $\tilde{\gamma}_0$, $\tau=(\{U_\alpha\},
\{U_\alpha^\prime\},\{f_{\alpha,s_0}\},\{\rho_{\beta\alpha,s_0}\})$ a system
representing $\tilde{f}_{s_0}$ such that the composition 
$\sigma^\prime=\tau\circ\sigma$, written 
$(\{I_i\},\{U_i^\prime\},\{f_{i,s_0}\circ\gamma_{i,0}\},
\{\rho_{ji,s_0}(\xi_{ji,0})\})$, represents $F(\tilde{\gamma}_0,s_0)$.
Here the index $i$ is running from $0$ to $n$.
Given $K^\prime=\{K_i^\prime\}$ where each $K_i^\prime\subset I_i$ is compact,
$V^\prime=\{V_i^\prime\}$ where each $V_i^\prime\subset\widehat{U_i^\prime}$ is
open such that $f_{i,s_0}\circ\gamma_{i,0}(K_i^\prime)\subset V_i^\prime$, and
$A^\prime=\{A^\prime_{ji}\}$ where each $A^\prime_{ji}$ is an open
neighborhood of $\rho_{ji,s_0}(\xi_{ji,0})$ in $Tran(U_i^\prime,U_j^\prime)$.
For each $1\leq i\leq n$, since $\gamma_{i,0}(K_i^\prime)$ is compact in 
$f_{i,s_0}^{-1}(V_i^\prime)$, it is easily seen that there is an open set
$V_i\subset \widehat{U_i}$, a neighborhood $O_i$ of $1_{G_{U_i}}$, such that
$O_i\cdot V_i\subset f_{i,s_0}^{-1}(V_i^\prime)$ and 
$\gamma_{i,0}(K_i^\prime)\subset V_i$. We set $V_0=f^{-1}_{0,s_0}(V_0^\prime)$
and $V=\{V_i\}$. We choose $O_n$ small enough so that there is a neighborhood
$A_{0n}$ of $\xi_{0n,0}$ such that $A_{0n}\circ O_n^{-1}
\subset\rho_{0n,s_0}^{-1}(A_{0n}^\prime)$. For other $A_{ji}$, we 
set $A_{10}=O_1^{-1}\circ\xi_{10,0}$, 
$A_{ji}=O_j^{-1}\circ\xi_{ji,0}\circ O_i$ for $2\leq j\leq n$, and 
set $A=\{A_{ji}\}$. We consider the neighborhood $\O(\sigma,K^\prime,V,A)$ of 
$\tilde{\gamma}_0$ in $\Omega(X,\underline{o})$. 

Set $K=\{K_i\}$ and $L=\{L_{ji}\}$ where $K_i:=\{\gamma_i(t)\}$ and 
$L_{ji}:=\{\xi_{ji}\}$, in which $t\in K^\prime_i$ and 
$(\{\gamma_i\},\{\xi_{ji}\})$ is in a canonical neighborhood of 
$\tilde{\gamma}_0$ associated to the system $\sigma$ (cf. Lemma 3.1.2) such
that $[(\{\gamma_i\},\{\xi_{ji}\})]\in {\cal{C}}$.
Then it is easily seen that each $K_i$ is a compact subset of 
$\widehat{U_i}$ and each $L_{ji}$ is a compact subset of $Tran(U_i,U_j)$.
Set $O_{ji}=A_{ji}^\prime$. Then by the nature of construction, we have 
$f_{i,s_0}(K_i)\subset V_i^\prime$ and $\rho_{ji,s_0}(L_{ji})\subset O_{ji}$. 
Let $I$ be a neighborhood of $s_0$ such that for all $s\in I$, the based 
morphisms $\tilde{f}_s$ are in $\O(\tau,K,V^\prime,L,O)$. 
Then it follows from the nature of construction that for any 
$(\tilde{\gamma},s)\in {\cal{C}}\times I$, we have $F(\tilde{\gamma},s)
\in\O(\sigma^\prime,K^\prime,V^\prime,A^\prime)$. 

By taking the compact subset $\cal{C}$ to be the image of a compact
neighborhood of $x_0$ in $S^k$ under the continuous map $u$, we prove
that $\widehat{F}(x,s)$ is continuous at $(x_0,s_0)$ for any
$(x_0,s_0)\in S^k\times [1,2]$. Hence the theorem.

\hfill $\Box$

Let $\underline{o}_1=(o_1,U_{o,1},\hat{o_1})$, 
$\underline{o}_2=(o_2,U_{o,2},\hat{o_2})$ be two base-point structures of 
orbispace $(X,\U)$. Let $\{I_i| i=0, 1, \cdots, n\}$ be a cover of $[0,1]$
by open intervals such that $0\in I_0$, $1\in I_n$ and 
$I_i\cap I_j\neq\emptyset$ only if $j=i+1$, and $\{U_i| i=0, 1, \cdots, n\}$ 
be a collection of basic open sets of $X$ where $U_0=U_{o,1}$ and
$U_n=U_{o,2}$. We shall consider systems such as  
$(\{I_i\},\{U_i\},\{f_i\},\{\rho_{ji}\})$, where each $f_i:I_i\rightarrow 
\widehat{U_i}$ is a continuous map and each $\rho_{ji}$ is in 
$Tran(U_i, U_j)$ where $j=i+1$. Furthermore, we 
require that on each $I_i\cap I_j$, equation $f_j=\rho_{ji}\circ f_i$ holds, 
and $f_0(0)=\hat{o_1}$ in $\widehat{U_{o,1}}$ and $f_n(1)=\hat{o_2}$ in
$\widehat{U_{o,2}}$. An isomorphism relation on these 
systems is defined as in Definition 2.2.1, but we require that the 
automorphisms $\delta_{o,1}$ and $\delta_{o,2}$ of the $G$-structure 
$(\widehat{U_{o,1}},G_{U_{o,1}},\pi_{U_{o,1}})$ and 
$(\widehat{U_{o,2}},G_{U_{o,2}},\pi_{U_{o,2}})$ be $1_{G_{U_{o,1}}}$ and 
$1_{G_{U_{o,2}}}$ respectively. A system 
$(\{J_k\},\{V_k\},\{g_k\},\{\xi_{lk}\})$, where $k$ is running from $0$ to 
$m$, is said to be a refinement of $(\{I_i\},\{U_i\},\{f_i\},\{\rho_{ji}\})$ 
if (1) $\{J_k\}$ is a refinement of $\{I_i\}$ given by a mapping 
$\theta:\{0,\cdots,m\}\rightarrow \{0,\cdots,n\}$ such that $\theta(0)=0$ and
$\theta(m)=n$, and (2) there are transition maps 
$\jmath_k\in Tran(V_k, U_{\theta(k)})$ such that 
$g_k=(\jmath_k)^{-1}\circ f_{\theta(k)}|_{J_k}$ and 
$\xi_{lk}=(\jmath_l)^{-1}\circ\rho_{\theta(l)\theta(k)}\circ\jmath_k$, where 
the transition maps $\jmath_0\in Tran(V_0,U_0)$ and 
$\jmath_m\in Tran(V_m,U_n)$ are required to be $1_{G_{U_{o,1}}}$ and 
$1_{G_{U_{o,2}}}$ respectively (note that $V_0=U_0=U_{o,1}$ and $V_m=U_n=
U_{o,2}$ as required). With these understood, an equivalence relation on 
these systems can be defined as in Definition 3.1.1.

\vspace{1.5mm}

\noindent{\bf Definition 3.2.4:}
{\em A based morphism from $[0,1]$ to $(X,\underline{o_1},\underline{o_2})$
is defined to be an equivalence class of the systems as described in the 
previous paragraph.
}

\hfill $\Box$

A based morphism from $[0,1]$ to $(X,\underline{o_1},\underline{o_2})$ will 
be called a {\it based path} from $\underline{o_1}$ to $\underline{o_2}$. 
The space of all based paths from $\underline{o_1}$ to $\underline{o_2}$ is
denoted by $P(X,\underline{o_1},\underline{o_2})$. Clearly 
$P(X,\underline{o},\underline{o})$ is just the based loop space 
$\Omega(X,\underline{o})$. The homotopy inverse $\nu$ and the homotopy 
associative multiplication $\#$ on a based loop space can be defined more
generally on the spaces of based paths. To be more concrete, there is a 
continuous map $\nu:P(X,\underline{o_1},\underline{o_2})\rightarrow
P(X,\underline{o_2},\underline{o_1})$, where $\nu(\tilde{\gamma})$ is 
obtained by pre-composing $\tilde{\gamma}$ with the self-homeomorphism of 
$[0,1]$ given by $t\mapsto 1-t$. The corresponding map 
$\#:P(X,\underline{o_1},\underline{o_2})\times 
P(X,\underline{o_2},\underline{o_3})\rightarrow 
P(X,\underline{o_1},\underline{o_3})$, written as $(\tilde{\gamma}_1,
\tilde{\gamma}_2)\mapsto \tilde{\gamma}_1\#\tilde{\gamma}_2$, can be
similarly defined, generalizing the map 
$\#:\Omega(X,\underline{o})\times\Omega(X,\underline{o})
\rightarrow \Omega(X,\underline{o})$. Observe that for any based paths 
$\tilde{\gamma}_1,\tilde{\gamma}_2\in P(X,\underline{o_1},\underline{o_2})$, 
$\tilde{\gamma}_1\#\nu(\tilde{\gamma}_2)$ and 
$\nu(\tilde{\gamma}_1)\#\tilde{\gamma}_2$ are elements of 
$\Omega(X,\underline{o_1})$ and $\Omega(X,\underline{o_2})$ respectively.

\vspace{1.5mm}

The following theorem is self-evident.

\vspace{1.5mm}

\noindent{\bf Theorem 3.2.5:}
{\em For any based morphism $\tilde{\gamma}:[0,1]\rightarrow 
(X,\underline{o_1},\underline{o_2})$, there is an isomorphism, written as
$\tilde{\gamma}_\ast:\pi_k(X,\underline{o_2})\rightarrow 
\pi_k(X,\underline{o_1})$, whose inverse $(\gamma_\ast)^{-1}$ is given by 
$\nu(\tilde{\gamma})_\ast:\pi_k(X,\underline{o_1})\rightarrow 
\pi_k(X,\underline{o_2})$. Moreover, for any 
$\tilde{\gamma}_1,\tilde{\gamma}_2$, we have $(\tilde{\gamma}_2)^{-1}_\ast\circ
(\tilde{\gamma}_1)_\ast=C([\nu(\tilde{\gamma}_1)\#\tilde{\gamma}_2])$ as an
element of $Aut(\pi_k(X,\underline{o_2}))$, where 
$C:\pi_1(X,\underline{o_2})\rightarrow Aut(\pi_k(X,\underline{o_2}))$ is the 
canonical action of $\pi_1(X,\underline{o_2})$ on 
$\pi_k(X,\underline{o_2})$, and $[\nu(\tilde{\gamma}_1)\#\tilde{\gamma}_2]$ 
stands for the class of $\nu(\tilde{\gamma}_1)\#\tilde{\gamma}_2$ in 
$\pi_1(X,\underline{o_2})$.
}

\hfill $\Box$

In the special case when $\underline{o_1}=(o,U_{o,1},\hat{o_1})$ and 
$\underline{o_2}=(o,U_{o,2},\hat{o_2})$ for the same point $o\in X$ 
(i.e. $o_1=o_2=o$), there is always a based morphism $\tilde{\gamma}: 
[0,1]\rightarrow (X,\underline{o_1},\underline{o_2})$, defined as the 
equivalence class of the following system 
$(\{I_0,I_1\},\{U_{o,1},U_{o,2}\},\{\gamma_0,\gamma_1\},\{\rho_{10}\})$, where 
$I_0=[0,\frac{2}{3})$, $I_1=(\frac{1}{3},1]$, $\gamma_0(t)=\hat{o_1}$, 
$\forall t\in I_0$, $\gamma_1(t)=\hat{o_2}$, $\forall t\in I_1$, and 
$\rho_{10}\in Tran(U_{o,1},U_{o,2})$ is a transition map sending $\hat{o_1}$ 
to $\hat{o_2}$ ($\rho_{10}$ exists because 
$\pi_{U_{o,1}}(\hat{o_1})=\pi_{U_{o,2}}(\hat{o_2})=o$). As a consequence, 
the isomorphism classes of homotopy groups $\pi_k(X,\underline{o})$ where 
$\underline{o}=(o,U_o,\hat{o})$ depend only on the point $o\in X$. We will 
denote the corresponding abstract groups by $\pi_k(X,o)$. Note that for any 
morphism $\tilde{f}:X\rightarrow X^\prime$, if setting $o^\prime=f(o)$, there 
is a homomorphism $f_{\#}:\pi_k(X,o)\rightarrow\pi_k(X^\prime,o^\prime)$. 

\vspace{1.5mm}

Observe that the natural continuous map $\Pi_X:(\Omega(X,\underline{o}),
\tilde{o})\rightarrow (\Omega(X_{top},o),o)$ induces a homomorphism 
$\Pi_X:\pi_k(X,o)\rightarrow \pi_k(X_{top},o)$, which is natural in the 
following sense: the commutativity $f_\ast\circ\Pi_X=\Pi_{X^\prime}
\circ f_{\#}$ holds for any morphism $\tilde{f}:X\rightarrow X^\prime$, 
where $f_{\#}:\pi_k(X,o)\rightarrow \pi_k(X^\prime, f(o))$ and 
$f_\ast:\pi_k(X_{top},o)\rightarrow \pi_k(X^\prime_{top}, f(o))$.

\vspace{1.5mm}

An orbispace $X$ is called {\it connected} if $X_{top}$ is connected. Two 
points $o_1,o_2$ of $X$ are called {\it path-connected in X} if there is a 
morphism $\tilde{u}:[0,1]\rightarrow X$ such that $u(0)=o_1$ and $u(1)=o_2$. 
Clearly being path-connected is an equivalence relation amongst the points in 
$X$. We denote by $\pi_0(X)$ the set of all equivalence classes of 
path-connected points in $X$, and by $\pi_0(X,o)$ the corresponding based set 
with the base point being the equivalence class of point $o\in X$. We call an 
orbispace $X$ {\it path-connected} if $\pi_0(X)$ is a set of only one point.  

\vspace{1.5mm}

For a path-connected orbispace $X$, the isomorphism class of $\pi_k(X,o)$ does 
not depend on the base point $o\in X$. We will denote the corresponding 
abstract group by $\pi_k(X)$, and call it the {\it k-th homotopy group} of 
orbispace $X$. Note that for any morphism $\tilde{f}:X\rightarrow X^\prime$ 
between path-connected orbispaces, there is a homomorphism 
$f_{\#}:\pi_k(X)\rightarrow \pi_k(X^\prime)$ for any $k\geq 1$. 

\vspace{1.5mm}

Recall that for any two orbispaces $X, X^\prime$, the cartesian product 
$X\times X^\prime$ is defined, where the orbispace structure is taken as the
cartesian product of the orbispace structures of $X$ and $X^\prime$. Given
base-point structures $\underline{o}$, $\underline{o}^\prime$ of $X$ and 
$X^\prime$ respectively, there is a canonical base-point structure 
$\underline{o}\times\underline{o}^\prime$ of $X\times X^\prime$. Observe that
the set of systems (resp. based systems) from $Y$ to $X\times X^\prime$ is
the cartesian product of the set of systems (resp. based systems) from $Y$
to $X$ with the set of systems (resp. based systems) from $Y$ to $X^\prime$,
and passing to equivalence classes, the set of morphisms (resp. 
based morphisms) from $Y$ to $X\times X^\prime$ is the cartesian product of 
the set of morphisms (resp. based morphisms) from $Y$ to $X$ with the set of 
morphisms (resp. based morphisms) from $Y$ to $X^\prime$. This implies the
following

\vspace{1.5mm}

\noindent{\bf Theorem 3.2.6:}
{\em  $\pi_0(X\times X^\prime,(o,o^\prime))=\pi_0(X,o)\times\pi_0(X^\prime, 
o^\prime)$, and for any $k\geq 1$, $\pi_k(X\times X^\prime,\underline{o}
\times\underline{o}^\prime)=\pi_k(X,\underline{o})\times 
\pi_k(X^\prime,\underline{o}^\prime)$. 
}

\hfill $\Box$

Given an orbispace $X$, we can form the cartesian product $X\times [1,2]$
where the interval $[1,2]$ is a trivial orbispace. There are two orbispace
embeddings $\tilde{i}_j:X\rightarrow X\times [1,2]$ into the sub-orbispaces
$X\times\{j\}$, $j=1,2$, respectively.

\vspace{1.5mm}

\noindent{\bf Definition 3.2.7:}
{\em Two morphisms $\tilde{f}_1,\tilde{f}_2:X\rightarrow X^\prime$ are called 
homotopic if there is a morphism 
$\widetilde{F}:X\times [1,2]\rightarrow X^\prime$ such that 
$\tilde{f}_j=\widetilde{F}\circ\tilde{i}_j$ for $j=1,2$. A morphism 
$\tilde{f}:X\rightarrow X^\prime$ is called a homotopy equivalence if there is 
a morphism $\tilde{g}:X^\prime\rightarrow X$ such that both 
$\tilde{f}\circ\tilde{g}$ and $\tilde{g}\circ\tilde{f}$ are 
homotopic to the identity morphism.
}

\hfill $\Box$

The following theorem is self-evident.

\vspace{1.5mm}

\noindent{\bf Theorem 3.2.8:}
{\em Let $\tilde{f}:X\rightarrow X^\prime$ be a homotopy equivalence. Suppose 
$X$ is path-connected, then $X^\prime$ is also path-connected, 
and $\tilde{f}$ induces an isomorphism $f_{\#}:\pi_k(X)\rightarrow
\pi_k(X^\prime)$ for all $k\geq 0$.
}

\hfill $\Box$

We insert two remarks concerning the Hurewicz homomorphism. 

\vspace{1.5mm}

\noindent{\bf Remark 3.2.9 a:}\hspace{2mm}
The {\it Hurewicz homomorphism} $\pi_k(X,o)\rightarrow
H_k(X_{top},\Z)$ is defined as the composition $h\circ \Pi_X$ where
$\Pi_X: \pi_k(X,o)\rightarrow \pi_k(X_{top},o)$ is the natural
homomorphism induced by
$(\Omega(X,\underline{o}),\tilde{o})\rightarrow
(\Omega(X_{top},o),o)$, and $h:\pi_k(X_{top},o)\rightarrow
H_k(X_{top},\Z)$ is the ordinary Hurewicz homomorphism.

\hfill $\Box$

\noindent{\bf Remark 3.2.9 b:}\hspace{2mm} 
Suppose $X$ is an \'{e}tale orbispace. We have shown that for any
element $u\in \pi_k(X,\underline{o})$ there is a based morphism 
$\tilde{f}:(S^k,\ast)\rightarrow (X,\underline{o})$ such that
$f_{\#}(1)=u$ where $1\in\pi_k(S^k,\ast)$ is the class of the identity
map $S^k\rightarrow S^k$. In this case, the Hurewicz homomorphism
$\pi_k(X,o)\rightarrow H_k(X_{top},\Z)$ sends the class $u$
to the class $f_\ast([S^k])$ where $f:S^k\rightarrow X_{top}$ is the 
continuous map induced by the morphism $\tilde{f}$. It seems that 
there should be a singular homology theory of \'{e}tale orbispaces 
constructed using chains of morphisms from simplexes into
the orbispace, to which the above Hurewicz homomorphism can be lifted. 

\hfill $\Box$

We end this section with a theorem relating the fundamental group of a
complex of groups as defined in \cite{Ha2} with the $\pi_1$ of the
associated orbihedron which is canonically regarded as an orbispace in
the sense of this paper. 

We first recall the definition of fundamental group of a complex
of groups in \cite{Ha2}. Let $X$ be a simplicial cell complex, and
$G(X)=(X,G_\sigma,\psi_a,g_{a,b})$ be a complex of groups defined over
$X$. Here $a\in E(X)$ is an edge of the barycentric subdivision of
$X$, and $\sigma\in V(X)$ is a cell of $X$. Each
$a\in E(X)$ is naturally oriented. We write $i(a)$ for the
initial point of $a$ and $t(a)$ for the terminal point of $a$. The
symbol $a^{-1}$ denotes an edge $a\in E(X)$ with the inversed
orientation so that $i(a^{-1})=t(a)$ and $t(a^{-1})=i(a)$. Let
$\sigma_0\in V(X)$ be a fixed cell of $X$. 

The fundamental group $\pi_1(G(X),\sigma_0)$ of the complex of groups
$G(X)$ with respect to the base point $\sigma_0$ is defined as the
group generated by all of the following words
$$
g_0e_1g_1e_2g_2\cdots e_ng_n, \leqno (3.2.4)
$$
where $g_0\in G_{\sigma_0}$, each $e_i$ is either $a$ or $a^{-1}$ for
some $a\in E(X)$ such that $i(e_1)=\sigma_0$, $t(e_i)=i(e_{i+1})$,
$t(e_n)=\sigma_0$,  and $g_i\in G_{t(e_i)}$. Besides the relations
inherited from each group $G_\sigma$, $\sigma\in V(X)$, the words 
in $(3.2.4)$ are subject to the following further relations
$$
aa^{-1}=a^{-1}a=1, \hspace{4mm}\psi_a(g)=a^{-1}ga, \forall g\in G_{i(a)}, 
\hspace{4mm} ab=bag_{a,b}, \leqno (3.2.5)
$$
where $a,b$ are any composable edges in $E(X)$. 

We equip $X$ with the orbihedron structure associated to $G(X)$ and
regard $X$ canonically as an orbispace in the sense of this paper.
For any given cell $\sigma_0\in V(X)$, we define a base-point structure
$\underline{\sigma_0}=(o,U_o,\hat{o})$ of $X$ as follows:
$U_o=St\sigma_0$, $o$ is the barycenter of $\sigma_0$, and
$\hat{o}=p_{\sigma_0}^{-1}(o)\in St\tilde{\sigma}_0$ (note that 
$p_{\sigma_0}:St\tilde{\sigma}_0\rightarrow St\sigma_0$ is one to one
restricted on $\tilde{\sigma}_0$). 

\vspace{2mm}

\noindent{\bf Theorem 3.2.10:}
{\em The fundamental group $\pi_1(G(X),\sigma_0)$ is isomorphic to 
$\pi_1(X,\underline{\sigma_0})$.
}

\vspace{2mm}

\noindent{\bf Proof:}
The orbihedron structure on $X$ given in \cite{Ha2} and the
corresponding orbispace structure on $X$ in the sense of this paper
were discussed in Example 2.1.3 c. Let us first briefly review the latter.

Given an orbihedron structure on a simplicial cell complex $X$ as
defined in \cite{Ha2}, an orbispace structure as defined in
Definition 2.1.2 was canonically constructed on $X$ as follows. We let
$\U=\{St\sigma,\sigma\in V(X)\}$. (Although $\U$ does not form a base
of the underlying topological space of the simplicial cell complex
$X$, $\U$ determines an orbispace structure on it in an obvious way.)
The $G$-structure of $St\sigma$ is
$(St\tilde{\sigma},G_\sigma,p_\sigma)$. For
any two different cells $\sigma$ and $\tau$, $St\sigma\cap
St\tau\neq\emptyset$ if and only if $\dim\sigma\neq \dim\tau$ and
there is an edge $a\in E(X)$ such that $i(a)=\sigma$ and $t(a)=\tau$
(assuming $\dim\sigma > \dim\tau$), and in the latter case,
$St\sigma\cap St\tau$ is the disjoint union of all the $Sta$'s where
$a\in E(X)$ satisfy $i(a)=\sigma$ and $t(a)=\tau$. The set of
transition maps $Tran(St\sigma,St\tau)$ is given by
$$
Tran(St\sigma,St\tau)=\{(g\circ f_a,g\psi_ag^{-1})|a\in E(X),
i(a)=\sigma, t(a)=\tau, g\in G_{t(a)}\}, \leqno (3.2.6)
$$
where $f_a$, which is a $\psi_a$-equivariant simplicial map
$p_\sigma^{-1}(Sta) \rightarrow p_\tau^{-1}(Sta)$ for any $a\in E(X)$
with $i(a)=\sigma$ and $t(a)=\tau$, is given as part of the data of
the orbihedron structure on $X$. 

We shall construct an isomorphism $\Psi:\pi_1(X,\underline{\sigma_0})
\rightarrow \pi_1(G(X),\sigma_0)$ as follows. Suppose that
$u=(\{I_i\},\{U_i\},\{\gamma_i\},\{\xi_{ji}\})$ is a representing
system of a based loop in $(X,\underline{\sigma_0})$, where the index
$i$ is running from $0$ to $n$, each $U_i=St\tau_i$ for some
$\tau_i\in V(X)$ with $\tau_0=\sigma_0$. As a notational convention,
we assume $j=i+1$ for $i<n$ and $j=0$ for $i=n$. It follows easily
that \\
(1) if $\dim\tau_j<\dim\tau_i$, then $\xi_{ji}=(g_i\circ f_{a_i},
g_i\psi_{a_i}g_i^{-1})$ for some $a_i\in E(X)$ with $i(a_i)=\tau_i$
and $t(a_i)=\tau_j$, and some $g_i\in G_{\tau_j}$, \\
(2) if $\dim\tau_j>\dim\tau_i$, then $\xi_{ji}^{-1}=(g_i\circ f_{a_i},
g_i\psi_{a_i}g_i^{-1})$ for some $a_i\in E(X)$ with $i(a_i)=\tau_j$
and $t(a_i)=\tau_i$, and some $g_i\in G_{\tau_i}$,\\
(3) if $\dim\tau_j=\dim\tau_i$, then $\tau_j=\tau_i$, and
$\xi_{ji}=g_i$ for some $g_i\in G_{\tau_i}$.\\
For case (1), we set $A_{ij}=a_ig_i$, for case (2), we set
$A_{ij}=g_i^{-1}a_i^{-1}$, and for case (3), we set $A_{ij}=g_i$. 
We define $\Psi(u)$ to be the class of the word
$$
A_{01}A_{12}\cdots A_{(n-1)n}A_{n0} \leqno (3.2.7)
$$
in the fundamental group $\pi_1(G(X),\sigma_0)$.  One can verify that
$\Psi(u)$ depends only on the equivalence class $[u]$ of $u$, with the
help of relations in $(3.2.5)$. On the other hand, appealing to the
description of a canonical neighborhood of a based loop as given in
Lemma 3.1.2 and the fact that the groups $\{G_\sigma,\sigma\in V(X)\}$
are discrete, we see that $\Psi(u)$ is independent of the homotopy
class of $[u]$, hence $\Psi$ descents to a map
$\pi_1(X,\underline{\sigma_0})\rightarrow \pi_1(G(X),\sigma_0)$, which
clearly is a homomorphism by the nature of construction. 

For the surjectivity of 
$\Psi:\pi_1(X,\underline{\sigma_0})\rightarrow \pi_1(G(X),\sigma_0)$,
given any word $w$ as in $(3.2.4)$, we can write it as a product of
$A_{ij}$'s as in $(3.2.7)$. Then we can easily construct a system $u$
in which the $\xi_{ji}$'s give rise to the $A_{ij}$'s from the word
$w$. It follows easily that $\Psi(u)=[w]$ in $\pi_1(G(X),\sigma_0)$.
Hence $\Psi$ is surjective.

For the injectivity of 
$\Psi:\pi_1(X,\underline{\sigma_0})\rightarrow \pi_1(G(X),\sigma_0)$,
we observe that if for some representing system $u$, the word
defined in $(3.2.7)$ is equivalent to the trivial word through a
sequence of cancellations using the relations in $(3.2.5)$, then there
is a corresponding sequence of representing systems of based loops, 
$\{u_i|i=0,\cdots m\}$, 
such that $u_0=u$, $[u_i]$ is homotopic to $[u_{i+1}]$, and $u_m$ is a
closed loop in $St\tilde{\sigma}_0$. The claim that 
$[u_i]$ is homotopic to $[u_{i+1}]$ for each $i$ also uses the fact
that each $St\tilde{\sigma}$, $\sigma\in V(X)$, and
$p_\sigma^{-1}(Sta)$ with $i(a)=\sigma$ is contractable. Now it is
easily seen that $[u]$ is null-homotopic, hence $\Psi$ is injective. 
This concludes the proof of the theorem.

\subsection{Relative homotopy groups}

\hspace{5mm}
This section concerns, for a given based pseudo-embedding, the construction 
of {\it relative homotopy groups} and establishment of the corresponding 
{\it exact homotopy sequence} $(1.7)$ associated to it. 

\vspace{1.5mm}

Given based orbispaces $(X,\underline{p})$ and $(Y,\underline{q})$ where the
base-point structures are specified as $\underline{p}=(p,U_o,\hat{p})$ and
$\underline{q}=(q,V_o,\hat{q})$, 
recall that a based pseudo-embedding $\tilde{i}: (Y,\underline{q})\rightarrow 
(X,\underline{p})$ is a based morphism which has a representing system 
$\sigma=(\{V_\alpha\},\{U_\alpha\},\{i_\alpha\},\{\rho_{\beta\alpha}\})$,
where each $i_\alpha:\widehat{V_\alpha}\rightarrow\widehat{U_\alpha}$ is an 
embedding and each $\rho_\alpha:G_{V_\alpha}\rightarrow G_{U_\alpha}$ is a 
monomorphism such that the natural projection $G_{U_\alpha}\rightarrow
G_{U_\alpha}/\rho_\alpha(G_{V_\alpha})$ is a weak fibration. 

Fixing a representing system 
$\sigma=(\{V_\alpha\},\{U_\alpha\},\{i_\alpha\},\{\rho_{\beta\alpha}\})$,
$\alpha\in\Lambda$, of the based pseudo-embedding, we will first construct
the {\it based relative loop space} $\Omega(X,Y,\sigma)$. An element of 
$\Omega(X,Y,\sigma)$, which is called a {\it based relative loop} in 
$(X,Y,\sigma)$, is an equivalence class of based systems from $[0,1]$ to
$X$ satisfying a certain boundary condition specified by $\sigma$. 
More concretely, let $I_i$, $i=0,1,\cdots,n$, be intervals such that 
$[0,1]=\cup_{i=0}^n I_i$ and $I_i\cap I_j\neq \emptyset$ iff $j=i+1$, and
let $\{U_i\}$ be a collection of basic open sets of $X$ such that 
$U_0=U_o$ and $U_n=U_\alpha$ for some $\alpha\in\Lambda$.
We consider system $(\{I_i\},\{U_i\},V_\alpha,\{f_i\},\{\xi_{ji}\})$ where
each $f_i:I_i\rightarrow \widehat{U_i}$ is continuous, and $\xi_{ji}\in
Tran(U_i,U_j)$ such that $f_j=\xi_{ji}\circ f_i$ on $I_i\cap I_j$. Furthermore,
we require that $f_0(0)=\hat{p}$ and $f_n(1)\in i_\alpha(\widehat{V_\alpha})
\subset\widehat{U_\alpha}=\widehat{U_n}$. Two systems 
$(\{f_{i,1}\},\{\xi_{ji,1}\})$, $(\{f_{i,2}\},\{\xi_{ji,2}\})$ 
defined over $(\{I_i\},\{U_i\},V_\alpha)$ are 
said to be {\it isomorphic} if there are automorphisms of 
$(\widehat{U_i},G_{U_i},\pi_{U_i})$, $\delta_i\in G_{U_i}$, $i=1,2,\cdots,n$, 
such that $f_{i,2}=\delta_i\circ f_{i,1}$, 
$\xi_{ji,2}=\delta_j\circ\xi_{ji,1}\circ\delta_i^{-1}$,  and 
$\delta_n\in\rho_\alpha(G_{V_\alpha})$. The refinement relation is defined 
as follows. A system $(\{J_k\},\{W_k\},V_\beta,\{g_k\},\{\eta_{lk}\})$ (here 
each $W_k$ is a basic open set of $X$), where $k$ is running from $0$ to $m$, 
is a refinement of $(\{I_i\},\{U_i\},V_\alpha,\{f_{i}\},\{\xi_{ji}\})$ if (1) 
$\{J_k\}$ is a refinement of $\{I_i\}$, given by an index mapping 
$\theta:\{0,\cdots,m\}\rightarrow \{0,\cdots,n\}$, such that $\theta(0)=0$ and 
$\theta(m)=n$, (2) there are transition maps 
$\jmath_k\in Tran(W_k,U_{\theta(k)})$ for $k=1,\cdots,m$ which satisfy 
the following conditions: $g_k=\jmath_k^{-1}\circ f_{\theta(k)}|_{J_k}$ and 
$\eta_{lk}=\jmath_l^{-1}\circ\xi_{\theta(l)\theta(k)}\circ\jmath_k$, and 
(3) there is a transition map $\imath\in Tran(V_\beta,V_\alpha)$ such that 
$\jmath_m=\rho_{\alpha\beta}(\imath)$ (note that 
$\jmath_m\in Tran(W_m,U_n)=Tran(U_\beta,U_\alpha)$). The isomorphism relation 
and the refinement relation together determines an equivalence relation on 
these systems. 

\vspace{1.5mm}

\noindent{\bf Definition 3.3.1:}
{\em The based relative loop space of $(X,Y,\sigma)$, denoted by 
$\Omega(X,Y,\sigma)$, is the set of equivalence classes of systems as 
described in the previous paragraph.
}

\hfill $\Box$

As in the case of based loop space, a natural compact-open topology can be 
given on $\Omega(X,Y,\sigma)$. Let $\tilde{f}_0$ be a based relative loop and 
$\tau_0=(\{I_i\},\{U_i\},V_\alpha,\{f_{i,0}\},\{\xi_{ji,0}\})$ be a 
representing system of $\tilde{f}_0$. Given $K=\{K_i\}$ where each $K_i$ is 
a compact subset of $I_i$, $W=\{W_i\}$ where each $W_i$ is an open subset of 
$\widehat{U_i}$, and $A=\{A_{ji}\}$ where each $A_{ji}$ is a neighborhood of 
$\xi_{ji,0}$ in $Tran(U_i ,U_j)$, such that $f_{i,0}(K_i)\subset W_i$, we 
similarly define a subset $\O(\tau_0, K,W,A)$ of $\Omega(X,Y,\sigma)$ by
$$
\O(\tau_0, K,W,A):=\{\tilde{f}|\tilde{f}=
[(\{I_i\},\{U_i\},V_\alpha,\{f_{i}\},\{\xi_{ji}\})], f_i(K_i)\subset W_i, 
\xi_{ji}\in A_{ji}\}. \leqno (3.3.1)
$$
The topology on $\Omega(X,Y,\sigma)$ is generated by the set of all 
$\O(\tau_0, K,W,A)$.

Given a representing system 
$\tau_0=(\{I_i\},\{U_i\},V_\alpha,\{f_{i,0}\},\{\xi_{ji,0}\})$ of 
$\tilde{f}_0\in\Omega(X,Y,\sigma)$, where $i$ is running from $0$ to $n$, 
a based relative loop $\tilde{f}$ in a neighborhood of $\tilde{f}_0$ also 
has a canonical form as in the case of based loops. More concretely, 
$\tilde{f}$ can be represented by a system 
$(\{I_i\},\{U_i\},V_\alpha,\{f_{i}\},\{\xi_{ji}\})$, where 
$\xi_{ji}=\xi_{ji,0}$ for $j=1,\cdots,n-1$, and the last transition map 
$\xi_{n(n-1)}=\delta\circ\xi_{n(n-1),0}$ for some 
$\delta\in G_{U_n}=G_{U_\alpha}$, which is determined only as an element 
in the space of cosets $G_{U_\alpha}/\rho_\alpha(G_{V_\alpha})$. 
However, the assumption that the natural projection 
$G_{U_\alpha}\rightarrow G_{U_\alpha}/\rho_\alpha(G_{V_\alpha})$ is a weak 
fibration bears the following important consequence: a continuous map from 
a Euclidean ball $D$ into a neighborhood of $\tilde{f}_0$ can be represented 
by a system $(\{I_i\},\{U_i\},V_\alpha,\{f_{i}\},\{\xi_{ji}\})$, 
where each $f_i: D\times I_i\rightarrow \widehat{U_i}$ is continuous, 
$\xi_{ji}=\xi_{ji,0}$ for $j=1,\cdots,n-1$, and 
$\xi_{n(n-1)}=\delta\circ\xi_{n(n-1),0}$ where 
$\delta: D\rightarrow G_{U_\alpha}$ is continuous. 

We point out that the homeomorphism class of based relative loop space 
$\Omega(X,Y,\sigma)$ only depends on the based pseudo-embedding 
$\tilde{i}: (Y,\underline{q})\rightarrow (X,\underline{p})$, not on the 
representing system $\sigma$. When there is no confusion, we simply denote 
the based relative loop space by $\Omega(X,Y,\tilde{i})$.

Finally, we observe that each based relative loop 
$\tilde{f}\in \Omega(X,Y,\tilde{i})$ induces a path 
$f:[0,1]\rightarrow X$ such that $f(0)=p$ and $f(1)\in i(Y)$. Moreover, 
$\tilde{f}$ determines a point in $Y$, which will be denoted by 
$\partial\tilde{f}$, as follows: let 
$(\{I_i\},\{U_i\},V_\alpha,\{f_{i}\},\{\xi_{ji}\})$, where $i$ is running 
from $0$ to $n$, be a representing system of $\tilde{f}$, we set 
$$
\partial\tilde{f}=\pi_{V_\alpha}(i_\alpha^{-1}(f_n(1))),\leqno (3.3.2)
$$
which is easily seen to be independent of the representing system.

There is a special point in $\Omega(X,Y,\tilde{i})$, the constant relative 
loop $\tilde{p}$, which is the equivalence class of the system 
$([0,1],U_o,V_o,f)$ where $f:[0,1]\rightarrow \widehat{U_o}$ is the constant 
map into $\hat{p}$. We fix $\tilde{p}$ as the base point of 
$\Omega(X,Y,\tilde{i})$.

\vspace{1.5mm}

\noindent{\bf Definition 3.3.2:}
{\em The relative homotopy groups $\pi_k(X,Y,\tilde{i})$ associated to a 
based pseudo-embedding 
$\tilde{i}: (Y,\underline{q})\rightarrow (X,\underline{p})$ are defined by 
$$
\pi_k(X,Y,\tilde{i}):=\pi_{k-1}(\Omega(X,Y,\tilde{i}),\tilde{p})
$$
for all $k\geq 1$. In fact, $\pi_1(X,Y,\tilde{i})$ is only a based set in 
general. 
}

\hfill $\Box$

Each element of $\Omega(X,\underline{p})$ naturally gives rise to an element
of $\Omega(X,Y,\sigma)$ as follows. Suppose 
$\tilde{\gamma}\in\Omega(X,\underline{p})$ is represented by a system 
$\tau=(\{I_i\},\{U_i\},\{\gamma_i\},\{\xi_{ji}\})$, where $i$ is running
from $0$ to $n-1$. We decompose $I_0$ as $I_{0,-}\cup I_{0,+}$ along the
base point $\ast\in S^1$, and name $J_0=I_{0,+}$, $J_1=I_1, \cdots, 
J_n=I_{0,-}$. Then the system $\tau$ gives rise to a system 
$j(\tau)=(\{J_k\},\{W_k\},V_o,\{\gamma_k^\prime\},\{\eta_{lk}\})$ where 
$W_k=U_k$ for $0\leq k\leq n-1$, $W_n=U_0$, $\gamma_0^\prime=
\gamma_0|_{I_{0,+}}$, $\gamma_k^\prime=\gamma_k$ for $0\leq k\leq n-1$, 
$\gamma_n^\prime=\gamma_0|_{I_{0,-}}$, and $\eta_{lk}=\xi_{lk}$ for $1\leq l
\leq n-1$, $\eta_{n(n-1)}=\xi_{0(n-1)}$. The map $j$ clearly descents to a 
continuous map $j:\Omega(X,\underline{p})\rightarrow\Omega(X,Y,\sigma)$.
The constant loop $\tilde{p}$ is sent to the constant relative loop $\tilde{p}$
in $\Omega(X,Y,\sigma)$. The geometric interpretation of going from a based
loop $\tilde{\gamma}$ to the based relative loop $j(\tilde{\gamma})$ is that
in the definition of isomorphism relation, the automorphism $\delta_o$ 
on the $I_{0,-}$ part is set free to be anything in $\rho_o(G_{V_o})$ 
(cf. $(2.2.2a-b)$), and similar things occur for the refinement relation. 

Now we consider the following sequence of continuous maps
$$
(\Omega(Y,\underline{q}),\tilde{q})\stackrel{i}{\rightarrow}
(\Omega(X,\underline{p}),\tilde{p})\stackrel{j}{\rightarrow}
(\Omega(X,Y,\sigma),\tilde{p}), \leqno (3.3.3)
$$
where $i$ is the continuous map induced by the based pseudo-embedding
$\tilde{i}:(Y,\underline{q})\rightarrow (X,\underline{p})$.

\vspace{1.5mm}

\noindent{\bf Theorem 3.3.3:}
{\em The sequence $(3.3.3)$ induces a long exact sequence
$$
\rightarrow \pi_{k+1}(X,Y,\tilde{i})\stackrel{\partial}{\rightarrow}
\pi_{k}(Y,\underline{q})\stackrel{i_{\#}}{\rightarrow}
\pi_k(X,\underline{p})\stackrel{j_{\#}}{\rightarrow}
\pi_k(X,Y,\tilde{i})\stackrel{\partial}{\rightarrow} 
\cdots \stackrel{\partial}{\rightarrow} \pi_0(Y,q)
\stackrel{i_{\#}}{\rightarrow} \pi_0(X,p).
\leqno (3.3.4)
$$
}

\vspace{2mm}

We shall first prove a technical lemma. Recall that the suspension 
$(SK,\ast)$ of a based space $(K,\ast)$ is the quotient space
$$
[0,1]\times K/\{0\}\times K\cup [0,1]\times\{\ast\}\cup\{1\}\times K. 
$$
There is a homotopy associative comultiplication $\mu^\prime:SK\rightarrow 
SK\vee SK$ defined by 
$$
\mu^\prime [t,x]=\left\{\begin{array}{cc}
([2t,x],\ast) & 0\leq t\leq \frac{1}{2}\\
(\ast, [2t-1,x]) & \frac{1}{2}\leq t \leq 1,
\end{array} \right.
$$
and a homotopy inverse $\nu^\prime:SK\rightarrow SK$ defined by $\nu^\prime 
[t,x]\rightarrow [1-t,x]$, under which $(SK,\ast)$ is an H-cogroup. For a 
based space $(K,\ast)$, its cone $(CK,\ast)$ is the based space defined as 
the quotient space $[0,1]\times K/\{0\}\times K\cup [0,1]\times\{\ast\}$. 
The space $(K,\ast)$ can be regarded as a subspace of its cone $(CK,\ast)$ 
via the embedding $x\mapsto [1,x]$.

\vspace{1.5mm}

\noindent{\bf Lemma 3.3.4:}
{\em Let $(K,\ast)$ be a compact space, which underlies a 
simplicial complex. Then there is a correspondence $\partial$, which assigns 
to each continuous map $f: (SK,\ast)\rightarrow (\Omega(X,Y,\tilde{i}),
\tilde{p})$ a continuous map $\partial f: (K,\ast)\rightarrow 
(\Omega(Y,\underline{q}),\tilde{q})$, such that there is a continuous map 
$F: (CK,\ast)\rightarrow (\Omega(X,\underline{p}),\tilde{p})$ whose 
restriction on $(K,\ast)\subset (CK,\ast)$ equals $i\circ \partial f$, where 
$i: (\Omega(Y,\underline{q}),\tilde{q})\rightarrow (\Omega(X,\underline{p}),
\tilde{p})$ is the map in $(3.3.3)$. Moreover, we have 
$\partial(\mu^\prime(f,g))=\partial f\#\partial g$ and 
$\partial(f\circ\nu^\prime)=\nu(\partial f)$, 
where $\#,\nu$ are the homotopy associative multiplication and homotopy 
inverse of the H-group $(\Omega(Y,\underline{q}),\tilde{q})$.
}

\vspace{2mm}

\noindent{\bf Proof:}
To begin with, we fix the notation for the base-point structures 
$\underline{p}=(p,U_o,\hat{p})$ of $X$ and $\underline{q}=(q,V_o,\hat{q})$ 
of $Y$, and fix a representing system 
$$
\sigma=(\{V_\alpha\},\{U_\alpha\},\{i_\alpha\},\{\rho_{\beta\alpha}\}),
\hspace{3mm} \alpha\in\Lambda, \leqno (3.3.5)
$$
of the pseudo-embedding $\tilde{i}:(Y,\underline{q})\rightarrow 
(X,\underline{p})$.

We regard $f:(SK,\ast)\rightarrow (\Omega(X,Y,\tilde{i}),\tilde{p})$ as a map
from $[0,1]\times K$ to $\Omega(X,Y,\tilde{i})$. We can subdivide $K$ into a
union of simplexes $\cup_{a\in A} K_a$, and the interval $[0,1]$ into a union
$[0,1]=\cup_{k=0}^{m-1}[s_k,s_{k+1}]$, such that the restriction 
$f|_{[s_k,s_{k+1}]\times K_a}$ for each $(k,a)$ lies in a canonical 
neighborhood of a based relative loop associated to a representing system.
More precisely, $f|_{[s_k,s_{k+1}]\times K_a}$ is represented by a system
$$
\tau_{(k,a)}=(\{[s_k,s_{k+1}]\times K_a\times I_i\},\{U_i^{(k,a)}\},
V_{\alpha(k,a)}, \{f_i^{(k,a)}\},
\{\xi_{ji}^{(k,a)}\}), \hspace{3mm} i=0,1,\cdots,n,\leqno (3.3.6)
$$
where $U_0^{(k,a)}=U_o$, $U_n^{(k,a)}=U_{\alpha(k,a)}$ for some index 
$\alpha(k,a)\in\Lambda$ such that there is an embedding $i_{\alpha(k,a)}:
\widehat{V_{\alpha(k,a)}}\rightarrow \widehat{U_{\alpha(k,a)}}$ given 
as part of the data of $(3.3.5)$. Each 
$f_i^{(k,a)}:[s_k,s_{k+1}]\times K_a\times I_i\rightarrow 
\widehat{U_i^{(k,a)}}$ is a continuous map such that 
$f_0^{(k,a)}(s,x,0)=\hat{p}$ and $f_n^{(k,a)}(s,x,1)\in i_{\alpha(k,a)}
\widehat{V_{(k,a)}}$ for any $(s,x)\in [s_k,s_{k+1}]\times K_a$, 
each $\xi_{ji}^{(k,a)}$ is a transition map
in $Tran(U_i^{(k,a)},U_j^{(k,a)})$ which is constant in $(s,x)\in
[s_k,s_{k+1}]\times K_a$ for $j=1,\cdots,n-1$, and $\xi_{n(n-1)}^{(k,a)}:
[s_k,s_{k+1}]\times K_a\rightarrow Tran(U_{n-1}^{(k,a)},U_n^{(k,a)})$ is a 
continuous map, and $f_i^{(k,a)}$, $\xi_{ji}^{(k,a)}$ satisfy the compatibility
condition $f_j^{(k,a)}=\xi_{ji}^{(k,a)}\circ f_i^{(k,a)}$ on 
$[s_k,s_{k+1}]\times K_a\times (I_i\cap I_j)$. Here since $[0,1]\times K$ is
compact, we can take $\{I_i\}$ to be independent of $(k,a)$, by passing to 
refined systems. If we set $K_o$ to be the simplex of $K$ which contains
the base point $\ast\in K$, then we can assume that 
$$
U_i^{(k,o)}=U_i^{(0,a)}=U_i^{(m-1,a)}=U_o \leqno (3.3.7)
$$
for any $a\in A, k\in\{0,\cdots,m-1\}$ and $i\in\{0,\cdots,n\}$, and
$$
f_i^{(k,o)}(s,\ast,t)=f_i^{(0,a)}(0,x,t)=f_i^{(m-1,a)}(1,x,t)=\hat{p}
\leqno (3.3.8)
$$
for any $s\in [s_k,s_{k+1}], x\in K_a$ and $t\in I_i$, and 
$$
\xi_{ji}^{(k,o)}=\xi_{ji}^{(0,a)}=\xi_{ji}^{(m-1,a)}=1_{G_{U_o}}
\leqno (3.3.9)
$$
for $j=1,2,\cdots,n-1$, and $\xi_{n(n-1)}^{(k,o)},
\xi_{n(n-1)}^{(0,a)}$ and $\xi_{n(n-1)}^{(m-1,a)}$ are in a prescribed 
neighborhood of $1_{G_{U_o}}$, with 
$$
\xi_{n(n-1)}^{(k,o)}(s,\ast)=\xi_{n(n-1)}^{(0,a)}(0,x)=
\xi_{n(n-1)}^{(m-1,a)}(1,x)=1_{G_{U_o}} \leqno (3.3.10)
$$
for any $s\in [s_k,s_{k+1}]$ and $x\in K_a$. 
Finally, we define for each $a\in A$ a collection of continuous maps 
indexed by $k\in\{0,\cdots,m-1\}$:
$$
g_k^a(x,s):=i^{-1}_{\alpha(k,a)}(f_n^{(k,a)}(s,x,1)), \hspace{2mm}
\forall x\in K_a, s\in [s_k,s_{k+1}], \leqno (3.3.11)
$$
each of which is from $K_a\times [s_k,s_{k+1}]$ to 
$\widehat{V_{\alpha(k,a)}}$. The $g_k^a$'s satisfy 
$$
g^o_k(\ast,s)=g^a_0(x,0)=g^a_{k-1}(x,1)=\hat{q}\in \widehat{V_o}
\leqno (3.3.12)
$$
for any $k\in\{0,\cdots,m-1\}$, $a\in A$,
where $\ast\in K_o$ is the base point of $K$, $s\in [s_k,s_{k+1}]$, 
and $x\in K_a$.

\vspace{1.5mm}

Now restricted to $\{s_k\}\times K_a$, since the systems $\tau_{(k-1,a)}$
and $\tau_{(k,a)}$ represent the same family of based relative loops
$f_{\{s_k\}\times K_a}$, there are maps $\eta_{k(k-1)}^{(i,a)}:K_a\rightarrow
Tran(U_i^{(k-1,a)},U_i^{(k,a)})$, and $\lambda_{k(k-1)}^a:K_a\rightarrow
Tran(V_{\alpha(k-1,a)},V_{\alpha(k,a)})$ such that 
$$
\eta_{k(k-1)}^{(0,a)}=1_{G_{U_o}}, \hspace{2mm}
\eta_{k(k-1)}^{(n,a)}=\rho_{\alpha(k,a)\alpha(k-1,a)}(\lambda_{k(k-1)}^a),
\leqno (3.3.13)
$$
and for any $x\in K_a$, we have 
$$
\xi_{ji}^{(k,a)}(s_k,x)=\eta_{k(k-1)}^{(j,a)}(x)\circ 
\xi_{ji}^{(k-1,a)}(s_k,x)\circ (\eta_{k(k-1)}^{(i,a)})^{-1}(x),
\leqno (3.3.14 a) 
$$
and 
$$
f_i^{(k,a)}(s_k,x,\cdot)=\eta_{k(k-1)}^{(i,a)}(x)\circ 
f_i^{(k-1,a)}(s_k,x,\cdot). \leqno (3.3.14 b)
$$
Since $\eta_{k(k-1)}^{(0,a)}=1_{G_{U_o}}$, $\xi_{ji}^{(k,a)}$ is independent
of $x\in K_a$ for $j=1,\cdots,n-1$ and $\xi_{n(n-1)}^{(k,a)}$ is continuous
in $x\in K_a$, we conclude that $\eta_{k(k-1)}^{(i,a)}$ is constant for
$i=0,1,\cdots,n-1$ and $\eta_{k(k-1)}^{(n,a)}$ and $\lambda_{k(k-1)}^a$ are 
continuous in $x\in K_a$. Finally, $(3.3.14 b)$ implies the following 
compatibility conditions on the $g^a_k$'s:
$$
g^a_k(x,s_k)=\lambda_{k(k-1)}^a(x)(g^a_{k-1}(x,s_k)). \leqno (3.3.15)
$$

\vspace{1.5mm}

Now we identify $0,1 \in [0,1]$ and take it as the base point, and then choose
small enough open intervals $J_0,J_1,\cdots,J_{m-1}$ such that $[0,s_1]\cup
[s_{m-1},1]\subset J_0$, $[s_1,s_2]\subset J_1, \cdots, [s_{m-2},s_{m-1}]
\subset J_{m-1}$. It is easily seen that the maps $\{g^a_k\}$ and transition
maps $\lambda_{k(k-1)}^a$ can be put together to define a collection of 
continuous families of systems indexed by $a\in A$ and parametrized by 
$x\in K_a$:
$$
\tau_a(x)=(\{J_k\},\{V_{\alpha(k,a)}\},\{g_k^a(x,\cdot)\},
\{\lambda_{k(k-1)}^a(x)\}), \hspace{2mm} x\in K_a \leqno (3.3.16)
$$
from $(S^1,\ast)$ to $(Y,\underline{q})$, such that $\tau_o(\ast)$, $\ast\in
K_o\subset K$, is the canonical representing system of the base point
$\tilde{q}$ in $\Omega(Y,\underline{q})$. However, in order to patch this
collection of parametrized systems to yield the desired map $\partial f: 
(K,\ast)\rightarrow (\omega(Y,\underline{q}),\tilde{q})$, we need the following

\vspace{1.5mm}

\noindent{\bf Sublemma 3.3.5:}
{\em Suppose $\theta_k(z)=(\{[t_i,t_{i+1}]\},\{U_i^k\},\{\gamma_i^k(z,\cdot)\},
\{\lambda_{ji}^k(z)\})$, $k=1,2$, are two continuous families of systems 
parametrized by $z\in Z$, where each $\gamma_i^k:Z\times [t_i,t_{i+1}]
\rightarrow \widehat{U_i^k}$ and $\lambda_{ji}^k:Z\rightarrow 
Tran(U_i^k,U_j^k)$ is continuous, and there are continuous maps 
$\xi^i:Z\times [t_i,t_{i+1}]\rightarrow Tran(U_i^1,U_i^2)$ such that 
$$
\gamma_i^2=\xi^i\circ\gamma_i^1,\hspace{2mm}
\lambda_{(i+1)i}^2=\xi^{(i+1)}(\cdot,t_{i+1})\circ\lambda_{(i+1)i}^1\circ 
(\xi^i)^{-1}(\cdot,t_{i+1}).
$$
Then there is a canonically constructed continuous family of systems 
$$
\Theta(z,s)=(\{[t_i,t_{i+1}]\},\{U_i^2\},\{\gamma_i(z,s,\cdot)\},
\{\lambda_{ji}(z,s)\})
$$ 
parametrized by $Z\times [1,2]$ such that $\Theta(z,2)=\theta_2(z)$ and
$\Theta(z,1)$ is isomorphic to $\theta_1(z)$, and each $\gamma_0(z,s,0)$ 
and $\gamma_{n}(z,s,1)$ is constant in $s\in [1,2]$. 
Here $[0,1]=\cup_{i=0}^{n}[t_i,t_{i+1}]$ is a subdivision of $[0,1]$.
}

\vspace{2mm}

\noindent{\bf Proof:}
For $0\leq i\leq n-1$, set 
$$
\gamma_{i,1}(z,t)=\xi^i(z,t_i)(\gamma_i^1(z,t)),\hspace{2mm}
\lambda_{i(i-1),1}(z)=\xi^{i}(z,t_{i})\circ \lambda_{i(i-1)}^1(z)\circ
(\xi^{(i-1)}(z,t_{(i-1)}))^{-1},
$$ 
and for the remaining case, set 
$$
\gamma_{n,1}(z,t)=\xi^n(z,1)(\gamma_n^1(z,t)),\hspace{2mm}
\lambda_{n(n-1),1}(z)=\xi^n(z,1)\circ \lambda_{n(n-1)}^1(z)
\circ (\xi^{(n-1)}(z,t_{n-1}))^{-1}.
$$
Then the parametrized systems 
$\theta^1(z)=(\{[t_i,t_{i+1}]\},\{U_i^2\},\{\gamma_{i,1}(z,\cdot)\},
\{\lambda_{ji,1}(z)\})$, $z\in Z$, are isomorphic to $\theta_1(z)$.
Similarly, we define for $0\leq i\leq n-1$, 
$$
\gamma_i(z,s,t)=\xi^i_s(z,t)(\gamma_i^1(z,t)),\hspace{2mm}
\lambda_{i(i-1)}(z,s)=\xi^{i}_s(z,t_{i})\circ \lambda_{i(i-1)}^1(z)\circ
(\xi^{(i-1)}_s(z,t_{i}))^{-1},
$$ 
and for the remaining case,
$$
\gamma_n(z,s,t)=\xi^n_s(z,t)(\gamma_n^1(z,t)),\hspace{2mm}
\lambda_{n(n-1)}(z,s)=\xi^n_s(z,t_n)\circ \lambda_{n(n-1)}^1(z)
(\xi^{(n-1)}_s(z,t_n))^{-1},
$$
where $\xi_s^i(z,t)=\xi^i(z,(2-s)t_i+(s-1)t)$ for $0\leq i\leq n-1$, and 
$\xi_s^n(z,t)=\xi^n(z,(2-s)1+(s-1)t)$, $s\in [1,2]$. Then one easily verifies
that $\gamma_{(i+1)}(z,s,t_{i+1})=\lambda_{(i+1)i}(z,s)(\gamma_i(z,s,t_{i+1}))$
so that we have a well-defined continuous family of systems
$$
\Theta(z,s):=(\{[t_i,t_{i+1}]\},\{U_i^2\},\{\gamma_i(z,s,\cdot)\},
\{\lambda_{ji}(z,s)\}),
$$ 
which satisfies that $\Theta(z,2)=\theta_2(z)$, $\Theta(z,1)=\theta^1(z)$ 
which is isomorphic to $\theta_1(z)$, and each $\gamma_0(z,s,0)$ 
and $\gamma_{n}(z,s,1)$ is constant in $s\in [1,2]$. 
\hfill $\Box$

\vspace{1.5mm}

Now back to the proof of Lemma 3.3.4. Suppose simplexes $K_a,K_b$ intersect at
a face $K_{ab}$. Then restricted to $[s_k,s_{k+1}]\times K_{ab}$, since the
systems $\tau_{(k,a)}$ and $\tau_{(k,b)}$ (cf. $(3.3.6)$) represent to the
same family of based relative loops, a similar argument shows that there are
continuous families of transition maps $\zeta_{ba}^k(x,s)\in 
Tran(V_{\alpha(k,a)},V_{\alpha(k,b)})$ parametrized by $(x,s)\in K_{ab}\times
[s_k,s_{k+1}]$ such that 
$$
g_k^b(x,s)=\zeta_{ba}^k(x,s)(g_k^a(x,s)), \hspace{2mm}
\lambda_{k(k-1)}^b(x)=\zeta_{ba}^k(x,s_k)\circ\lambda_{k(k-1)}^a(x)\circ
\zeta_{ba}^{(k-1)}(x,s_k), \leqno (3.3.17)
$$
(cf. $(3.3.16)$ for $g_k^a,\lambda_{k(k-1)}^a$, etc.). Moreover, 
$\{\zeta_{ba}^k\}$ also satisfies
$$
\xi_{n(n-1)}^{(k,b)}(s,x)=\rho_{\alpha(k,b)\alpha(k,a)}(\zeta_{ba}^k(x,s))
\circ\xi_{n(n-1)}^{(k,a)}(s,x)\circ (\epsilon_{ba}^k)^{-1}\leqno (3.3.18)
$$
for some $\epsilon_{ba}^k\in Tran(U_{n-1}^{(k,a)},U_{n-1}^{(k,b)})$ which is
independent of $(x,s)\in K_{ba}\times [s_k,s_{k+1}]$. In the case when each
$\zeta_{ba}^k(x,s)$, $k\in\{0,\cdots,m-1\}$, is constant in 
$s\in [s_k,s_{k+1}]$, the systems $\tau_a(x)$ and $\tau_b(x)$ are isomorphic
along $x\in K_{ba}$, hence define a family of based loops in 
$(Y,\underline{q})$ parametrized by $K_a\cup_{K_{ba}} K_b$. In general, we
apply Sublemma 3.3.5 and obtain a family of based loops in 
$(Y,\underline{q})$ parametrized by $K_a\cup_{K_{ba}\times [1,2]} K_b$.
We then identify $K_a\cup_{K_{ba}\times [1,2]} K_b$ with 
$K_a\cup_{K_{ba}} K_b$ by a deformation retract in $K_b$ to obtain a family
of based loops parametrized by $K_a\cup_{K_{ba}} K_b$. Using $(3.3.18)$, this
amounts to a change of $\xi_{n(n-1)}^{(k,b)}$ by homotopy, with the map 
$f|_{K_a\cup K_b}$ remaining unchanged, only the representing systems 
$\{\tau_{k,b},k\in\{0,\cdots,m-1\}\}$ changed. Now it is easily seen that 
we can construct a continuous map $(K,\ast)\rightarrow 
(\Omega(Y,\underline{q}),\tilde{q})$ by patching the systems $\tau_a$ 
together simplex by simplex, which is defined to be the map $\partial
f$. 

In summary, given any representing systems $\tau_{k,a}$'s of 
$f:(SK,\ast)\rightarrow (\Omega(X,Y,\tilde{i})$ as in
$(3.3.6)$, we can modify them through homotopy to isomorphic systems 
without changing the map $f$, so that the systems in
$(3.3.16)$ can be patched together to defined a map $\partial f:
(K,\ast)\rightarrow (\Omega(Y,\underline{q}),\tilde{q})$. 
The equations $\partial(\mu^\prime(f,g))=\partial f\#\partial g$ and 
$\partial(f\circ\nu^\prime)=\nu(\partial f)$ are clear from the 
construction. If we start with different (but equivalent) 
representing systems of $f$, we get equivalent representing systems
for $\partial f$. In other words, the map $\partial f$ is well-defined.

It remains to construct a continuous map $F:(CK,\ast)\rightarrow
(\Omega(X,\underline{p}),\tilde{p})$ such that $i\circ \partial f=F|_K$ where
$K$ is regarded as a subset of $CK$ via embedding $x\mapsto [1,x]$, and
$i:\Omega(Y,\underline{q})\rightarrow\Omega(X,\underline{p})$ is induced by
the pseudo-embedding $\tilde{i}$. For each $a\in A$, the systems
$$
\tau^{(i,a)}=(\{[s_k,s_{k+1}]\times K_a\times I_i\},\{U_i^{(k,a)}\},
\{f_i^{(k,a)}\},
\{\eta_{k(k-1)}^{(i,a)}\}), \hspace{3mm} k=0,1,\cdots,m-1, \leqno (3.3.19)
$$
where $i$ is running from $0$ to $n$, almost patch together to give a map
$(CK_a,\ast)\rightarrow (\Omega(X,\underline{p}),\tilde{p})$, except that
the transition maps $\xi_{n(n-1)}^{(k,a)}(s,x)$ is not constant in $s\in
[s_k,s_{k+1}]$. But clearly we can overcome this problem by using Sublemma
3.3.5, and obtain the desired map $F$. This concludes the proof of Lemma 3.3.4.
\hfill $\Box$

\vspace{1.5mm}

\noindent{\bf Proof of Theorem 3.3.3:}

The sequence of continuous maps $(3.3.3)$ induces the maps 
$\pi_0(Y,q)\stackrel{i_{\#}}{\rightarrow}\pi_0(X,p)$ and 
$\pi_k(Y,\underline{q})\stackrel{i_{\#}}{\rightarrow}\pi_k(X,\underline{p})
\stackrel{j_{\#}}{\rightarrow}\pi_k(X,Y,\tilde{i})$ for $k\geq 1$, which 
clearly satisfies $j_{\#}\circ i_{\#}=0$. What remains is the definition of 
the connecting homomorphisms $\partial:\pi_{k+1}(X,Y,\tilde{i})\rightarrow 
\pi_k(Y,\underline{q})$ for $k\geq 0$. For the case of $k=0$, given the class 
$[\tilde{\gamma}]$ of any based relative loop $\tilde{\gamma}$, we define
$\partial[\tilde{\gamma}]:=[\partial\tilde{\gamma}]\in\pi_0(Y,q)$, where 
$\partial\tilde{\gamma}$ is defined by $(3.3.2)$. A baby version of 
Lemma 3.3.4 shows that such defined map $\partial:\pi_1(X,Y,\tilde{i})
\rightarrow \pi_0(Y,q)$ is well-defined. Moreover, it is easily seen that the 
sequence $(3.3.4)$ is exact at $\pi_1(X,Y,\tilde{i})$ and $\pi_0 (Y,q)$. For 
the case of $k\geq 1$, given the class $[f]$ of a continuous map 
$f:(S^k,\ast)\rightarrow (\Omega(X,Y,\tilde{i}),\tilde{p})$, we think of $f$ 
as a map from $(SS^{k-1},\ast)$, and define $\partial[f]:=[\partial f]$ where 
$\partial f:(S^{k-1},\ast)\rightarrow (\Omega(Y,\underline{q}),\tilde{q})$ 
is constructed in Lemma 3.3.4.  The properties of $\partial f$ established in 
Lemma 3.3.4 imply that the map $\partial:\pi_{k+1}(X,Y,\tilde{i})\rightarrow 
\pi_k(Y,\underline{q})$ just defined is a well-defined homomorphism which 
satisfies $i_{\#}\circ\partial=0$ and $\partial\circ j_{\#}=0$. In summary, we
have shown the existence of sequence $(3.3.4)$, and proved that the 
composition of any two consecutive homomorphisms is zero.

It remains to show that $(3.3.4)$ is exact at $\pi_k(X,\underline{p})$,
$\pi_k(Y,\underline{q})$ and $\pi_{k+1}(X,Y,\tilde{i})$ for any $k\geq 1$.

\vspace{1.5mm}

\noindent{(1)} Exactness at $\pi_k(Y,\underline{q})$: Suppose a class $[u]\in
\pi_k(Y,\underline{q})$ is represented by a continuous map $u:(S^{k-1},\ast)
\rightarrow (\Omega(Y,\underline{q}),\tilde{q})$ such that there is a 
continuous map $H:(CS^{k-1},\ast)\rightarrow 
(\Omega(X,\underline{p}),\tilde{p})$ whose restriction on the subspace
$(S^{k-1},\ast)\subset (CS^{k-1},\ast)$ equals $i\circ u$. We shall construct
a continuous map 
$f:(SS^{k-1},\ast)\rightarrow (\Omega(X,Y,\tilde{i}),\tilde{p})$ from $H$
such that $\partial f=u$. The construction is parallel to the one given in
Lemma 3.3.4, so we will only sketch the main steps here. We think of $H$ as
a map defined over $[0,1]\times S^{k-1}$, and take a subdivition $[0,1]=
\cup_{k=0}^{m-1}[s_k,s_{k+1}]$ and a subdivition $S^{k-1}=\cup_{a\in A}K_a$
of $S^{k-1}$ into simplexes, such that each $H|_{[s_k,s_{k+1}]\times K_a}$
lies in a canonical neighborhood of a based loop in $(\Omega(X,\underline{p})$
(cf. Lemma 3.1.2). Hence there is a collection of continuous families of 
systems
$$
\tau_{(k,a)}=(\{[s_k,s_{k+1}]\times K_a\times I_i\},\{U_i^{(k,a)}\},
\{h_i^{(k,a)}\},\{\xi_{ji}^{(k,a)}\}), \hspace{3mm} i=0,1,\cdots,n,\leqno 
(3.3.20)
$$
which is indexed by $\{(k,a)\}$ and each family is parametrized by 
$[s_k,s_{k+1}]\times K_a$. The transition maps $\xi_{ji}^{(k,a)}(s,x)
\in Tran(U_i^{(k,a)},U_j^{(k,a)})$ are constant in $(s,x)$ for $j=1,\cdots,n$,
and $\xi_{0n}^{(k,a)}(s,x)\in Tran(U_n^{(k,a)},U_0^{(k,a)})$ is continuous in 
$(s,x)$. Unlike the situation in Lemma 3.3.4, the representing system
of a based loop in a canonical neighborhood is unique, as opposed to the case
of based relative loops where there is ambiguity coming from various choices of
liftings with respect to the weak fibrations $G_{U_\alpha}\rightarrow
G_{U_\alpha}/\rho_{\alpha}(G_{V_\alpha})$ associated to the pseudo-embedding
$\tilde{i}$.  Hence for each $i\in \{0,\cdots,n\}$, the systems $(3.3.20)$ 
give rise to a map $f_i$ from $S^{k-1}\times I_i$ to $\Omega(X,Y,\tilde{i})$.
These maps $f_i$ can be almost patched together to give rise to a continuous
map $f:(CS^{k-1},\ast)\rightarrow (\Omega(X,Y,\tilde{i}),\tilde{p})$, except
for the fact that each $\xi_{0n}^{(k,a)}(s,x)\in Tran(U_n^{(k,a)},U_0^{(k,a)})$
may not be constant in $s$. But certainly we can get around of it by using
Sublemma 3.3.5, and the resulting map $f$ satisfies the condition that
$\partial f=u$.

\noindent{(2)} Exactness at $\pi_{k+1}(X,Y,\tilde{i})$: Let $[f]$ be a class in
$\pi_{k+1}(X,Y,\tilde{i})$, which is represented by a continuous map $f: 
(S^k,\ast)\rightarrow (\Omega(X,Y,\tilde{i}),\tilde{p})$ such that $\partial f:
(S^{k-1},\ast)\rightarrow (\Omega(Y,\underline{q}),\tilde{q})$ is 
null-homotopic. We need to construct a continuous map $g:(S^k,\ast)\rightarrow
(\Omega(X,\underline{p}),\tilde{p})$ such that $j_{\#}[g]=[f]$. This goes as
follows. Let $H:(CS^{k-1},\ast)\rightarrow (\Omega(Y,\underline{q}),\tilde{q})$
be the homotopy between $\partial f$ and $\tilde{q}$. As we did in (1), we can
construct a map $h:(SS^{k-1},\ast)\rightarrow (\Omega(Y,Y,Id),\tilde{q})$
such that $\partial h=\partial f$. We define $g:(SS^{k-1},\ast)\rightarrow
(\Omega(X,\underline{p}),\tilde{p})$ by $g(x)=f(x)\# i(h(x))$, i.e., by 
composing the two families of based relative loops along $\partial (f(x))=
\partial (i(h(x)))$ (here $\partial$ is defined as in $(3.3.2)$). It is clear
from the construction that $j\circ g$ is homotopic to $f$ in 
$(\Omega(X,Y,\tilde{i}),\tilde{p})$.

\vspace{1.5mm}

\noindent{(3)} Exactness at $\pi_k(X,\underline{p})$. 
Let $[u]$ be a class in $\pi_k(X,\underline{p})$ represented by a continuous 
family of based loops $u:(S^{k-1},\ast)\rightarrow (\Omega(X,\underline{p}),
\tilde{p})$ such that $j\circ u$ is homotopic to $\tilde{p}$ in 
$\Omega(X,Y,\tilde{i})$ through map $H:(CS^{k-1},\ast)
\rightarrow \Omega(X,Y,\tilde{i})$. As in the proof of Lemma 3.3.4, we
can subdivide $[0,1]$ into 
$\cup_{k=0}^{m-1}[s_k,s_{k+1}]$ and $S^{k-1}$ into a union of simplexes 
$\cup_{a\in A} K_a$, and represent $H$ by a collection of
systems indexed by $(k,a)$:
$$
\tau_{(k,a)}=(\{[s_k,s_{k+1}]\times K_a\times I_i\},\{U_i^{(k,a)}\},
V_{\alpha(k,a)},\{h_i^{(k,a)}\},\{\xi_{ji}^{(k,a)}\}),
\hspace{3mm} i=0,1,\cdots,n, \leqno (3.3.21)
$$
where $\xi_{ji}^{(k,a)}(s,x)\in Tran(U_i^k,U_j^k)$ is constant in 
$(s,x)\in [s_k,s_{k+1}]\times K_a$ for $j=1,\cdots,n-1$, and 
$\xi_{n(n-1)}^k(s,x)\in Tran(U_{n-1}^k,U_n^k)$ is
continuous in $(s,x)$. Restricting $\tau_{k-1},\tau_k$ to 
$\{s_k\}\times K_a$, since they represent the same based relative loop, 
there are continuous families of transition maps 
$\lambda_{k(k-1)}^{(i,a)}(x)\in Tran(U_i^{(k-1,a)},U_i^{(k,a)}$, $x\in K_a$,
such that 
$$
\begin{array}{c}
h^{(k,a)}_i(s_k,x,\cdot)=\lambda_{k(k-1)}^i(x)\circ 
h^{(k-1,a)}_i(s_k,x,\cdot),\\
\xi_{ji}^k(s_k,x)=\lambda_{k(k-1)}^j(x)\circ\xi_{ji}^{(k-1)}(s_k,x)\circ
(\lambda_{k(k-1)}^i(x))^{-1}.
\end{array} \leqno (3.3.22)
$$
In particular, the systems $(\{[s_k,s_{k+1}]\times \{x\}\},\{U_n^{(k,a)}\},
\{h^{(k,a)}_n(\cdot,x,1)\}, \{\lambda_{k(k-1)}^n(x)\})$, $x\in K_a$,
give rise to a continuous family of based paths of $(X,\underline{p})$. These
based paths are actually based loops and can be patched together using 
Sublemma 3.3.5, to yield a map $f:(S^{k-1},\ast)\rightarrow 
(\Omega(X,\underline{p}),\tilde{p})$ which equals $i\circ {g}$ for some 
${g}:(S^{k-1},\ast)\rightarrow (\Omega(Y,\underline{q}),\tilde{q})$ . 
On the other hand, the systems
$$
\tau^{(i,a)}=(\{[s_k,s_{k+1}]\times K_a\times I_i\},\{U_i^{(k,a)}\},
\{h_i^{(k,a)}\}, \{\lambda_{k(k-1)}^{(i,a)}\}), 
\hspace{2mm} k=0,1,\cdots, m-1, \leqno (3.3.23)
$$
are almost patched together to give a family of paths ${f}_{(x,t)}$, 
$(x,t)\in S^{k-1}\times [0,1]$,  with ${f}_{(x,1)}={f}$ and ${f}_{(x,0)}
=\tilde{p}$, except that $\xi_{n(n-1)}^k(s,x)\in Tran(U_{n-1}^{(k,a)},
U_n^{(k,a)})$ may depend on $s$. But we
can again use Sublemma 3.3.5 to overcome this problem, without changing $u$.
Suppose this is done. Then $u\# f: (S^{k-1},\ast)\rightarrow 
(\Omega(X,\underline{p}),\tilde{p})$ 
defined by $(u\# f)(x):=u(x)\# f(x)$ is homotopic to $\tilde{p}$ in 
$(\Omega(X,\underline{p}),\tilde{p})$ through a 1-family of maps 
$u_t\# {f}_{(x,t)}$ from $(S^{k-1},\ast)$ to $(\Omega(X,\underline{p}),
\tilde{p})$ where $u_t$ is the $S^{k-1}$-family of paths obtained by 
restricting $u(x)$, $x\in S^{k-1}$, to the interval $[0,t]$ 
(at the level of systems). Now recall that ${f}=i\circ {g}$ for some 
${g}:(S^{k-1},\ast)\rightarrow (\Omega(Y,\underline{q}),\tilde{q})$. Hence
$[u]=i_{\#}([g]^{-1})$ and $(3.3.4)$ is exact at $\pi_k(X,\underline{p})$.

\subsection{Relation to Borel construction}

\hspace{6mm}
In the equivariant category, there is a so-called Borel construction,
which associates each $G$-space $(Y,G)$ with a fiber bundle $Y_G$ over the
classifying space $BG$ of $G$ with fiber $Y$. The fiber bundle $Y_G$ is
called the Borel space, and is defined by $Y_G:=EG\times_G Y$ where $pr:EG
\rightarrow BG$ is the universal principal $G$-bundle
(cf. \cite{Bo}). 
We will fix Milnor's construction of $EG$ as our choice throughout
(cf. \cite{Mil}). 
Recall that each element of Milnor's $EG$ is denoted $\langle
x,t\rangle$ and written 
$$
\langle g,t\rangle=(t_0g_0,t_1g_1,\cdots,t_kg_k,\cdots )
$$
where each $g_i\in G$ and $t_i\in [0,1]$ such that only a finite number of
$t_i\neq 0$ and $\sum_{i\geq 0}t_i=1$. In the set $EG$, $\langle g,t\rangle=
\langle g^\prime,t^\prime\rangle$ if and only if $t_i=t_i^\prime$ for each $i$
and $g_i=g^\prime_i$ for all $i$ with $t_i=t_i^\prime>0$. There is a right
action of $G$ on $EG$ given by $\langle g,t\rangle h=\langle
gh,t\rangle$ (cf. \cite{H}).
Observe that for any homomorphism $\rho:G\rightarrow G^\prime$, there is a 
canonical $\rho$-equivariant map $\rho_!:EG\rightarrow EG^\prime$ given by
$\rho_!\langle g,t\rangle=\langle\rho(g),t\rangle$.

\vspace{1.5mm}

The equivariant topology of a $G$-space $(Y,G)$ is studied through the
associated Borel space $Y_G$. For example, the equivariant cohomology of
$(Y,G)$ is defined to be the ordinary cohomology of $Y_G$. 
In this section, we shall establish a natural isomorphism between the
homotopy groups of the global orbispace $X:=Y/G$ canonically defined by
$(Y,G)$ and the (ordinary) homotopy groups of the Borel space $Y_G$. Suppose
$(f,\rho):(Y,G)\rightarrow (Y^\prime,G^\prime)$ is a pair of maps where
$f$ is $\rho$-equivariant. Then there is a fiber preserving map $(f,\rho)_!:
Y_G\rightarrow Y^\prime_{G^\prime}$ induced from $(\rho_!,f):EG\times Y
\rightarrow EG^\prime\times Y^\prime$. We denote the corresponding
homomorphisms between homotopy groups by $(f,\rho)_\ast:\pi_k(Y_G,\ast)
\rightarrow \pi_k(Y^\prime_{G^\prime},\ast)$. On the other hand, the pair
$(f,\rho)$ determines a morphism between the corresponding global orbispaces
$(\tilde{f},\tilde{\rho}):X\rightarrow X^\prime$. We set $(f,\rho)_{\#}:
\pi_k(X,\ast)\rightarrow \pi_k(X^\prime,\ast)$ for the induced homomorphisms.

\vspace{1.5mm}

\noindent{\bf Theorem 3.4.1:}
{\em For any $k\geq 1$, there is a natural isomorphism $\theta_k^X:
\pi_k(X,\ast)\rightarrow \pi_k(Y_G,\ast)$ satisfying the commutativity
condition: 
$$
\theta_k^{X^\prime}\circ (f,\rho)_{\#}
=(f,\rho)_\ast\circ \theta_k^X \leqno (3.4.1)
$$ 
for any equivariant map $(f,\rho):(Y,G)
\rightarrow (Y^\prime,G^\prime)$.
}

\hfill $\Box$

We start with a different characterization of the based loop space 
$\Omega(X,\underline{o})$ of a global orbispace $X$ arising from a $G$-space 
$(Y,G)$. The base-point structure $\underline{o}=(o,U_o,\hat{o})$ where we 
identify $\hat{o}$ as a point in $Y$. We denote by $P(Y,\hat{o})$ the space 
of all paths $\gamma:[0,1]\rightarrow Y$ such that $\gamma(0)=\hat{o}$, and
by $P(Y,G,\hat{o})$ the subspace of $P(Y,\hat{o})\times G$ which consists of
$(\gamma,g)$ satisfying $\gamma(1)=g\cdot\hat{o}$. The space $P(Y,\hat{o})$ is
given the compact-open topology, and $P(Y,G,\hat{o})$ given the relative
topology as a subspace of $P(Y,\hat{o})\times G$. The space $P(Y,G,\hat{o})$ 
is naturally a based space with base point $\tilde{o}=(\gamma_o,1_{G})$ 
where $\gamma_o$ is the constant path into $\hat{o}$. Any equivariant map 
$(f,\rho):(Y,G)\rightarrow (Y^\prime,G^\prime)$ induces a continuous map
$P(f,\rho):P(Y,G,\hat{o})\rightarrow P(Y^\prime,G^\prime,f(\hat{o}))$ given
by $(\gamma,g)\mapsto (f\circ\gamma,\rho(g))$.

\vspace{1.5mm}

\noindent{\bf Lemma 3.4.2:}
{\em There is a natural homeomorphism $\phi_X:(\Omega(X,\underline{o}),
\tilde{o})\cong (P(Y,G,\hat{o}),\tilde{o})$ such that for any equivariant
map $(f,\rho):(Y,G)\rightarrow (Y^\prime,G^\prime)$, we have commutativity:
$$
P(f,\rho)\circ\phi_X=\phi_{X^\prime}\circ \Omega((\tilde{f},\tilde{\rho})).
\leqno (3.4.2)
$$
where
$\Omega((\tilde{f},\tilde{\rho})):\Omega(X,\underline{o})\rightarrow 
\Omega(X^\prime,\underline{o}^\prime)$ is induced from 
$(\tilde{f},\tilde{\rho}):X\rightarrow X^\prime$. 
}

\noindent{\bf Proof:}
We shall first construct $\phi_X:(\Omega(X,\underline{o}),\tilde{o})
\rightarrow (P(Y,G,\hat{o}),\tilde{o})$ as a continuous map.
Without loss of generality, we assume that in the base-point structure 
$\underline{o}=(o,U_o,\hat{o})$, $U_o$ is the connected component of $X$
containing $o$. Suppose $\tilde{\gamma}\in\Omega(X,\underline{o})$ is 
represented by $(\{I_i\},\{U_i\},\{\gamma_i\},\{\rho_{ji}\})$, 
$i=0,1,\cdots,n$. Since the image of $\tilde{\gamma}$ in $X$ lies in $U_o$, 
each $U_i$ is an open subset of $U_o$, each $\widehat{U_i}$ is an open subset
of $\widehat{U_o}\subset Y$, and each $\rho_{ji}\in Tran(U_i,U_j)$
can be identified with an element of $G_{U_o}$. We define a path $\gamma:
[0,1]\rightarrow Y$ as follows: $I_0$ is decomposed as $I_{0,-}\cup I_{0,+}$
along $\ast\in S^1$. On $I_{0,+}$, $\gamma=\gamma_0$,
on $I_1$, $\gamma=\rho_{10}^{-1}\circ \gamma_1, \cdots,$ on $I_i$, $\gamma=
\rho_{10}^{-1}\circ\cdots\circ \rho_{i(i-1)}^{-1}\circ \gamma_i, \cdots,$
and finally on $I_{0,-}$, $\gamma=
\rho_{10}^{-1}\circ\cdots\circ \rho_{0n}^{-1}\circ \gamma_0$. We define $g:=
\rho_{10}^{-1}\circ\cdots\circ \rho_{0n}^{-1}\in G_{U_o}$. Clearly we have
$(\gamma,g)\in P(Y,G,\hat{o})$. We define $\phi_X$ by setting 
$\phi_X(\tilde{\gamma}):=(\gamma,g)$, which is independent on the choice of
representatives of the based loop $\tilde{\gamma}$. The map $\phi_X$ is
continuous by checking with canonical neighborhoods of a based loop. 
Furthermore, $(3.4.2)$ is clear from the construction, and $\phi_X$ is a
base-point preserving map.

It remains to show that $\phi_X$ is a homeomorphism. We construct a map
$\psi_X:P(Y,G,\hat{o})\rightarrow \Omega(X,\underline{o})$ as follows.
Given $(\gamma,g)\in P(Y,G,\hat{o})$, it is easily seen that $\gamma(t)\in
\widehat{U_o}$ for all $t\in [0,1]$, and $g\in G_{U_o}$ by connectivity.
We associate $(\gamma,g)$ with a canonical based path in $P(X,\underline{o},
\underline{o})$, which is the equivalence class of the following system
$(\{[0,\frac{2}{3}),(\frac{1}{3},1]\},\{U_o,U_o\},\{\gamma_0,\gamma_1\},
\rho_{10})$ where $\gamma_0=\gamma|_{[0,\frac{2}{3})}$, 
$\gamma_1=g^{-1}\circ\gamma|_{(\frac{1}{3},1]}$, and $\rho_{10}=g^{-1}$.
Now we identify $P(X,\underline{o},\underline{o})$ with 
$\Omega(X,\underline{p})$, and the map $\psi_X$ is resulted. Clearly
$\psi_X$ is continuous and $\phi_X\circ\psi_X=Id$,
$\psi_X\circ\phi_X=Id$.
Hence $\phi_X$ is a homeomorphism.

\hfill $\Box$

Now we consider the continuous map $\pi:(P(Y,G,\hat{o}),\tilde{o})
\rightarrow (G,1_G)$ by sending $(\gamma,g)$ to $g$. The fiber at $1_G$ 
is identified with the based loop space $\Omega(Y,\hat{o})$ of $Y$, via the 
embedding $\Omega(Y,\hat{o})\hookrightarrow P(Y,G,\hat{o})$ sending $\gamma$ 
to $(\gamma,1_G)$.

\vspace{1.5mm}

\noindent{\bf Lemma 3.4.3:}
{\em The continuous map $\pi:(P(Y,G,\hat{o}),\tilde{o})\rightarrow (G,1_G)$ 
is a fibration. As a consequence, we have a long exact sequence
$$
\cdots
\rightarrow \pi_k(G,1_{G})\stackrel{\partial}{\rightarrow}\pi_k(Y,\hat{o})
\stackrel{i_\ast}{\rightarrow}\pi_k(X,\underline{o})\stackrel{\pi_\ast}
{\rightarrow}\pi_{k-1}(G,1_{G})\rightarrow\cdots \leqno (3.4.3)
$$
for $k\geq 1$. 
}

\vspace{2mm}

\noindent{\bf Proof:}
In order to show that $\pi:(P(Y,G,\hat{o}),\tilde{o})\rightarrow 
(G,1_G)$ is a fibration, we only need to show that for any continuous maps
$f:T\rightarrow P(Y,G,\hat{o})$ and $H:T\times [0,1]\rightarrow G$
such that $\pi\circ f=H(\cdot,0)$, there is a continuous map $F:T\times [0,1]
\rightarrow P(Y,G,\hat{o})$ such that $\pi\circ F=H$ and 
$F(\cdot,0)=f$. Set $f(t)=(\gamma_{t},g_t)$. Then $\pi\circ f=H(\cdot,0)$
implies that $g_t=H(t,0)$ for all $t\in T$. We define $F:T\times [0,1]
\rightarrow P(Y,G,\hat{o})$ by 
$F(t,s)=(\gamma_{t,s},H(t,s))$ where $\gamma_{t,s}$ is defined by
$$
\gamma_{t,s}(\tau)=\left\{\begin{array}{cc}
\gamma_{t}((1+s)\tau) & \tau\leq (1+s)^{-1}\\
H(t,(1+s)\tau-1)\hat{o} & (1+s)^{-1}\leq\tau\leq 1.
\end{array} \right.
$$

The fibration $\pi$ induces a long exact sequence
$$
\cdots\rightarrow \pi_k(G,1_{G})\stackrel{\partial}{\rightarrow}\pi_{k-1}
(\Omega(Y,\hat{o}),\hat{o})\stackrel{i_\ast}{\rightarrow}
\pi_{k-1}(P(Y,G,\hat{o}),\tilde{o})\stackrel{\pi_\ast}{\rightarrow}
\pi_{k-1}(G,1_{G})\rightarrow\cdots,
$$
which gives rise to $(3.4.3)$ after identifying $(\Omega(X,\underline{o}),
\tilde{o})$ with $(P(Y,G,\hat{o}),\tilde{o})$ by the natural homeomorphism 
$\phi_X$.

\hfill $\Box$

We now look at the Borel space $Y_G=EG\times_G Y$. We fix a base point
$\ast\in EG$ by $\ast=(11_{G},0,\cdots)$, which then gives rise to a base 
point $[(\ast,\hat{o})]\in Y_G$, and a base point $\ast=pr(\ast)\in BG$. 
For simplicity we denote $[(\ast,\hat{o})]$ by $\check{o}$.
As a fiber bundle $\pi:Y_G\rightarrow BG$ over the classifying space $BG$ of 
$G$, with fiber $Y$, we have a long exact sequence
$$
\cdots
\rightarrow \pi_k(BG,\ast)\stackrel{\partial}{\rightarrow}\pi_{k-1}(Y,\hat{o})
\stackrel{i_\ast}{\rightarrow}\pi_{k-1}(Y_G,\check{o})\stackrel{\pi_\ast}
{\rightarrow}\pi_{k-1}(BG,\ast)\rightarrow\cdots\leqno (3.4.4)
$$
for $k\geq 1$. 

On the other hand, the universal principal $G$-bundle 
$pr:EG\rightarrow BG$ induces  for each $k\geq 1$ an isomorphism 
$\delta: \pi_k(BG,\ast)\rightarrow \pi_{k-1}(G,1_{G})$. 
A geometric representation of the isomorphism $\delta$ can be obtained as
follows. There is a continuous map $\Theta:G\rightarrow P(EG,\ast)$ where
$P(EG,\ast)$ is the space of paths $\gamma$ in $EG$ with $\gamma(0)=\ast$:
for each $g\in G$, the path $\Theta(g):[0,1]\rightarrow EG$ is given by
$$
\Theta(g)(t)=\left\{\begin{array}{cc}
((1-2t)1_{G},2t1_{G},0,\cdots), & 0\leq t\leq \frac{1}{2}\\
((2t-1)g,(2-2t)1_{G},0,\cdots), & \frac{1}{2}\leq t\leq 1.
\end{array} \right.
$$
Then for each class $[u]\in\pi_{k-1}(G,1_{G})$, represented by $u:
(S^{k-1},\ast)\rightarrow (G,1_{G})$, the class $\delta^{-1}([u])$ in
$\pi_k(BG,\ast)$ is represented by $pr\circ\Theta\circ u:(S^{k-1},\ast)
\rightarrow (\Omega(BG,\ast),\ast)$. 

Finally, we observe that the map $\Theta$ is natural with respect to 
homomorphisms $\rho:G\rightarrow G^\prime$, i.e., $\Theta(\rho(g))=
\rho_!\circ\Theta(g)$.

\vspace{2mm}

\noindent{\bf Proof of Theorem 3.4.1:}

Let $\theta_\ast: (P(Y,G,\hat{o}),\tilde{o})\rightarrow 
(\Omega(Y_G,\check{o}),\check{o})$ be the continuous map defined by
$(\gamma,g)\mapsto [(\Theta(g^{-1}),\gamma)]$, and let 
$\theta:\pi_k(X,\underline{o})\rightarrow\pi_k(Y_G,\check{o})$ be the 
homomorphism induced by $\theta_\ast$. Set $\Delta:\pi_{k-1}(G,1_{G})
\rightarrow\pi_k(BG,\ast)$ to be the composition $\delta^{-1}\circ\nu$
where $\nu: \pi_{k-1}(G,1_{G})\rightarrow \pi_{k-1}(G,1_{G})$ is induced by
$g\mapsto g^{-1}$.

We will show that the long exact sequences $(3.4.3)$ and $(3.4.4)$ fit into 
the following commutative diagram:
$$
\begin{array}{ccccccccc}
\cdots
\rightarrow & \pi_{k}(G,1_{G}) & \stackrel{\partial}{\rightarrow} & 
\pi_{k}(Y,\hat{o})
\stackrel{i_\ast}{\rightarrow} & \pi_{k}(X,\underline{o}) & \stackrel{\pi_\ast}
{\rightarrow} & \pi_{k-1}(G,1_{G}) & \rightarrow \cdots\\
            & \Delta\downarrow   &            & Id\downarrow
                               & \theta \downarrow &                   
              & \Delta\downarrow &             \\
\cdots\rightarrow & \pi_{k+1}(BG,\ast) & \stackrel{\partial}{\rightarrow} & 
\pi_k(Y,\hat{o})\stackrel{i_\ast}{\rightarrow} & \pi_k(Y_G,\check{o}) 
& \stackrel{\pi_\ast}{\rightarrow} & \pi_{k}(BG,\ast) & \rightarrow \cdots
\end{array}
$$
where $k\geq 1$, which implies that $\theta$ is an isomorphism. The
commutativity $(3.4.1)$ also follows since every thing involved is natural
with respect to equivariant maps $(Y,G)\rightarrow (Y^\prime,G^\prime)$.

The commutativities $\Delta\circ\pi_\ast=\pi_\ast\circ\theta$ and
$\theta\circ i_\ast=i_\ast\circ Id$ follows immediately from the definitions.
As for $Id\circ\partial=\partial\circ\Delta$, we first observe that for a
class $[u]\in \pi_k(G,1_{G})$ represented by a map $u:(S^k,\ast)
\rightarrow (G,1_{G})$, the image $\partial[u]$ in $\pi_k(Y,\hat{o})$ is
represented by the map $x\mapsto u(x)\cdot\hat{o}, \forall x\in S^k$. On the
other hand, $\Delta [u]$ is represented by the map from $(S^k,\ast)
\rightarrow (\Omega(BG,\ast),\ast)$ given by 
$x\mapsto pr\circ\Theta(u(x)^{-1})$, whose image under $\partial$ is 
represented by the map $x\mapsto u(x)\cdot\hat{o}$. Hence 
$Id\circ\partial=\partial\circ\Delta$, and the theorem is proved.

\hfill $\Box$

\noindent{\bf Remark 3.4.4:}\hspace{2mm}
Let $X_G$ be the orbispace which consists of one point with a trivial action
of $G$ on it. Then by Theorem 3.4.1, we have $\pi_k(X_G,\ast)=\pi_k(BG,\ast)$.
Let $H\subset G$ be a subgroup. Then it is easily seen that $X_H\rightarrow
X_G$ is a pseudo-embedding if and only if $G\rightarrow G/H$ is a weak
fibration. The based relative loop space is identified with $G/H$. Hence this
example shows that in order to have the homotopy exact sequence $(3.3.4)$
associated to a pseudo-embedding, it is necessary to incorporate the weak
fibration condition on the maps $G_{U_\alpha}\rightarrow 
G_{U_\alpha}/\rho_\alpha(G_{V_\alpha})$
into the definition.

\subsection{Free loop space and twisted sectors}

\hspace{5mm}
Let $X$ be an orbispace. The {\it free loop space} of $X$, written as
$\L X$, is the space of all morphisms from $S^1$ into $X$ equipped
with the compact-open topology, where $S^1$ is regarded as an
orbispace with trivial orbispace structure. The group of
self-homeomorphisms of $S^1$, $Homeo(S^1)$, acts on the free loop
space $\L X$ by re-parameterization. In particular, there is a
canonical $S^1$ action on $\L X$ by rotating the domain of each free
loop. Before we identify the fixed-point set of the actions, let us
first give a more convenient description of $\L X$ when $X$ is a
global orbispace defined by a $G$-space $(Y,G)$. 

Let $(Y,G)$ be a $G$-space. We denote by $P(Y)$ the space of all paths
$\gamma:[0,1]\rightarrow Y$ equipped with the compact-open
topology. Let $P(Y,G)$ be the subspace of $P(Y)\times G$ consisting of
pairs $(\gamma,g)$ satisfying $\gamma(1)=g\cdot\gamma(0)$. There is a
canonical $G$-action on $P(Y,G)$ given by $h\cdot
(\gamma,g)=(h\circ\gamma, hgh^{-1})$. Denote the quotient space by
$P(Y,G)/G$. 

\vspace{1.5mm}

\noindent{\bf Lemma 3.5.1:}
{\em Let $X$ be a global orbispace defined by a $G$-space
$(Y,G)$. Then the free loop space $\L X$ is naturally homeomorphic to
$P(Y,G)/G$. 
}

\vspace{2mm}

\noindent{\bf Proof:} We define a map $\phi:\L X\rightarrow P(Y,G)/G$
as follows. Let $\tilde{\gamma}$ be a free loop in $X$. We represent
$\tilde{\gamma}$ by a system $\sigma$. Then as in the proof of Lemma
3.4.2, $\sigma$ determines an element $(\gamma,g)\in P(Y,G)$. We
define $\phi(\tilde{\gamma})=[(\gamma,g)]$ where $[(\gamma,g)]$ is the
image of $(\gamma,g)$ in $P(Y,G)/G$. In order to verify that $\phi$ is
well-defined, we observe that if we think of $\sigma$ is a based
system, then $(\gamma,g)$ does not depend on the equivalent class of
$\sigma$. On the other hand, if we change the base-point structure,
then it amounts to change $(\gamma,g)$ to $(h\circ\gamma, hgh^{-1})$
for some $h\in G$. Hence $\phi$ is well-defined. 

Similarly, $\phi$ has an inverse $\psi:P(Y,G)/G\rightarrow \L X$
defined as follows. As in the proof of Lemma 3.4.2, each element
$(\gamma,g)$ canonically determines a based loop in $X$. For any class
$[(\gamma,g)]\in P(Y,G)/G$, we choose a representative $(\gamma,g)$,
and define $\psi([(\gamma,g)])$ to be the free loop corresponding to
the based loop determined canonically by $(\gamma,g)$. It is clear
that $\psi$ is well-defined, and is the inverse of $\phi$. 

Finally, we observe that both $\phi$ and $\psi$ are continuous with
respect to the compact-open topology. 

\hfill $\Box$

\noindent{\bf Remark 3.5.2:}\hspace{2mm}
When $Y$ is a smooth manifold, $G$ is a finite group, the $G$-space
$(P(Y,G),G)$ is precisely the space of strings satisfying twisted
boundary conditions in the physics literature (cf. \cite{DHVW}). 

\hfill $\Box$

Now we are ready to identify the fixed-point set of the canonical
$S^1$-action on $\L X$. First, we extend the space $\widetilde{X}$
defined by $(1.5)$ to any orbispace. For any orbispace $X$, we define
$$
\widetilde{X}:=\{(p,(g)_{G_p})|p\in X, g\in G_p\}, \leqno (3.5.1)
$$
where $G_p$ is the isotropy group of $p$ (cf. Definition 2.1.2), and
$(g)_{G_p}$ is the conjugacy class of $g$ in $G_p$.  We remark that
when $X$ is an almost complex orbifold (here the local group actions
are allowed to be non-effective), a new cup product was constructed in
\cite{CR1} on the cohomology groups of $\widetilde{X}$ with suitable
degree shifting. The new cup product was defined using the degree
zero, genus zero ``Gromov-Witten invariants'' of $X$. When $X$ is
defined from the space of one point with a trivial action by a finite
group $G$, the new cohomology ring was shown to be
isomorphic to the center of the group algebra of $G$
(cf. \cite{CR1}). The components of $\widetilde{X}$ with $(g)\neq (1)$
are called {\it twisted sectors} in \cite{CR1}.

\vspace{2mm}

\noindent{\bf Proposition 3.5.3:}
{\em Let $X$ be an orbispace such that for any basic open set $U$ of
$X$, the group $G_U$ is finite and the space $\widehat{U}$ is a $T_1$
space. Then the fixed-point set of the canonical $S^1$ action on $\L
X$ can be identified with $\widetilde{X}$ as defined in $(3.5.1)$.
}

\vspace{2mm}

\noindent{\bf Proof:} Each free loop $\tilde{\gamma}$ in $X$
determines a free loop $\gamma$ in $X_{top}$. If $\tilde{\gamma}$ is
invariant under the canonical $S^1$ action on $\L X$, so is $\gamma$
under the canonical $S^1$ action on the free loop space of
$X_{top}$. In particular, if $\tilde{\gamma}$ is a fixed point,
$\gamma$ is a constant loop in $X_{top}$. Let $U$ be a basic open set
of $X$ containing $\gamma$. By Lemma 3.5.1, $\tilde{\gamma}$ can be
identified as an element $(\gamma^\prime, g)$ in
$P(\widehat{U},G_U)/G_U$. In this context, it is easy to see that
$\tilde{\gamma}$ being a fixed point implies that there is a $g_t\in
G_U$ for any $t\in [0,1]$ such that
$\gamma^\prime(t)=g_t\circ\gamma^\prime(0)$. Since $G_U$ is a finite
group by assumption, the image of $\gamma^\prime$ in $\widehat{U}$ is
a finite subset. Moreover, $\widehat{U}$ is a $T_1$ space, so that any
finite subset must be closed. It follows that $\gamma^\prime$ must be
a constant map. On the other hand, any element $(\gamma^\prime,g)$ in 
$P(\widehat{U},G_U)/G_U$ where $\gamma^\prime$ is a constant map
defines a fixed point of the $S^1$ action on $\L X$. From here, it is
easily seen that the fixed-point set is identified with $\widetilde{X}$.

\hfill $\Box$

We remark that orbifolds obviously satisfy the hypothesis in
Proposition 3.5.3. 

\vspace{1.5mm}

Lemma 3.5.1 shows that the free loop space of a global orbispace is
the orbit space of a $G$-space. Although for a general orbispace, it
may be difficult to show any orbispace structure on the free loop
space, we shall next prove that the space of free smooth loops in an
orbifold has a natural pre-Hilbert orbifold structure. 

Let $X$ be an orbifold. For any representing system
$\sigma=(\{I_i\},\{U_i\},\{\gamma_i\},\{\xi_{ji}\})$ of a smooth free
loop $\tilde{\gamma}$ in $X$, there is a pull-back bundle $E_\sigma$
over $S^1$, which is defined as 
$$
E_\sigma:=\cup_i\gamma_i^\ast TV_i/\sim, \leqno (3.5.2)
$$ 
where each $U_i$ is uniformized by $(V_i,G_i,\pi_i)$, and   
the identification $\sim$ is made by the transition maps
$\{\xi_{ji}\}$. There is a natural bundle morphism
$\bar{\sigma}:E_\sigma\rightarrow TX$ covering the morphism
$\tilde{\gamma}$. If $\tau$ is another representative of
$\tilde{\gamma}$, then there is an isomorphism
$\theta_{\tau\sigma}:E_\sigma\rightarrow E_\tau$ such that
$\bar{\tau}\circ \theta_{\tau\sigma}=\bar{\sigma}$. We define
$G_\sigma$ to be the group of automorphisms
$\theta:E_\sigma\rightarrow E_\sigma$ such that
$\bar{\sigma}\circ\theta=\bar{\sigma}$. Clearly the isomorphism class
of $G_\sigma$ depends only on the free loop $\tilde{\gamma}$. We
denote the abstract group by $G_{\tilde{\gamma}}$, and call it the
{\it isotropy group} of $\tilde{\gamma}$.

\vspace{1.5mm}

\noindent{\bf Lemma 3.5.4:} {\em The isotropy group
$G_{\tilde{\gamma}}$ of a free loop $\tilde{\gamma}$ is finite. In
fact, if let $\gamma$ be the loop in $X_{top}$ induced by
$\tilde{\gamma}$, then for any $p\in Im\gamma$, there is a
monomorphism $G_{\tilde{\gamma}}\rightarrow G_p$. 
}

\vspace{2mm}

\noindent{\bf Proof:} Let
$\sigma=(\{I_i\},\{U_i\},\{\gamma_i\},\{\xi_{ji}\})$ be a
representative of $\tilde{\gamma}$. Then any element $\delta$ of
$G_\sigma$ is realized by a collection of $\{\delta_i\}$, where each
$\delta_i$ is an automorphism of the uniformizing system of $U_i$
which fixes each point in the image of $\gamma_i$ and satisfies the
compatibility condition
$\xi_{ji}=\delta_j\circ\xi_{ji}\circ\delta_i^{-1}$. For any $p\in
Im\gamma$, suppose $p$ is in $U_i$ for some $i$. Then there is a
homomorphism $G_\sigma\rightarrow G_p$ given by $\delta\mapsto
\delta_i$. It is a monomorphism by the compatibility condition 
$\xi_{ji}=\delta_j\circ\xi_{ji}\circ\delta_i^{-1}$. As a consequence,
$G_\sigma$ hence $G_{\tilde{\gamma}}$ is finite since each $G_p$ is finite.

\hfill $\Box$

\noindent{\bf Theorem 3.5.5:} {\em The space of free smooth loops of
an orbifold has a pre-Hilbert orbifold structure such that for any
free loop $\tilde{\gamma}$, a neighborhood of $\tilde{\gamma}$ is
uniformized by $(\O(\Gamma(E_\sigma)),G_\sigma)$, where $\sigma$ is
any representing system of $\tilde{\gamma}$, $\O(\Gamma(E_\sigma))$ is
a $G_\sigma$-invariant neighborhood of the zero section in the space
of smooth sections of $E_\sigma$. 
}

\vspace{2mm}

\noindent{\bf Proof:}
We fix a Riemannian metric on $X$. We assume that each basic open set
$U$ of $X$ is a convex geodesic neighborhood (cf. Example 2.1.3 b) so
that its uniformizing system $(V,G,\pi)$ is a geodesically convex
ball. The exponential map is defined on each $V$ which is
$G$-equivariant.

Let $\sigma=(\{I_i\},\{U_i\},\{\gamma_{0,i}\},\{\xi_{ji}\})$ be a
representing system of a free loop $\tilde{\gamma}_0$. Without loss of
generality, we may assume that $\gamma_{0,i}$ is defined over the
closure of $I_i$. The bundle $E_\sigma$ over $S^1$ defined by
$(3.5.2)$ inherits a natural metric from $TX$, so that we can speak of 
$C^\infty$ norm $|\cdot|$ of its smooth sections. It is easy to see that for 
sufficiently small $r>0$,
any smooth section $s=\{s_i\}$ of $E_\sigma$ with $|s|<r$ defines a
system $(\{I_i\},\{U_i\},\{\gamma_{i}\},\{\xi_{ji}\})$ where
$\gamma_i=exp_{\gamma_{0,i}} s_i$. We set 
$$
\O_r(\Gamma(E_\sigma))=\{s\in
\Gamma(E_\sigma)||s|<r\}. \leqno (3.5.3)
$$
There is an induced effective $G_\sigma$ action on
$\O_r(\Gamma(E_\sigma))$. Through exponential maps, we can define a
$G_\sigma$-equivariant continuous map
$\pi_\sigma:\O_r(\Gamma(E_\sigma))\rightarrow \L X$, which is easily
seen to be surjective onto a neighborhood of $\tilde{\gamma}$ in $\L
X$. It remains to show that we can choose $r$ sufficiently small, so
that for any two smooth sections $s_1,s_2\in \O_r(\Gamma(E_\sigma))$,
if $\pi_\sigma(s_1)=\pi_\sigma(s_2)$, there is a $\delta\in G_\sigma$
such that $s_2=\delta\circ s_1$.  

Let $\pi_\sigma(s_1),\pi_\sigma(s_2)$ be given by
$\sigma_1=(\{I_i\},\{U_i\},\{\gamma_{i}^1\},\{\xi_{ji}\})$ and 
$\sigma_2=(\{I_i\},\{U_i\},\{\gamma_{i}^2\},\{\xi_{ji}\})$
respectively. The assumption $\pi_\sigma(s_1)=\pi_\sigma(s_2)$ means
that there is a common refinement of $\sigma_1$ and $\sigma_2$. This
gives rise to a collection of $\{\delta_i\}$, where each $\delta_i$ is
an automorphism of the uniformizing system of $U_i$ satisfying the
compatibility condition
$\xi_{ji}=\delta_j\circ\xi_{ji}\circ\delta_i^{-1}$, such that
$\gamma_i^2=\delta_i\circ\gamma_i^1$. We will show that when $r$ is 
chosen sufficiently small, each $\delta_i$ will leave the image of
$\gamma_{0,i}$ fixed. The collection $\{\delta_i\}$ then defines an
element $\delta\in G_\sigma$ which satisfies $s_2=\delta\circ s_1$.

We need a digression on the structure of the set of singular points $\Sigma X$
of the orbifold $X$. Recall the following resolution of $\Sigma X$
from \cite{Ka}.
Let $(1)=(H_p^0), (H_p^1),..., (H_p^{n_p})$ be all the orbit types of a 
geodesic uniformizing system $(V_p,G_p,\pi_p)$ at $p$. 
For $q\in U_p=\pi_p(V_p)$, we may take $U_q$ small enough so that
$U_q\subset U_p$. Then any injection $\phi: V_q\rightarrow V_p$ induces a 
unique homomorphism $\lambda_\phi:G_q\rightarrow G_q$, which gives a 
correspondence $(H_q^i)\rightarrow \lambda_\phi(H_q^i)=(H_p^j)$. This 
correspondence is independent of the choice of $\phi$. Consider the set of 
pairs:
$$
\widetilde{\Sigma U_p}=\{(q,(H_q^i))|q\in \Sigma U_p, i\neq 0\}.
\leqno (3.5.4)
$$
Take one representative $H_q^i\in (H_q^i)$. Then the pair $(q,(H_q^i))$
determines exactly one orbit $[\tilde{q}]$ in the fixed point set 
$V_p^{H_p^j}$ by the action of the normalizor $N_{G_p}(H_p^j)$, where
$\tilde{q}=\phi(q)$, $(H_p^j)=\lambda_\phi(H_q^i)$. The correspondence 
$(q,(H_q^i))\rightarrow [\tilde{q}]$ gives a homeomorphism
$$
\widetilde{\Sigma U_p}\cong \coprod_{j=1}^{n_p}V_p^{H_p^j}/N_{G_p}(H_p^j),
\hspace{2mm} \mbox {(disjoint union)}, \leqno (3.5.5)
$$
which gives $\widetilde{\Sigma X}=\{(p,(H_p^j))|p\in\Sigma X, j\neq 0\}$ 
an orbifold structure
$$
\{\pi_{p,j}:(V_p^{H_p^j}, N_{G_p}(H_p^j)/H_p^j)\rightarrow  
V_p^{H_p^j}/N_{G_p}(H_p^j): p\in X, j=1,\cdots,n_p.\}. \leqno (3.5.6)
$$
The canonical map $\pi:\widetilde{\Sigma X}\rightarrow \Sigma X$ defined by 
$(p,(H_p^j))\rightarrow p$ is surjective, which has a $C^\infty$ lifting 
$\tilde{\pi}$ given locally by embeddings $V_p^{H_p^j}\rightarrow V_p$.
$\pi:\widetilde{\Sigma X}\rightarrow \Sigma X$ is called the {\it canonical 
resolution} of the singular set $\Sigma X$.

A point $(p,(H_p^j))$ in $\widetilde{\Sigma X}$ is called
generic if $G_p=H_p^j$. The set $\widetilde{\Sigma X}_{gen}$ of all generic 
points is open dense in $\widetilde{\Sigma X}$, and the map 
$\pi|\widetilde{\Sigma X}_{gen}:\widetilde{\Sigma X}_{gen}
\rightarrow \Sigma X$ is bijective. Denote $X_{reg}$ the set of smooth
points in $X$. Then  we have a partition of $X$ 
into a disjoint union of smooth manifolds:
$$
X=X_{reg}\cup\widetilde{\Sigma X}_{gen}, \leqno (3.5.7)
$$
which is called the {\it canonical stratification} of $X$. Each
connected stratum $X_i$ of the canonical stratification is assigned
with a finite group $G_i$ such that if any point $p\in X$ lies in
$X_i$, then $G_i$ is isomorphic to $G_p$. We can introduce a partial
order $\prec$ amongst the strata $\{X_i\}$.  We say that $X_i\prec
X_j$ for $i\neq j$ if $X_j$ is contained in the closure of $X_i$ in
$X$. The condition $X_i\prec X_j$ implies that $G_i$ is a proper
subgroup of $G_j$. The canonical stratification $(3.5.7)$ is locally
finite. For any uniformizing system $(V,G,\pi)$, the inverse image of
each stratum $X_i$ consists of finitely many disjoint components, each
of which is mapped diffeomorphic to $\pi(V)\cap X_i$. It is clear that 
the canonical stratification $(3.5.7)$ induces a stratification of
$V$, which has the property that for any stratum $V_k$, if an element
$g\in G$ fixes a point $x\in V_k$, then $g$ fixes the entire stratum
$V_k$. End of digression.

We observe that for any point $p\in X$, there is a $r_p>0$ with the
following significance: for any uniformizing system $(V,G,\pi)$ such
that $p\in \pi(V)$, and any $\hat{p}\in V$ such that $\pi(\hat{p})=p$,
if $x,y\in V$ are connected to $\hat{p}$ by geodesics of length less
than $r_p$, and there is a $g\in G$ such that $g\cdot x=y$, then
$g\cdot \hat{p}=\hat{p}$. Let $(V_i,G_i,\pi_i)$ be the uniformizing
system of $U_i$ for each $i$. For each $U_i$, let
$\{V_{i,\alpha_j}|\alpha_j\in \Lambda(i)\}$ be the set of strata in
the stratification of $V_i$ induced by $(3.5.7)$ where each
$V_{i,\alpha_j}$ contains points in $\gamma_{0,i}(I_i)$. The local
finiteness of $(3.5.7)$ implies that each $\Lambda(i)$ is a finite
set. For each $\alpha_j\in\Lambda(i)$, choose a $t_{\alpha_j}\in I_i$
such that $\gamma_{0,i}(t_{\alpha_j})$ is contained in
$V_{i,\alpha_j}$. Denote
$p_{\alpha_j}=\pi_i(\gamma_{0,i}(t_{\alpha_j}))$ the image of
$\gamma_{0,i}(t_{\alpha_j})$ in $U_i$.  Let
$r_i=\min{r_{p_{\alpha_j}}}$ and $r=\min{r_i}$.  We take this $r$ for
the neighborhood $\O_r(\Gamma(E_\sigma))$
defined by $(3.5.3)$. If the two smooth
sections $s_1,s_2$ satisfying $\pi_\sigma(s_1)=\pi_\sigma(s_2)$ are in
$\O_r(\Gamma(E_\sigma))$, then each $\delta_i$ must fix
$\gamma_{0,i}(t_{\alpha_j})$ for every $\alpha_j\in \Lambda(i)$. 
The assumption that $\gamma_{0,i}(t_{\alpha_j})$ is contained in
$V_{i,\alpha_j}$ implies that $\delta_i$ fixes the entire stratum 
$V_{i,\alpha_j}$. Now we see that $\delta_i$ fixes every point in
$\gamma_{0,i}(I_i)$ because $\gamma_{0,i}(I_i)$ is contained in the union  
$\cup_{\alpha_j\in\Lambda(i)} V_{i,\alpha_j}$.  
Hence $\delta=\{\delta_i\}$ is an element of $G_\sigma$, and
the proof of Theorem 3.5.5 is completed.

\hfill $\Box$

\noindent{\bf Remark 3.5.6:}\hspace{2mm} 
Let $X$ be an orbifold. As a subspace of the free loop space $\L X$, the
space $\widetilde{X}$ defined in $(1.5)$ (or $(3.5.1)$) inherits a
canonical orbispace structure from the pre-Hilbert orbifold structure
of the free loop space $\L X$ specified in Theorem 3.5.5. For any
point $(p,(g)_{G_p})$ in $\widetilde{X}$, viewed as a constant free
loop, the corresponding pull-back bundle $E_\sigma$ defined in
$(3.5.2)$ is (with $\sigma=(\gamma_p,g)$, $\gamma_p$ being the
constant map into $p$)
$$
E_\sigma=(TV_p)_p\times [0,1]/\sim \leqno (3.5.8)
$$
where $\sim$ identifies $(TV_p)_p\times\{0\}$ with
$(TV_p)_p\times\{1\}$ by the action of $g\in G_p$. The isotropy 
group $G_\sigma$ is clearly the centralizer $C_{G_p}(g)$ of $g$ in
$G_p$. The constant loops in a neighborhood of $\sigma$ are given by
the ``constant'' sections of $E_\sigma$ of $(3.5.8)$, which can be
identified with a ball of the fixed point set $V_p^g$ of $g$ in
$V_p$. Hence the orbispace structure of $\widetilde{X}$ inherited from
the free loop space is given by
$$
\{(V_p^g,C_{G_p}(g))|p\in X, g\in G_p\}. \leqno (3.5.9)
$$
Note that the action of $C_{G_p}(g)$ on $V_p^g$ is not effective in
general, hence $\widetilde{X}$ equipped with $(3.5.9)$ is not an
orbifold in the classical sense. But we remark that this is precisely
the orbispace structure we need in the proof of associativity of the
new cup product in \cite{CR1}.

\hfill $\Box$

There are geometrical-analytical constructions on the free loop space 
of a smooth manifold, which uses the pre-Hilbert manifold structure of the free
loop space. Theorem 3.5.5 allows us to mimic these constructions on
the free loop space of an orbifold. For example, Bismut's work on
Atiyah-Singer index theorem (cf. \cite{Bi}), Floer homology
(cf. \cite{F}), and Witten's interpretation of elliptic genera
(cf. \cite{W, T, BT}). There is also an interesting
construction of Chas and Sullivan on the free loop space of a smooth
manifold, the so-called string topology (cf. \cite{CS}). 

Let us specially look at Witten's interpretation of elliptic genera. If
successfully carried out in the orbifold case, it would give a general
definition of {\it orbifold elliptic genera} and prove the
corresponding rigidity theorem\footnote{relevant work has been
done recently, cf. \cite{Liu}}. There is the so-called elliptic
cohomology theory behind the scene (cf. \cite{Se2}), which has an
equivariant version (cf. \cite{De}). Thus it would be interesting to
study the corresponding generalized cohomology theories of orbispaces, 
which generalize the corresponding equivariant theories, and
ultimately establish a corresponding elliptic cohomology theory of 
orbispaces. 

\section{Basic Properties of Homotopy Groups}

\subsection{Orbispace covering}
\hspace{5mm}
This section is devoted to an orbispace covering theory, which generalizes 
all the basic results in the topological category. It also recovers Thurston's 
orbifold covering theory (cf. \cite{Th}). In particular, it is shown that the 
fundamental group of an orbifold defined in this paper coincides with 
Thurston's orbifold fundamental group.

In order to motivate the definition of orbispace covering, let us consider the 
following simple situation. Suppose $Y$ is a locally connected, path-connected 
topological space with a discrete group action of $G$. Denote the quotient 
space $Y/G$ by $X$. If $G$ acts on $Y$ properly discontinuously
without fixed points, we have the 
following exact sequence
$$
1\rightarrow\pi_1(Y,\ast)\rightarrow\pi_1(X,\ast)\rightarrow G\rightarrow 1. 
\leqno (4.1.1)
$$
For a general action of $G$, we still have an exact sequence of $(4.1.1)$, 
provided that $X$ is regarded as an orbispace canonically obtained from the 
$G$-space $(Y,G)$, and $(4.1.1)$ is just a special case of the long exact 
sequence $(3.4.3)$. The natural projection $Y\rightarrow X=Y/G$ is an example 
of covering map when $G$ acts properly discontinuously without fixed
points. This suggests that in the orbispace category, the natural
projection $Y\rightarrow X=Y/G$, for any action of a discrete group
$G$, should be considered as an orbispace covering morphism, since in 
a good covering theory, one should be able to associate an 
exact sequence such as $(4.1.1)$ to a covering map. 

\vspace{2mm}

\noindent{\bf Definition 4.1.1:}
{\em Let $\tilde{\pi}: Y\rightarrow X$ be a morphism between two orbispaces, 
and $\pi: Y_{top}\rightarrow X_{top}$ be the induced continuous map. We call 
$\tilde{\pi}$ an orbispace covering morphism, if there exists a representing 
system $(\{V_\alpha\},\{U_\alpha\},\{\pi_\alpha\},\{\rho_{\beta\alpha}\})$ of 
$\tilde{\pi}$ such that 
\begin{itemize}
\item [{a)}] for any $U\in\{U_\alpha\}$, the set of connected components of 
$\pi^{-1}(U)$ is contained in $\{V_\alpha\}$,
\item [{b)}] each homomorphism $\rho_\alpha: G_{V_\alpha}\rightarrow 
G_{U_\alpha}$ is monomorphic, and if $g,h\in G_{U_\alpha}$ is path-connected 
in $G_{U_\alpha}$, then $gh^{-1}$ lies in the image 
$\rho_{\alpha}(G_{V_\alpha})$, 
\item [{c)}] each $\pi_\alpha:\widehat{V_\alpha}\rightarrow\widehat{U_\alpha}$ 
is a $\rho_\alpha$-equivariant homeomorphism, and 
\item [{d)}] the map $\pi:Y_{top}\rightarrow X_{top}$ is surjective.
\end{itemize}

The orbispace $Y$ together with the orbispace covering morphism 
$\tilde{\pi}:Y\rightarrow X$ is called a covering space of $X$. Each basic 
open set $U\in\{U_\alpha\}$, which all together forms a cover of $X$ by $d)$, 
is called an elementary neighborhood of $X$ with respect to the 
covering morphism $\tilde{\pi}$.
}

\hfill $\Box$

\noindent{\bf Remark 4.1.2 a:}\hspace{2mm}
Any connected open subset of an elementary neighborhood is an elementary 
neighborhood in the following sense. Let $W\subset U\in\{U_\alpha\}$ be a 
connected open subset. Then each connected component of $\pi^{-1}(W)$, being a 
connected open subset of $\pi^{-1}(U)$, is a basic open set of $Y$ by $a)$ of 
Definition 4.1.1. We add to $\{U_\alpha\}$ all of the connected open subsets 
$W\subset U$ for each $U\in\{U_\alpha\}$, and to $\{V_\alpha\}$ all of the 
connected components of $\pi^{-1}(W)$ for each $U\in\{U_\alpha\}$ and each 
$W\subset U$. Let the resulting collections be $\{V_a\},\{U_a\}$, which 
clearly form a refinement of $\{V_\alpha\},\{U_\alpha\}$. Hence by 
Lemma 2.2.2, there is an induced system defined over $(\{V_a\},\{U_a\})$, 
which also represents $\tilde{\pi}$. One can easily verify that conditions 
$a),b),c)$ and $d)$ in Definition 4.1.1 are satisfied for the induced system.

\hfill $\Box$

\noindent{\bf Remark 4.1.2 b:}\hspace{2mm}
In the representing system 
$(\{V_\alpha\},\{U_\alpha\},\{\pi_\alpha\},\{\rho_{\beta\alpha}\})$ of the  
orbispace covering morphism $\tilde{\pi}$, 
each $\rho_{\beta\alpha}: Tran(V_\alpha,V_\beta)\rightarrow 
Tran(U_\alpha,U_\beta)$ is injective by the axioms $b),c)$ of
Definition 4.1.1. Moreover, for any two path connected transition maps
$\xi_1,\xi_2\in Tran(U_\alpha,U_\beta)$, there is a $g\in G_{V_\beta}$
such that $\xi_2=\rho_{\beta}(g)\circ\xi_1$ by the axiom $b)$ of 
Definition 4.1.1. 

\hfill $\Box$

\noindent{\bf Remark 4.1.2 c:}
Clearly the composition of two orbispace covering morphisms is again an
orbispace covering morphism.

\hfill $\Box$

\noindent{\bf Remark 4.1.2 d:}
There is a based version of orbispace covering morphisms. Let 
$(Y,\underline{q})$ and $(X,\underline{p})$ be two based orbispaces, where
$\underline{q}=(q,V_o,\hat{q})$ and $\underline{p}=(p,U_o,\hat{p})$.
Given any orbispace covering morphism $\tilde{\pi}:Y\rightarrow X$, we may
assume that $V_o$, $U_o$ are contained in $\{V_\alpha\}$ and $\{U_\alpha\}$ 
respectively, and the homeomorphism $\pi_o:\widehat{V_o}\rightarrow 
\widehat{U_o}$ sends $\hat{q}$ to $\hat{p}$.

\hfill $\Box$ 

Let us look at the path-lifting property of orbispace covering first.
Let $\tilde{\pi}:(Y,\underline{q})\rightarrow (X,\underline{p})$ be a based 
orbispace covering morphism, with a representing system 
$\Pi=(\{V_\alpha\},\{U_\alpha\},\{\pi_\alpha\},\{\rho_{\beta\alpha}\})$ being 
fixed throughout. Let $\underline{p^\prime}=(p^\prime,U_{o^\prime},
\hat{p^\prime})$ be a base-point structure of $X$. Then for any connected
component $V_{o^\prime}$ of $\pi^{-1}(U_{o^\prime})$, the system $\Pi$ 
canonically determines a base-point structure $\underline{q^\prime}=(q^\prime,
V_{o^\prime},\hat{q^\prime})$ of $Y$ by setting $\hat{q^\prime}=
\pi^{-1}_{o^\prime}(\hat{p^\prime})$ and 
$q^\prime=\pi_{V_{o^\prime}}(\hat{q^\prime})\in Y$. For any $g\in 
G_{U_{o^\prime}}$, we denote the base-point structure $(p^\prime,U_{o^\prime},
g\cdot\hat{p^\prime})$ by $g\cdot\underline{p^\prime}$.

\vspace{2mm}

\noindent{\bf Lemma 4.1.3:}
{\em Suppose $\tilde{\pi}:(Y,\underline{q})\rightarrow (X,\underline{p})$ is 
a based orbispace covering morphism, with a fixed representing system $\Pi$.
For any based path $\tilde{\gamma}\in 
P(X,\underline{p},\underline{p^\prime})$, there is a unique component 
$V_{o^\prime(\tilde{\gamma})}$ of $\pi^{-1}(U_{o^\prime})$ and a unique coset 
$\delta(\tilde{\gamma})\in G_{U_{o^\prime}}/\rho_{o^\prime(\tilde{\gamma})}
(G_{V_{o^\prime(\tilde{\gamma})}})$ such that after fixing a representative
$g\in G_{U_{o^\prime}}$ of $\delta(\tilde{\gamma})$, there is a unique based 
path $\ell(\tilde{\gamma})\in P(Y,\underline{q},\underline{q^\prime})$, 
called the lifting of $\tilde{\gamma}$, such that $\tilde{\pi}\circ 
\ell(\tilde{\gamma})=g\circ\tilde{\gamma}$, where $g\circ\tilde{\gamma}$ 
stands for the image of $\tilde{\gamma}$ under the map 
$P(X,\underline{p},\underline{p^\prime})\rightarrow 
P(X,\underline{p},g\cdot\underline{p^\prime})$ induced by $g$, and 
$\underline{q^\prime}$ is the base-point structure canonically determined by
the system $\Pi$ for the component $V_{o^\prime(\tilde{\gamma})}$. For a 
continuous family of $\tilde{\gamma}_s\in 
P(X,\underline{p},\underline{p^\prime})$ parametrized by a locally 
path-connected space, the component 
$V_{o^\prime(\tilde{\gamma}_s)}$ and the coset $\delta(\tilde{\gamma}_s)$
are locally constant with respect to the parameter $s$, and the family of 
liftings $\ell(\tilde{\gamma}_s)$ is continuous in $s$.
}

\vspace{2mm}

\noindent{\bf Proof:}
Suppose $\tilde{\gamma}$ is represented by a system 
$(\{I_i\},\{U_i\},\{\gamma_i\},\{\xi_{ji}\})$, $i=0,1,\cdots,n$, 
where each $U_i$ is an 
elementary neighborhood of $X$ with respect to the orbispace covering 
morphism $\tilde{\pi}$. We define a system 
$(\{I_i\},\{V_i\},\{\ell(\gamma_i)\},\{\ell(\xi_{ji})\})$ from $[0,1]$
into $(Y,\underline{q})$ inductively as follows. For $i=0$, we simply let 
$\ell(\gamma_0):=\pi_o^{-1}\circ\gamma_0:I_0\rightarrow \widehat{V_o}$, which
clearly satisfies $\ell(\gamma_0)(0)=\hat{q}$. Now suppose that 
$\ell(\gamma_i)$ is done. We will define $\ell(\xi_{(i+1)i})$ and 
$\ell(\gamma_{i+1})$ as follows. First of all, there is a unique connected 
component $V_{i+1}$ of $\pi^{-1}(U_{i+1})$ with which  
$\pi_{V_i}(\ell(\gamma_i))$ intersects. Secondly, there is 
$g_{i+1}\in G_{U_{i+1}}$ such that $g_{i+1}\circ\xi_{(i+1)i}$ lies in the 
image of $\rho_{(i+1)i}: Tran(V_i,V_{i+1})\rightarrow Tran(U_i,U_{i+1})$. 
We change $\gamma_{i+1}$ and $\xi_{(i+1)i}$ to $g_{i+1}\circ\gamma_{i+1}$ and
$g_{i+1}\circ\xi_{(i+1)i}$ respectively, and change $\xi_{(i+2)(i+1)}$
to $\xi_{(i+2)(i+1)}\circ g_{i+1}^{-1}$ if $i+1<n$. This amounts to change the 
representing system of $\tilde{\gamma}$ to an isomorphic one. 
We define $\ell(\xi_{(i+1)i})=\rho_{(i+1)i}^{-1}(g_{i+1}\circ\xi_{(i+1)i})$, 
and $\ell(\gamma_{i+1})=\pi_{i+1}^{-1}(g_{i+1}\circ\gamma_{i+1})$. A different
choice of $g_{i+1}$ differs by a post-composition by $\rho_{i+1}(g)$
for some $g\in G_{V_{i+1}}$, which results in a post-composition of 
$\ell(\xi_{(i+1)i})$ and $\ell(\gamma_{i+1})$ by $g$. Hence by induction
a system $(\{I_i\},\{V_i\},\{\ell(\gamma_i)\},\{\ell(\xi_{ji})\})$ is obtained,
and a unique component $V_{o^\prime(\tilde{\gamma})}$ of 
$\pi^{-1}(U_{o^\prime})$ and a unique coset 
$\delta(\tilde{\gamma})\in G_{U_{o^\prime}}/\rho_{o^\prime(\tilde{\gamma})}
(G_{V_{o^\prime(\tilde{\gamma})}})$ is determined. One can easily
verify that $V_{o^\prime(\tilde{\gamma})}$ and $\delta(\tilde{\gamma})$ are 
independent of the choice of the representing
system $(\{I_i\},\{U_i\},\{\gamma_i\},\{\xi_{ji}\})$ we started with.
From the construction, it is also obvious that once a representative of the 
coset $\delta(\tilde{\gamma})$ is chosen, the system 
$(\{I_i\},\{V_i\},\{\ell(\gamma_i)\},\{\ell(\xi_{ji})\})$ determines a based
path which is independent of the system 
$(\{I_i\},\{U_i\},\{\gamma_i\},\{\xi_{ji}\})$ we started with.
We define $\ell(\tilde{\gamma})$ to be the based path as the equivalence 
class of $(\{I_i\},\{V_i\},\{\ell(\gamma_i)\},\{\ell(\xi_{ji})\})$, which is
clearly unique.

For the parametrized case of $\tilde{\gamma}_s$, the component $V_{i+1}$ 
and the coset of $g_{i+1}$ is locally constant in $s$ in each induction step.
The former is obvious and the latter is by condition $b)$ of Definition 4.1.1.
Hence the component $V_{o^\prime(\tilde{\gamma}_s)}$ and the coset 
$\delta(\tilde{\gamma}_s)$ are locally constant with respect to the parameter 
$s$, and the family of liftings $\ell(\tilde{\gamma}_s)$ is continuous in $s$.

\hfill $\Box$

This lifting property of based paths has the following consequences as in the 
case of topological spaces.

\vspace{1.5mm}

\noindent{\bf Theorem 4.1.4:}
{\em Let $\tilde{\pi}:(Y,\underline{q})\rightarrow (X,\underline{p})$ be 
a based orbispace covering morphism. Then the induced homomorphism $\pi_{\#}:
\pi_k(Y,\underline{q})\rightarrow \pi_k(X,\underline{p})$ is an isomorphism 
for all $k\geq 2$ and a monomorphism for $k=1$.
}

\vspace{2mm}

\noindent{\bf Proof:}
Let $u:(S^{k-1},\ast)\rightarrow (\Omega(Y,\underline{q}),\tilde{q})$ be
a continuous map such that $\Omega(\tilde{\pi})\circ u$ is homotopic to the
constant map into $\tilde{p}$ in $(\Omega(X,\underline{p}),\tilde{p})$, with
homotopy $H:(CS^{k-1},\ast)\rightarrow (\Omega(X,\underline{p}),\tilde{p})$.
We set $H(s)=\tilde{\gamma}_s$, $s\in CS^{k-1}$. Then by Lemma 4.1.3, we
have $V_{o^\prime(\tilde{\gamma}_s)}=V_o$ and $\delta(\tilde{\gamma}_s)$
equals the coset of $1_{G_{V_o}}$, which implies that we can choose 
$1_{G_{V_o}}$ as the representative $g$ in the statement of Lemma 4.1.3.
Hence the liftings $\ell(\tilde{\gamma}_s)$ define a continuous map from
$(CS^{k-1},\ast)$ to $(\Omega(Y,\underline{q}),\tilde{q})$, which provides a 
homotopy between $u$ and the constant map into $\tilde{q}$. This proves the
injectivity of homomorphism $\pi_{\#}:
\pi_k(Y,\underline{q})\rightarrow \pi_k(X,\underline{p})$ for all $k\geq 1$.

On the other hand, for any continuous map $u:(S^{k-1},\ast)\rightarrow
(\Omega(X,\underline{p}),\tilde{p})$ with $k\geq 2$, we set 
$u(s)=\tilde{\gamma}_s$, $s\in S^{k-1}$, and conclude that the component 
$V_{o^\prime(\tilde{\gamma}_s)}$ equals $V_o$ and $\delta(\tilde{\gamma}_s)$
equals the coset of $1_{G_{V_o}}$, because $u(\ast)=\tilde{p}$ and 
$V_{o^\prime(\tilde{p})}=V_o$, and $\delta(\tilde{p})$ equals the coset of
$1_{G_{V_o}}$. Hence we obtain a continuous map $\ell(u):(S^{k-1},\ast)
\rightarrow (\Omega(Y,\underline{q}),\tilde{q})$ defined by the liftings
$\ell(\tilde{\gamma}_s)$, and it is easily seen that $\pi_{\#}([\ell(u)])=
[u]$. This proves the surjectivity of $\pi_{\#}$ for $k\geq 2$.

\hfill $\Box$

\noindent{\bf Definition 4.1.5:}
{\em An orbispace $(X,\U)$ is called {\it locally strongly path-connected} if 
for any basic open set $U\in\U$, the space $\widehat{U}$ in its 
$G$-structure $(\widehat{U},G_{U},\pi_{U})$ is locally path-connected.
}

\hfill $\Box$

\noindent{\bf Theorem 4.1.6:}
{\em Let $\tilde{\pi}:(Y,\underline{q})\rightarrow (X,\underline{p})$ be 
a based orbispace covering morphism. Given any based morphism 
$\tilde{\phi}:(Z,\underline{z})\rightarrow (X,\underline{p})$ from a 
connected, locally strongly path-connected orbispace $Z$, there exists a 
unique based morphism $\ell(\tilde{\phi}):(Z,\underline{z})\rightarrow 
(Y,\underline{q})$ such that 
$\tilde{\pi}\circ \ell(\tilde{\phi})=\tilde{\phi}$, if and only if 
$\phi_{\#}(\pi_1(Z,\underline{z}))\subset \pi_{\#}(\pi_1(Y,\underline{q}))$ in 
$\pi_1(X,\underline{p})$.
}

\vspace{2mm}

\noindent{\bf Proof:}
The ``only if'' part is trivial as usual. We shall prove the
``if'' part next.

Suppose that $\phi_{\#}(\pi_1(Z,\underline{z}))\subset 
\pi_{\#}(\pi_1(Y,\underline{q}))$ in $\pi_1(X,\underline{p})$.
First we shall introduce some notations. Let $\Pi=(\{V_\alpha\},
\{U_\alpha\},\{\pi_\alpha\},\{\rho_{\beta\alpha}\})$ be a representing system 
of $\tilde{\pi}$, which will be fixed throughout the proof. Let 
$\sigma=(\{W_i\},\{U_i\},\{\phi_i\},\{\eta_{ji}\})$ be a representing system 
of $\tilde{\phi}$, where each $U_i$ is an elementary neighborhood of $X$ with 
respect to the orbispace covering morphism $\tilde{\pi}$. We shall also fix 
our notations for the base-point structures: $\underline{q}=(q,V_o,\hat{q})$, 
$\underline{p}=(p,U_o,\hat{p})$ and $\underline{z}=(z,W_o,\hat{z})$. 

The strategy of the proof is to find a collection of automorphisms 
$\{\delta_i\}$, where each $\delta_i\in G_{U_i}$, such that the system 
$\delta(\sigma):=(\{W_i\},\{U_i\},\{\delta_i\circ\phi_i\},
\{\delta_j\circ\eta_{ji}\circ\delta_i^{-1}\})$, which is isomorphic to 
$\sigma$, can be lifted to $(Y,\underline{q})$ as a system through the system 
$\Pi$.

The automorphisms $\{\delta_i\}$ are obtained as follows. For each $i$, we 
pick a point $z_i\in W_i$ and a $\hat{z_i}\in\widehat{W_i}$ such that 
$\pi_{W_i}(\hat{z_i})=z_i$, with $z_o=z$ and $\hat{z_o}=\hat{z}$ where
$z,\hat{z}$ are given in the base-point structure 
$\underline{z}=(z,W_o,\hat{z})$. This gives rise to a collection of 
base-point structures $\underline{z_i}=(z_i,W_i,\hat{z_i})$ of $Z$, and a 
corresponding collection of base-point structures of $X$: 
$\underline{p_i}=(p_i,U_i,\hat{p_i})$ where each $p_i=\phi(z_i)$ and 
$\hat{p_i}=\phi_i(\hat{z_i})$. Now for any $W_i$, since $Z$ is connected and 
locally strongly path-connected, $Z$ is path-connected and there
exists a based 
path $\tilde{\gamma}_i\in P(Z,\underline{z},\underline{z_i})$, and therefore a 
based path $\tilde{\gamma}_i^\prime:=\tilde{\phi}\circ\tilde{\gamma}_i\in 
P(X,\underline{p},\underline{p_i})$. We apply Lemma 4.1.3 to 
$\tilde{\gamma}_i^\prime$, and we obtain a connected component 
$V_{i(\tilde{\gamma}_i^\prime)}$ of $\pi^{-1}(U_i)$ and a coset 
$\delta(\tilde{\gamma}_i^\prime)$ such that after choosing a 
representative $g_i$ of $\delta(\tilde{\gamma}_i^\prime)$, the based path 
$g_i\cdot\tilde{\gamma}_i^\prime$ in 
$P(X,\underline{p},g_i\cdot\underline{p_i})$ has a unique lifting to 
$(Y,\underline{q})$. The upshot is that the component 
$V_{i(\tilde{\gamma}_i^\prime)}$ and the coset 
$\delta(\tilde{\gamma}_i^\prime)$ are independent of the choice on 
the based path $\tilde{\gamma}_i$, because of the assumption that
$\phi_{\#}(\pi_1(Z,\underline{z}))\subset \pi_{\#}(\pi_1(Y,\underline{q}))$ 
in $\pi_1(X,\underline{p})$. Now we set $V_i=V_{i(\tilde{\gamma}_i^\prime)}$,
and choose a representative $\delta_i\in G_{U_i}$ of the coset 
$\delta(\tilde{\gamma}_i^\prime)$ for each $i$. We apply these
automorphisms $\{\delta_i\}$ to the system $\sigma$ and still write the 
resulting isomorphic system as $\sigma$ for simplicity. With this
adjustment, the system $\sigma$ has the following property: for any
based path $\tilde{\gamma}_i\in P(Z,\underline{z},\underline{z}_i)$,
the push-forward based path $\sigma\circ\tilde{\gamma}_i$ by the
system $\sigma$ can be directly lifted to a based path in
$(Y,\underline{q}, \underline{q}_i)$ where the base-point structure 
$\underline{q}_i=(q_i,V_i,\hat{q}_i)$ with
$\hat{q}_i=\pi_i^{-1}(\delta_i(\hat{p}_i))$ and $q_i=\pi_{V_i}(\hat{q}_i)$
(here $\pi_i:\widehat{V}_i\rightarrow \widehat{U}_i$ is the
homeomorphism given in the representing system $\Pi$ of the covering
morphism $\tilde{\pi}$).

We shall next show that each 
$\eta_{ji}:Tran(W_i,W_j)\rightarrow Tran(U_i,U_j)$ has its image contained in 
$\rho_{ji}(Tran(V_i,V_j))$, in particular, each 
$\eta_i:G_{W_i}\rightarrow G_{U_i}$ has its image contained in 
$\rho_i(G_{V_i})$, so that we can define a system 
$(\{W_i\},\{V_i\},\{\ell(\phi_i)\},\{\ell(\eta_{ji})\})$ by setting 
$\ell(\phi_i)=\pi_i^{-1}\circ\phi_i$ and 
$\ell(\eta_{ji})=\rho_{ji}^{-1}\circ\eta_{ji}$, which defines the desired 
based morphism 
$\ell(\tilde{\phi}): (Z,\underline{z})\rightarrow (Y,\underline{q})$. 

Given any $\xi\in Tran(W_i,W_j)$, suppose it is an isomorphism between induced 
$G$-structures of a connected component $W$ of $W_i\cap W_j$. We pick a point 
$\hat{z^\prime}$ in the induced $G$-structure of $W$ from $W_i$. Since $Z$ is 
connected and locally strongly path-connected, there is a based path in 
$P(Z,\underline{z},\underline{z_i})$ and a path in $\widehat{W_i}$ connecting 
$\hat{z_i}$ and $\hat{z^\prime}$, and a path in $\widehat{W_j}$ connecting 
$\xi(\hat{z^\prime})$ and $\hat{z_j}$. These paths altogether define a based 
path $\tilde{\gamma}\in P(Z,\underline{z},\underline{z_j})$. The push-forward 
path $\sigma\circ\tilde{\gamma}\in 
P(X,\underline{p},\delta_j\cdot\underline{p_j})$ 
has the property that it has a representing system whose last transition map 
is $\eta_{ji}(\xi)$, and the system can be lifted up directly 
to $(Y,\underline{q},\underline{q}_j)$. 
Hence $\eta_{ji}(\xi)$ lies in the image of $\rho_{ji}$. 

The uniqueness of $\ell(\tilde{\phi})$ follows from the nature of the
construction. This concludes the proof.

\hfill $\Box$

Here is a consequence of Theorem 4.1.6. 
Let $\tilde{\pi}_1:(Y_1,\underline{q_1})\rightarrow 
(X,\underline{p})$ and $\tilde{\pi}_2:(Y_2,\underline{q_2})\rightarrow 
(X,\underline{p})$ be two connected, locally strongly
path-connected orbispace covering 
of $(X,\underline{p})$. Suppose that 
$(\pi_1)_{\#}(\pi_1(Y_1,\underline{q_1}))$ is contained in 
$(\pi_2)_{\#}(\pi_1(Y_2,\underline{q_2}))$. Then the orbispace covering 
$\tilde{\pi}_1$ factors through the orbispace covering $\tilde{\pi}_2$ with 
a based morphism 
$\tilde{\pi}:=\ell(\tilde{\pi}_1):(Y_1,\underline{q_1})\rightarrow 
(Y_2,\underline{q_2})$, which is obviously also an orbispace covering. This 
justifies the following 

\vspace{1.5mm}

\noindent{\bf Definition 4.1.7:}
{\em A connected, locally strongly path-connected orbispace covering 
$\tilde{\pi}: Y\rightarrow X$ is called universal if $\pi_1(Y)$ is trivial.
}

\hfill $\Box$ 

We shall next consider the question of existence of universal covering for
a given orbispace. We first introduce the following  

\vspace{1.5mm}

\noindent{\bf Definition 4.1.8:}
{\em An orbispace $X$ is called semilocally 1-connected if for any $p\in X$, 
there is a basic open set $U$ containing $p$, such that the composition of 
homomorphisms $\pi_1(\widehat{U},\hat{p})\rightarrow \pi_1(U,\underline{p})
\rightarrow\pi_1(X,\underline{p})$ has a trivial image for any base-point
structure $\underline{p}=(p,U,\hat{p})$. 
}

\hfill $\Box$

\noindent{\bf Theorem 4.1.9:}
{\em Suppose orbispace $X$ is connected, locally strongly path-connected and 
semilocally 1-connected. Fixing a base-point structure 
$\underline{p}=(p,U_o,\hat{p})$. Then for any subgroup $H$ of
$\pi_1(X,\underline{p})$, there is an orbispace covering morphism 
$\tilde{\pi}:(Y,\underline{q})\rightarrow (X,\underline{p})$ such that 
$\pi_{\#}(\pi_1(Y,\underline{q}))=H$ for some connected orbispace $Y$.
}

\vspace{2mm}

\noindent{\bf Proof:}
We first construct the underlying topological space $Y_{top}$. Consider the 
set of all based morphisms
$\tilde{\gamma}:([0,1],0)\rightarrow (X,\underline{p})$, which is given the
compact-open topology (completely parallel to the case of based loop space).  
We introduce an equivalence relation $\sim_H$ between the based 
morphisms as follows: $\tilde{\gamma}_1\sim_H \tilde{\gamma}_2$ if 
$\gamma_1(1)=\gamma_2(1)$ and there are representatives $(\{I_{1,k}\},
\{U_{1,k}\},\{\gamma_{1,k}\},\{\xi_{1,lk}\})$ for $0\leq k\leq n_1$ and 
$(\{I_{2,k}\},\{U_{2,k}\},\{\gamma_{2,k}\},\{\xi_{2,lk}\})$ for 
$0\leq k\leq n_2$ of $\tilde{\gamma}_1$ and $\tilde{\gamma}_2$ respectively, 
such that $\tilde{\gamma}_1$ and $\nu(\tilde{\gamma}_2)$  
can be put together by adding a transition map 
$\xi\in Tran(U_{1,n_1},U_{2,n_2})$ satisfying 
$\gamma_{2,n_2}(1)=(\xi\circ\gamma_{1,n_1})(1)$ to form a representing system
of an element in $H$. It is easy to verify that $\sim_H$ is indeed an 
equivalence relation. We define $Y_{top}$ to be the quotient space of
the space of all based morphisms from $([0,1],0)$ to $(X,\underline{p})$ 
under $\sim_H$, which is obviously path-connected, hence connected. 
The space $Y_{top}$ has a natural base point $[\tilde{p}]$ -- the equivalence 
class of the constant morphism into $\hat{p}$. There is a natural surjective 
continuous map $\pi:Y_{top}\rightarrow X_{top}$ sending each equivalence 
class of based morphisms to its terminal point in $X_{top}$ (surjectivity 
relies on the fact that $X$ is path-connected). Clearly we have 
$\pi([\tilde{p}])=p$ and we set $q=[\tilde{p}]$.

Next we put an orbispace structure on $Y_{top}$. (The resulting 
orbispace will be taken for $Y$.)  For any point $y\in Y_{top}$, 
take a representing system
$\sigma=(\{I_k\},\{U_k\},\{\gamma_k\},\{\xi_{lk}\})$ of $y$, where the index
$k$ is running from $0$ to $n$. We set $U_\sigma=U_n$ and 
$\hat{\sigma}=\gamma_n(1)\in \widehat{U_\sigma}$. 
We may assume that the semilocally 1-connectedness holds for
$U_\sigma$ without loss of generality.
We define a map $\pi_\sigma:\widehat{U_\sigma}\rightarrow Y_{top}$ as 
follows. For each $z\in\widehat{U_\sigma}$, we connect $\hat{\sigma}$ 
to $z$ by a path $\gamma_{z}$ in $\widehat{U_\sigma}$ (the existence
of $\gamma_z$ is ensured by the assumption that $X$ is locally
strongly path-connected), and define $\pi_\sigma(z)$ to be the
equivalence class under $\sim_H$ of the based morphism 
obtained from extending the system $\sigma$ by adding the path
$\gamma_z$. The assumption that $X$ is semilocally 1-connected ensures 
that the map $\pi_\sigma$ is well-defined, and the assumption that $X$
is locally strongly path-connected implies that $\pi_\sigma$ is
continuous. We set $V_\sigma=\pi_\sigma(\widehat{U_\sigma})$ in
$Y_{top}$, which obviously satisfies
$\pi(V_\sigma)=U_\sigma$. Moreover, $V_\sigma$ is connected since
$\widehat{U_\sigma}$ is.

We will show
(1) $V_\sigma$ is an open neighborhood of $y$ in $Y_{top}$, (2) There is a 
subgroup $G_{V_\sigma}$ of $G_{U_\sigma}$ such that 
$\widehat{U_\sigma}/G_{V_\sigma}$ is homeomorphic to $V_\sigma$ under 
$\pi_\sigma$, (3) By setting $\widehat{V_\sigma}:=\widehat{U_\sigma}$, 
$\pi_{V_\sigma}:=\pi_\sigma$ and define $(\widehat{V_\sigma}, G_{V_\sigma},
\pi_{V_\sigma})$ to be the $G$-structure of $V_\sigma\subset Y_{top}$, we 
actually obtain an orbispace structure on $Y_{top}$, (4) The identification 
$\widehat{V_\sigma}:=\widehat{U_\sigma}$ and inclusion 
$G_{V_\sigma}\subset G_{U_\sigma}$ can be fitted together to 
define an orbispace covering morphism $\tilde{\pi}:Y\rightarrow X$, 
and (5) We have $\pi_{\#}(\pi_1(Y,\underline{q}))=H$ for the base-point 
structure $\underline{q}=(q,V_o,\hat{q})$ where $\hat{q}:=\hat{p}$ in 
$\widehat{V_o}:=\widehat{U_o}$. 

For (1), suppose $\tilde{\gamma}_0:([0,1],0)\rightarrow (X,\underline{p})$ 
is a based morphism such that $\tilde{\gamma}_0\sim_H \pi_\sigma(z_0)$ for 
some $z_0\in\widehat{U_\sigma}$, we need to show that for any 
based morphism $\tilde{\gamma}$ sufficiently close to
$\tilde{\gamma}_0$ in the compact-open topology, 
there is a $z\in\widehat{U_\sigma}$ such that $\pi_\sigma(z)\sim_H
\tilde{\gamma}$. We take a representative $\tau_0$ for $\tilde{\gamma}_0$. 
Without loss of generality, we can assume that $\tau_0=(\{I_i\},\{U_i\},
\{\gamma_{0,i}\},\{\xi_{0,ji}\})$, where $i$ is running from $0$ to $m$,
such that $U_m=U_\sigma$ and $\gamma_{0,m}(1)=z_0$. As in the case of based
loops, we can show that a neighborhood of $\tilde{\gamma}_0$ under the
compact-open topology can be represented by a set of systems 
$\tau=(\{I_i\},\{U_i\},\{\gamma_{i}\},\{\xi_{0,ji}\})$. We pick $t_i\in
I_i\cap I_{i+1}$ for $i=0,\cdots,m-1$. Suppose each transition map $\xi_{0,ji}$
is defined over an open subset $W_i\subset\widehat{U_i}$, which is 
path-connected. We take a path $u_i$ in $W_i$ running from $\gamma_i(t_i)$ to
$\gamma_{0,i}(t_i)$ for each $i=0,\cdots,m-1$, and take a path $u$ in 
$\widehat{U_\sigma}$ running from $\gamma_m(1)$ to $\gamma_{0,m}(1)=z_0$.
We define for each $\tau$ a new system 
$\tau^\prime=(\{I_i\},\{U_i\},\{\gamma_{i}^\prime\},\{\xi_{0,ji}\})$ where
for each $0\leq i\leq m-1$, $\gamma_i^\prime$ is obtained from precomposing 
$\gamma_i$ by $\xi_{0,i(i-1)}\circ \nu(u_{i-1})$ and post-composing 
$\gamma_i$ by $u_i$, and $\gamma_m^\prime$ is obtained from precomposing 
$\xi_{0,m(m-1)}\circ \nu(u_{m-1})$ and post-composing $u\# \nu(u)$. The system
$\tau^\prime$ is clearly homotopic to $\tau$. On the other hand, by semilocally
1-connectedness of $X$, $\tau^\prime$ is equivalent under $\sim_H$ to the 
system $(\{I_i\},\{U_i\},\{\gamma^{(i)}\},\{\xi_{0,ji}\})$ where each
$\gamma^{(i)}=\gamma_{0,i}$ for $0\leq i\leq m-1$, and $\gamma^{(m)}=
\gamma_{0,m}\# \nu(u)$. This means exactly that $\pi_\sigma(z)=[\tau]$ where
$z=\gamma_m(1)\in \widehat{U_\sigma}$. This concludes the proof of (1).

For (2), we obtain the subgroup $G_{V_\sigma}$ as follows. We denote by
$\underline{\sigma}$ the base-point structure 
$(\pi(y),U_\sigma,\hat{\sigma})$ and
consider the isomorphism $\sigma_\ast:\pi_1(X,\underline{p})\rightarrow 
\pi_1(X,\underline{\sigma})$ induced by the based path $[\sigma]\in 
P(X,\underline{p},\underline{\sigma})$ determined by the system $\sigma$.
Set $H_\sigma=\sigma_\ast(H)$. Now for each $g\in G_{U_\sigma}$, we take
a path $\gamma_g$ in $\widehat{U}_\sigma$ 
connecting $\hat{\sigma}$ and $g\cdot\hat{\sigma}$. 
Then the pair $(\gamma_g,g)$ determines an element $[g]$ in 
$\pi_1(U_\sigma,\underline{\sigma})$ whose image in 
$\pi_1(X,\underline{\sigma})$ is independent of the choice on the path
$\gamma_g$ by the semilocally 1-connectedness of $X$. We simply put 
$G_{V_\sigma}:=\{g\in G_{U_\sigma}|[g]\in H_\sigma\}$. 
The fact that $(\gamma_g,g)\# (\gamma_h,h)=(\gamma_{gh},gh)$ for some
path $\gamma_{gh}$ connecting $\hat{\sigma}$ to $gh\cdot \hat{\sigma}$
in $\widehat{U}_\sigma$ shows that $G_{V_\sigma}$ is a subgroup of
$G_{U_\sigma}$. It remains to show
that $\pi_\sigma$ induces a homeomorphism between 
$\widehat{U_\sigma}/G_{V_\sigma}$ and $V_\sigma$. First of all, $\pi_\sigma:
\widehat{U_\sigma}\rightarrow V_\sigma$ is $G_{V_\sigma}$-invariant. This can
be seen as follows. Given any $z\in\widehat{U_\sigma}$, we take a path $\gamma$
running from $\hat{\sigma}$ to $z$, and let $\sigma_z$ be the system 
obtained from
post-composing $\sigma$ with $\gamma$, which defines $\pi_\sigma(z)$. For
$g\cdot z$, we use the path $\gamma_g\# g\circ\gamma$ and obtain a system
$\sigma_{g\cdot z}$ correspondingly to define $\pi_\sigma(g\cdot z)$. Now it
is easily seen that the systems $\sigma_{g\cdot z}$ and $\nu(\sigma_z)$ can be
put together by adding the transition map $g^{-1}\in Tran(U_\sigma,U_\sigma)$
to define a system, which represents the element $\sigma_\ast^{-1}([g])$ in
$H\subset \pi_1(X,\underline{p})$. Hence $\pi_\sigma(g\cdot z)=\pi_\sigma(z)$.
Secondly, suppose $\pi_\sigma(z)=\pi_\sigma(z^\prime)$. Take paths $\gamma_z$,
$\gamma_{z^\prime}$ connecting $\hat{\sigma}$ to $z$ and $z^\prime$ 
respectively.
Let $\sigma_z$ and $\sigma_{z^\prime}$ be the resulting systems to be used
to define $\pi_\sigma(z)$ and $\pi_\sigma(z^\prime)$. Then there is a $g\in
G_{U_\sigma}$ such that by adding $g$ the systems $\sigma_z$ and 
$\nu(\sigma_{z^\prime})$ can be  put together to form a system representing an
element $h\in H$. Clearly $z^\prime=g\cdot z$. If we join $\hat{\sigma}$ and
$g\cdot\hat{\sigma}$ by the path 
$\gamma_g:=\gamma_{z^\prime}\# \nu(g\circ\gamma_z)$,
we see that the class $[(\gamma_g,g)]$ in 
$\pi_1(U_\sigma,\underline{\sigma})$ equals $\sigma_\ast(h^{-1})\in H_\sigma$.
Hence we have $g\in G_{V_\sigma}$. From here it is easy to see that 
$\pi_\sigma$ induces a homeomorphism between 
$\widehat{U_\sigma}/G_{V_\sigma}$ and $V_\sigma$.

For (3), suppose 
$(\widehat{V_{\sigma_1}},G_{V_{\sigma_1}},\pi_{V_{\sigma_1}})$, 
$(\widehat{V_{\sigma_2}},G_{V_{\sigma_2}},\pi_{V_{\sigma_2}})$ 
are two $G$-structures on $Y_{top}$ constructed from based systems 
$\sigma_1$, $\sigma_2$ respectively, such that $V_{\sigma_1}\cap V_{\sigma_2}
\neq \emptyset$. We need to define the set of transitions $Tran(V_{\sigma_1},
V_{\sigma_2})$, and verify the axioms of Definition 2.1.2 ($Y_{top}$ is
locally connected since each $V_\sigma$ is as the orbit space of locally
connected space $\widehat{U_\sigma}$ by the action of $G_{V_\sigma}$).
The assumption $V_{\sigma_1}\cap V_{\sigma_2}\neq \emptyset$ implies that
$U_{\sigma_1}\cap U_{\sigma_2}\neq \emptyset$ since $\pi(V_\sigma)=U_\sigma$.
The set of transition maps $Tran(V_{\sigma_1},V_{\sigma_2})$ will be defined
as a subset of $Tran(U_{\sigma_1},U_{\sigma_2})$ as follows. Given $\xi\in
Tran(U_{\sigma_1},U_{\sigma_2})$, which is defined over the subset $W_1$ of
$\widehat{U_{\sigma_1}}$, if there is $z\in W_1$, with 
$\xi(z)\in W_2:=\xi(W_1)\subset
\widehat{U_{\sigma_2}}$ such that $\pi_{\sigma_1}(z)=\pi_{\sigma_2}(\xi(z))$,
then we put $\xi$ in $Tran(V_{\sigma_1},V_{\sigma_2})$. It is easy to see
that $\pi_{\sigma_1}(z)=\pi_{\sigma_2}(\xi(z))$ for some $z\in W_1$ implies
that $\pi_{\sigma_1}=\pi_{\sigma_2}\circ\xi$ on $W_1$, and 
$\pi_{\sigma_1}(W_1)$ is a connected component of 
$V_{\sigma_1}\cap V_{\sigma_2}$. Let $\Gamma_1=\{g\in G_{U_{\sigma_1}}|g\cdot
W_1=W_1\}$ and $\Gamma_2=\{g\in G_{U_{\sigma_2}}|g\cdot W_2=W_2\}$. It remains
to show that $\xi:\Gamma_1\cap G_{V_{\sigma_1}}\rightarrow 
\Gamma_2\cap G_{V_{\sigma_2}}$ is an isomorphism, so that $\xi$ is indeed an
isomorphism between the $G$-structure of $\pi_{\sigma_1}(W_1)$ induced
from $(\widehat{V_{\sigma_1}},G_{V_{\sigma_1}},\pi_{V_{\sigma_1}})$
and the $G$-structure of $\pi_{\sigma_2}(W_2)$ induced from
$(\widehat{V_{\sigma_2}},G_{V_{\sigma_2}},\pi_{V_{\sigma_2}})$.
This can be seen as follows. We pick a $z_1\in W_1$ and let $z_2=
\xi(z_1)\in W_2$. We take a path $\gamma_i$ connecting $\hat{\sigma_i}$ to
$z_i$, and let $\tau_i$ be the system $\sigma_i\#\gamma_i$ for $i=1,2$.
Let $h\in H$ be the element obtained from composing $\tau_1$, 
$\xi$ and $\nu(\tau_2)$, and $h_{\tau_1}=(\tau_1)_\ast(h)$ in 
$\pi_1(X,\underline{\tau_1})$. Denote by $\xi_\ast$ the isomorphism between 
$\pi_1(X,\underline{\tau_1})$ and $\pi_1(X,\underline{\tau_2})$ induced by 
$\xi$. Then it can be easily verified that for any 
$a\in\pi_1(X,\underline{p})$, we have 
$(\tau_2)_\ast(a)=\xi_\ast(h_{\tau_1}^{-1}\cdot (\tau_1)_\ast(a)\cdot 
h_{\tau_1})$. The required isomorphism $\xi_\ast:\Gamma_1\cap 
G_{V_{\sigma_1}}\rightarrow \Gamma_2\cap G_{V_{\sigma_2}}$ follows trivially 
from the above equation. The axioms in Definition 2.1.2 can be easily checked
for $Tran(V_{\sigma_1},V_{\sigma_2})$ because they are satisfied by 
$Tran(U_{\sigma_1},U_{\sigma_2})$, and $Tran(V_{\sigma_1},V_{\sigma_2})$ is
a subset of $Tran(U_{\sigma_1},U_{\sigma_2})$. This concludes the proof of (3).

For (4), there is indeed a morphism $\tilde{\pi}:Y\rightarrow X$ given by
the system $(\{V_\sigma\},\{U_\sigma\},\{\pi^\sigma\},\{\rho_{\sigma\tau}\})$,
where each $\pi^\sigma:\widehat{V_\sigma}\rightarrow \widehat{U_\sigma}$ is
the identity map, and each $\rho_{\sigma\tau}:Tran(V_\tau,V_\sigma)\rightarrow
Tran(U_\tau,U_\sigma)$ is the natural inclusion. The conditions $b),c)$ and 
$d)$ in Definition 4.1.1 are clearly satisfied by $\tilde{\pi}$, by the nature
of construction. As for $a)$, we need to show that for any $U\in\{U_\sigma\}$,
a connected component of $\pi^{-1}(U)$ is of form $V_\sigma$. Let $V$ be a
connected component of $\pi^{-1}(U)$. Then for any $y\in V$, since $\pi(y)\in
U$, there is a system $\sigma_y$ representing $y$, with a canonically 
constructed neighborhood $V_{\sigma_y}$ where 
$\widehat{V_{\sigma_y}}=\widehat{U}$. Since $V_{\sigma_y}$ is connected and 
$\pi(V_{\sigma_y})=U$, we have $V_{\sigma_y}\subset V$. Hence 
$V=\cup_{y\in V}V_{\sigma_y}$. On the other hand, it is easy to see that if
$V_\sigma$ and $V_\tau$ have non-empty intersection and 
$\pi(V_\sigma)=\pi(V_\tau)$, then $V_\sigma=V_\tau$. Hence $V=V_{\sigma_y}$
for any $y\in V$, and $a)$ is verified.

For (5), it suffices to show that for any based morphism 
$\tilde{\gamma}:([0,1],0)
\rightarrow (X,\underline{p})$, if we denote its equivalence class under
$\sim_H$ by $[\tilde{\gamma}]_H$, which is a point in $Y_{top}$, we have
$[\tilde{\gamma}]_H=\ell(\gamma)(1)$ where $\ell(\gamma)$ stands for the
induced map $[0,1]\rightarrow Y_{top}$ of the unique lifting 
$\ell(\tilde{\gamma})$ of $\tilde{\gamma}$ to $(Y,\underline{q})$, which was
shown to exist in Lemma 4.1.3. From the construction of $Y$, we see that if
$[\tilde{\gamma}]_H=\ell(\gamma)(1)$ is true for some $\tilde{\gamma}$ with
$[\tilde{\gamma}]_H\in V_\sigma$ for some $V_\sigma$, then 
$[\tilde{\gamma}]_H=\ell(\gamma)(1)$ is true for all $\tilde{\gamma}$ such
that $[\tilde{\gamma}]_H\in V_\sigma$. With this understood, we now appeal to
the fact that $Y_{top}$ is connected and $Y_{top}=\cup V_\sigma$, and the
claim, hence (5) follows.

\hfill $\Box$

\noindent{\bf Definition 4.1.10:}
{\em Given an orbispace covering morphism $\tilde{\pi}:Y\rightarrow X$, a 
deck transformation of $\tilde{\pi}$ is an invertible morphism
$\tilde{\phi}:Y\rightarrow Y$ which satisfies  
$\tilde{\pi}\circ\tilde{\phi}=\tilde{\pi}$. 
The set of all deck transformations of $\tilde{\pi}$ forms a group acting on
$Y$, which is called the {\it group of deck transformations} of the orbispace
covering $\tilde{\pi}$, and is denoted by $Deck(\tilde{\pi})$.
}

\hfill $\Box$

Let $\tilde{\pi}:Y\rightarrow X$ be an orbispace covering, 
$\Pi=(\{V_\alpha\},\{U_\alpha\},\{\pi_\alpha\},\{\rho_{\beta\alpha}\})$,
$\alpha\in\Lambda$, be a fixed representing system of $\tilde{\pi}$.
For each $\alpha$, set $K_{U_\alpha}=\{g\in G_{U_\alpha}|g\cdot x=x,
\forall x\in\widehat{U_\alpha}\}$.
We consider the set 
$$
K_\Pi=\{\langle g \rangle|\langle g
\rangle=\{g_\alpha|\alpha\in\Lambda\}, g_\alpha\in K_{U_\alpha}\} 
\leqno (4.1.2)
$$
where each $\langle g \rangle=\{g_\alpha\}$ satisfies the 
following compatibility conditions 
$$
g_\beta\circ\rho_{\beta\alpha}(\xi)=\rho_{\beta\alpha}(\xi)\circ g_\alpha,
\hspace{2mm} \forall \xi\in Tran(V_\alpha,V_\beta). \leqno (4.1.3)
$$
These compatibility conditions are quite strong in nature. In
particular, the following relations are contained in $(4.1.3)$:
$$
g_\alpha\rho_\alpha(h)=\rho_\alpha(h)g_\alpha, \forall h\in
G_{V_\alpha}, \hspace{3mm} \varphi_{\beta\alpha}(g_\alpha)=g_\beta
\leqno (4.1.4)
$$
where $\varphi_{\beta\alpha}$ is the group isomorphism in the
transition map $\rho_{\beta\alpha}$. 
The set $K_\Pi$ becomes a group under the multiplication 
$\langle g\rangle\cdot\langle h\rangle:=\langle gh\rangle$, where 
$\langle gh\rangle:=\{g_\alpha h_\alpha|\alpha\in\Lambda\}$ if 
$\langle g\rangle=\{g_\alpha|\alpha\in\Lambda\}$ and 
$\langle h\rangle=\{h_\alpha|\alpha\in\Lambda\}$.
Let $C_\Pi$ be the subgroup of $K_\Pi$ which consists of elements 
$\langle g\rangle=\{g_\alpha|\alpha\in\Lambda\}$ where each 
$g_\alpha\in \rho_\alpha(G_{V_\alpha})$. Then $C_\Pi$ is clearly contained
in the center of $K_\Pi$ by $(4.1.4)$, hence a normal subgroup of $K_\Pi$.
Let $\underline{q}=(q,V_o,\hat{q})$ and $\underline{p}=(p,U_o,\hat{p})$ be any 
base-point structures of $Y$ and $X$ respectively, with respect to
which $\Pi$ is a system representing 
$\tilde{\pi}:(Y,\underline{q})\rightarrow (X,\underline{p})$. Let
$\tilde{p}$ and $\tilde{q}$ be the constant map from $[0,1]$ onto
$\hat{p}$ and $\hat{q}$ respectively. Each element $g\in K_{U_o}$ or 
$h\in K_{V_o}$ defines an element $[g]\in
\pi_1(X,\underline{p})$ or $[h]\in\pi_1(Y,\underline{q})$ by setting 
$[g]=[(\tilde{p},g)]$ and $[h]=[(\tilde{q},h)]$. Clearly $\pi_\#([h])=
[\rho_o(h)]$. There is a homomorphism 
$$
\Phi_{p,q}:K_\Pi/C_\Pi\rightarrow N(\pi_{\#}(\pi_1(Y,\underline{q})))/
\pi_{\#}(\pi_1(Y,\underline{q})) \leqno (4.1.5)
$$
defined by $\langle g\rangle\mapsto [g_o]$, where $g_o$ is the component of
$\langle g\rangle$ corresponding to index $o$. Here
$N(\pi_{\#}(\pi_1(Y,\underline{q})))$ denotes the normalizer of 
$\pi_{\#}(\pi_1(Y,\underline{q}))$ in $\pi_1(X,\underline{p})$.

\vspace{2mm}

\noindent{\bf Theorem 4.1.11:}
{\em Let $\tilde{\pi}:Y\rightarrow X$ be a connected, locally strongly 
path-connected orbispace covering. For any base-point structures 
$\underline{q}=(q,V_o,\hat{q})$ and $\underline{p}=(p,U_o,\hat{p})$ of $Y$ 
and $X$ respectively, we denote $H=\pi_{\#}(\pi_1(Y,\underline{q}))$ and 
$N(H)$ the normalizer of $H$ in $\pi_1(X,\underline{p})$. Then there is a 
homomorphism $\Theta_{p,q}:N(H)/H\rightarrow Deck(\tilde{\pi})$, such that
$\Phi_{p,q}$ and $\Theta_{p,q}$ fit into a short exact sequence
$$
1\rightarrow K_{\tilde{\Pi}}/C_{\tilde{\Pi}}\stackrel{\Phi_{p,q}}{\rightarrow}
N(H)/H \stackrel{\Theta_{p,q}}{\rightarrow} Deck(\tilde{\pi})\rightarrow 1.
\leqno (4.1.6)
$$
}

\vspace{2mm}

\noindent{\bf Proof:}
Let $\tilde{\gamma}$ be a based loop in $(X,\underline{p})$ such that its
class $[\tilde{\gamma}]\in\pi_1(X,\underline{p})$ lies in the normalizer
$N(H)$ of $H$. By Lemma 4.1.3, $\tilde{\gamma}$ is associated with a pair 
$(V_{\tilde{\gamma}},\delta_{\tilde{\gamma}})$ where $V_{\tilde{\gamma}}$ is
a connected component of $\pi^{-1}(U_o)$, $\delta_{\tilde{\gamma}}$ is a coset
in $G_{U_o}/\rho_o(G_{V_o})$, such that for any chosen representative 
$g_{\tilde{\gamma}}$ of $\delta_{\tilde{\gamma}}$, there is a unique lifting
$\ell(\tilde{\gamma})\in P(Y,\underline{q},\underline{q_g})$ satisfying
$\tilde{\pi}\circ\ell(\tilde{\gamma})=g_{\tilde{\gamma}}\circ\tilde{\gamma}$,
where $\underline{q_g}$ is the base-point structure canonically determined 
from $g_{\tilde{\gamma}}\cdot\underline{p}$, and 
$g_{\tilde{\gamma}}\circ\tilde{\gamma}$ stands for the image of 
$\tilde{\gamma}$ in $P(X,\underline{p},g_{\tilde{\gamma}}\cdot\underline{p})$
induced by $g_{\tilde{\gamma}}$. We consider two based versions of 
$\tilde{\pi}$, $\tilde{\pi}^{(1)}:(Y,\underline{q})\rightarrow 
(X,g_{\tilde{\gamma}}\cdot\underline{p})$ and $\tilde{\pi}^{(2)}:
(Y,\underline{q_g})\rightarrow (X,g_{\tilde{\gamma}}\cdot\underline{p})$.
Then the fact that $[\tilde{\gamma}]\in N(H)$ implies that 
$\pi_{\#}^{(1)}(\pi_1(Y,\underline{q}))=\pi_{\#}^{(2)}
(\pi_1(Y,\underline{q_g}))$ in $\pi_1(X,g_{\tilde{\gamma}}\cdot\underline{p})$.
Hence by Theorem 4.1.6, there is a unique based morphism $\tilde{\phi}:
(Y,\underline{q})\rightarrow (Y,\underline{q_g})$ such that 
$\tilde{\pi}^{(2)}\circ\tilde{\phi}=\tilde{\pi}^{(1)}$. It is easy to see
that the corresponding morphism $\tilde{\phi}:Y\rightarrow Y$ is
independent of the choice of the representative $g_{\tilde{\gamma}}$ and 
depends only on the class of $\tilde{\gamma}$ in $N(H)/H$, and is
obviously a deck transformation of $\tilde{\pi}$. The map 
$\Theta_{p,q}:N(H)/H\rightarrow Deck(\tilde{\pi})$
is defined by $[\tilde{\gamma}]\mapsto \tilde{\phi}$.

We first verify that $\Theta_{p,q}$ is a surjective homomorphism. Let
$\tilde{\gamma}_1$, $\tilde{\gamma}_2$ be two based loops in 
$(X,\underline{p})$ such that both $[\tilde{\gamma}_1]$ and 
$[\tilde{\gamma}_2]$ are in $N(H)$. By the uniqueness of path lifting 
in Lemma 4.1.3, and the construction of $\Theta_{p,q}$, we have 
$$
\Theta_{p,q}([\tilde{\gamma}_1])\circ\ell(\tilde{\gamma}_2)=
\ell(\tilde{\gamma}_1\#\tilde{\gamma}_2). \leqno (4.1.7)
$$ 
Then the uniqueness of the factorization in Theorem 4.1.6 and the identity
$(4.1.7)$ imply that 
$$
\Theta_{p,q}([\tilde{\gamma}_1\#\tilde{\gamma}_2])=
\Theta_{p,q}([\tilde{\gamma}_1])\circ\Theta_{p,q}([\tilde{\gamma}_2]),
\leqno (4.1.8)
$$
which shows that $\Theta_{p,q}$ is a homomorphism. As for the surjectivity,
given any $\tilde{\phi}\in Deck(\tilde{\pi})$, we take a based path 
$\tilde{u}\in P(Y,\underline{q},\underline{q^\prime})$ where 
$\underline{q^\prime}$ is the base-point structure of $Y$ as the image of 
$\underline{q}$ for a choice of based versions of $\tilde{\phi}$. Then 
$\tilde{\pi}\circ\tilde{u}$ determines a based loop $\tilde{\gamma}$
in $(X,\underline{p})$. Now we observe that 
$$
\pi_{\#}(\pi_1(Y,\underline{q}))=\pi_{\#}(\pi_1(Y,\underline{q^\prime}))=
\pi_{\#}(\nu(\tilde{u})_\ast(\pi_1(Y,\underline{q})))=[\tilde{\gamma}]\cdot
\pi_{\#}(\pi_1(Y,\underline{q}))\cdot [\tilde{\gamma}]^{-1}, \leqno (4.1.9)
$$
which implies that $[\tilde{\gamma}]$ is in $N(H)$. The uniqueness 
property in Theorem 4.1.6 then asserts that 
$\Theta_{p,q}([\tilde{\gamma}])=\tilde{\phi}$. Hence $\Theta_{p,q}$ is
surjective.

Finally, we determine the kernel of $\Theta_{p,q}$. Suppose 
$\Theta_{p,q}([\tilde{\gamma}])=1$ for some $\tilde{\gamma}\in 
\Omega(X,\underline{p})$. Then from the construction of $\Theta_{p,q}$,
we see that the pair $(V_{\tilde{\gamma}},\delta_{\tilde{\gamma}})$ associated
to $\tilde{\gamma}$ in Lemma 4.1.3 must satisfy the following condition:
$V_{\tilde{\gamma}}=V_o$ and we can find a representative $g_{\tilde{\gamma}}$
of $\delta_{\tilde{\gamma}}$ such that $g_{\tilde{\gamma}}\in K_{U_o}$ and
$g_{\tilde{\gamma}}\cdot h\cdot g_{\tilde{\gamma}}^{-1}=h$ for any $h\in
\rho_o(G_{V_o})$. Moreover, $[\tilde{\gamma}]=[g_{\tilde{\gamma}}]$ in
$\pi_1(X,\underline{p})$. We need to show that there is a $\langle g\rangle\in
K_{\Pi}$ such that the component of $\langle g\rangle$ corresponding to $o$
is $g_{\tilde{\gamma}}$. This goes as follows. Recall the construction in
Theorem 4.1.6 for the morphism $\Theta_{p,q}([\tilde{\gamma}])$. For each
$V_\alpha$, we take a base-point structure $\underline{q_\alpha}=(q_\alpha,
V_\alpha,\hat{q_\alpha})$ and take a based path $\tilde{u}_\alpha\in
P(Y,\underline{q},\underline{q_\alpha})$. Then we apply Lemma 4.1.3 to
$\tilde{\gamma}_\alpha:=\tilde{\pi}^{(1)}\circ\tilde{u}_\alpha$ 
with respect to $\tilde{\pi}^{(2)}$,
where $\tilde{\pi}^{(1)}:(Y,\underline{q})\rightarrow 
(X,g_{\tilde{\gamma}}\cdot\underline{p})$ and $\tilde{\pi}^{(2)}:
(Y,\underline{q_g})\rightarrow (X,g_{\tilde{\gamma}}\cdot\underline{p})$.
The condition $\Theta_{p,q}([\tilde{\gamma}])=1$ implies that
the associated $(V_{\tilde{\gamma}_\alpha},\delta_{\tilde{\gamma}_\alpha})$
satisfies the conditions that $V_{\tilde{\gamma}_\alpha}=V_\alpha$
and after choosing representatives $g_{\alpha}$ for each 
$\delta_{\tilde{\gamma}_\alpha}$, the constructed system representing 
$\Theta_{p,q}([\tilde{\gamma}])$ must be isomorphic to the identity system.
This means that after modifications on $\{g_\alpha\}$, we obtain an element
of $K_\Pi$, whose component corresponding to index $o$ is 
$g_{\tilde{\gamma}}$. Hence we have shown that $\ker\;\Theta_{p,q}\subset
Im\;\Phi_{p,q}$. On the other hand, it is obvious that $Im\;\Phi_{p,q}\subset
\ker\;\Theta_{p,q}$. Hence $\ker\;\Theta_{p,q}=Im\;\Phi_{p,q}$. Finally,
the homomorphism $\Phi_{p,q}$ is injective because if 
$\langle g\rangle\mapsto [g_o]\in H$, then $g_o\in Im\;\rho_o$ and $(4.1.4)$
implies $g_\alpha\in Im\;\rho_\alpha$ for all other $\alpha\in\Lambda$,
and $\langle g\rangle$ lies in $C_\Pi$. This concludes the proof of $(4.1.6)$,
hence the theorem.

\hfill $\Box$

As an application, we compare our construction in the case of
orbifolds with Thurston's orbifold covering theory (cf. \cite{Th}).  
Recall that according to Thurston's
definition, a {\it covering orbifold} of an orbifold $X$ is an orbifold $Y$
together with a projection $\pi:Y\rightarrow X$ satisfying the following
condition: for each $p\in X$ there is a neighborhood $U$ uniformized by
$(V,G)$ such that for each connected component $U_i$ of $\pi^{-1}(U)$ in
$Y$, the uniformizing system of $U_i$ is $(V,G_i)$ for some subgroup $G_i$
of $G$. The {\it universal covering} of a connected orbifold $X$ is a
connected covering orbifold $\pi_0:X_0\rightarrow X$ such that for any 
connected covering orbifold $Y$ of $X$ with $\pi:Y\rightarrow X$, $X_0$ is
a covering orbifold of $Y$, with $pr: X_0\rightarrow Y$ factoring $\pi_0:
X_0\rightarrow X$ through $\pi:Y\rightarrow X$. Thurston has shown that 
universal covering always exists, and defined the fundamental group of a
connected orbifold to be the group of deck transformations of its universal
covering.

\vspace{2mm}

\noindent{\bf Theorem 4.1.12:}
{\em For a connected orbifold $X$, $\pi_1(X)$ is isomorphic to Thurston's
orbifold fundamental group of $X$.
}

\vspace{2mm}

\noindent{\bf Proof:}
An orbispace covering in the sense of Definition 4.1.1 is a Thurston's 
orbifold covering in the case of orbifolds. On the other hand, an orbifold
is locally strongly path-connected and semilocally 1-connected, hence by
Theorem 4.1.9, a connected orbifold has a universal covering in the sense
of Definition 4.1.7, which is a Thurston's universal covering by Theorem
4.1.6. Now the short exact sequence $(4.1.6)$ implies that $\pi_1$ of the
orbifold is isomorphic to the group of deck transformations since orbifolds
are normal and reduced. This concludes the proof.

\hfill $\Box$

We give an example to illustrate our orbispace covering theory.

\vspace{1.5mm}

\noindent{\bf Example 4.1.13:}\hspace{2mm}
Let $X$ be the orbispace $(S^1,\U_\tau)$ in Remark 2.1.4 e, where for 
simplicity we assume that the group $G$ has a discrete topology. We will
show in the next section that $\pi_1(X)$ is the semi-direct product of
$G$ by $\Z$ with respect to the homomorphism $\Z\rightarrow Aut(G)$ given by
$1\mapsto \tau$ (cf. Example 4.2.9). 
Let $H$ be a subgroup of $G$, which naturally becomes a
subgroup of $\pi_1(X)$. The orbispace $X$ is obviously connected, locally
strongly path-connected and semilocally 1-connected. Hence by Theorem 4.1.9,
there is an orbispace covering morphism $\tilde{\pi}:Y\rightarrow X$ such that
$\pi_{\#}(\pi_1(Y))=H$ in $\pi_1(X)$. It is easy to see that $Y$ is the
global orbispace defined by $(\R,H)$ where $H$ acts on $\R$ trivially. If
we use the notations introduced in Remark 2.1.4 e, then $X$ is covered by
$U_1$ and $U_2$ which are elementary neighborhoods with respect to 
$\tilde{\pi}$, and $\pi^{-1}(U_1)=\{V_{n,1}\}$ where $V_{n,1}=(n-\frac{1}{8},
n+\frac{5}{8})$, $n\in\Z$, and $\pi^{-1}(U_2)=\{V_{n,2}\}$, where
$V_{n,2}=(n+\frac{3}{8},n+\frac{9}{8})$, $n\in\Z$. A representing system
of $\tilde{\pi}$ defined over these open sets is given by $\Pi=(\{\pi_{n,1},
\pi_{n,2}\},\{\rho_{n,1},\rho_{n,2}\})$, where 
$\pi_{n,j}:V_{n,j}\rightarrow U_j$, $\rho_{n,j}:H\rightarrow G$, $j=1,2$, are
defined by $\pi_{n,j}(t)=t-n$, and $\rho_{n,j}(h)=\tau^n(h)$. Let $N(H)$ be
the normalizer of $H$ in $G$. Then the normalizer of $H$ in $\pi_1(X)$ is
either the semi-direct product of $N(H)$ by $\Z$ with respect to
the homomorphism $1\mapsto\tau$, or $N(H)$, depending on whether $H$ is 
invariant under $\tau$ or not. On the other hand, $K_\Pi/C_\Pi$ is the 
subgroup of $N(H)/H$ consisting of classes $[g]$ where $\tau^n(g)$
lies in the center
of $H$ in $G$ for all $n\in \Z$. The group of deck transformations 
$Deck(\tilde{\pi})$ can be determined using the short exact sequence
$(4.1.6)$ accordingly. In the case when $\tilde{\pi}:Y\rightarrow X$ is the
universal covering, the exact sequence $(4.1.6)$ reduces to the canonical
short exact sequence associated to the semi-direct product structure of
$\pi_1(X)$, with $Deck(\tilde{\pi})=\Z$ and $K_\Pi=G$, $C_\Pi=\{1\}$.

\hfill $\Box$

We end this section with a criterion as to when a universal covering of an
orbispace has a trivial orbispace structure.

\vspace{1.5mm}

\noindent{\bf Theorem 4.1.14:}
{\em Let $\tilde{\pi}:Y\rightarrow X$ be a universal covering of a connected,
locally strongly path-connected orbispace $X$ with
a representing system $\Pi=(\{V_\alpha\},\{U_\alpha\},\{\pi_\alpha\},
\{\rho_{\beta\alpha}\})$. For each elementary neighborhood $U\in\{U_\alpha\}$
of $X$, pick a base-point structure $\underline{o}=(o,U,\hat{o})$, then each
pair $(\gamma,g)$ determines a class $[(\gamma,g)]$ in 
$\pi_1(X,\underline{o})$, where $g\in G_U$ and 
$\gamma:[0,1]\rightarrow\widehat{U}$ is a path such that $\gamma(0)=\hat{o}$,
$\gamma(1)=g\cdot\hat{o}$. The criterion for $Y$ to have a trivial orbispace
structure is that $[(\gamma,g)]\neq 1$ if $g\neq 1$ for any
$(U,\underline{o})$. 
}

\vspace{2mm}

\noindent{\bf Proof:}
Suppose there is a $(\gamma,g)$ with $g\neq 1$ and $[(\gamma,g)]=1$ for
some elementary neighborhood $U$ and base-point structure $\underline{o}$.
We pick a connected component $V$ of $\pi^{-1}(U)$ and fix a 
base-point structure there. Then $[(\gamma,g)]=1$ implies that the
lifting of $(\gamma,g)$ is of form $(\gamma^\prime,g^\prime)$ with
$g^\prime\in G_V$ and $\rho(g^\prime)=g$, where $\rho:G_V\rightarrow
G_U$ is the injective homomorphism given in $\Pi$. Now we see that 
$G_V$ is not trivial since $g^\prime=\rho^{-1}(g)$ and $g\neq 1$, 
which implies that $Y$ has a non-trivial orbispace structure.

On the other hand, if $Y$ has a non-trivial orbispace structure, there is
a $V\in\{V_\alpha\}$ such that $G_V$ is non-trivial. We choose a 
base-point structure $\underline{o^\prime}=(o^\prime,V,\hat{o^\prime})$,
and a pair $(\gamma^\prime,g^\prime)$ where $g^\prime\in G_V$, 
$g^\prime\neq 1$ and $\gamma^\prime:[0,1]\rightarrow\widehat{V}$ is a path 
such that $\gamma^\prime(0)=\hat{o^\prime}$, 
$\gamma^\prime(1)=g^\prime\cdot\hat{o^\prime}$. The class 
$[(\gamma^\prime,g^\prime)]=1$ in $\pi_1(Y,\underline{o^\prime})$ since 
$Y$ is universal. Now look at the image $(\gamma,g)$ of 
$(\gamma^\prime,g^\prime)$ under $\tilde{\pi}$. We have $[(\gamma,g)]
=\pi_{\#}([(\gamma^\prime,g^\prime)])=1$ and $g\neq 1$ (since
$g^\prime\neq 1$) 
for the elementary neighborhood $U:=\pi(V)$ and base-point structure 
$\underline{o}:=\tilde{\pi}(\underline{o^\prime})$. This concludes the proof.

\hfill $\Box$

As for an example illustrating this criterion, we consider the orbifold
$X$ where the underlying topological space $X_{top}=S^2$, and the orbifold
structure has three singular points $\{z_1,z_2,z_3\}$ at which the 
uniformizing system is given by a ramified covering of degree $n_1,n_2,n_3$
respectively. By the generalized Seifert-Van Kampen theorem which is proved
in section 4.3, we have 
$$
\pi_1(X)=\{\lambda_1,\lambda_2,\lambda_3|\lambda_i^{n_i}=1,\lambda_1\lambda_2
\lambda_3=1\}. \leqno (4.1.10)
$$
The classes of form $[(\gamma,g)]$ are given by $\lambda_i^{m_i}$ where
$1\leq m_i\leq n_i-1$, which are obviously non-zero in $\pi_1(X)$. Hence
by Theorem 4.1.14, the universal covering of $X$ has a trivial orbifold
structure (which must be either $S^2$ or $\R^2$). In Thurston's
terminology (cf. \cite{Th}), $X$ is called a good orbifold, as being
the quotient of a smooth manifold by a proper discrete action.

\vspace{2mm}

\noindent{\bf Remark 4.1.15:}\hspace{2mm} A complex of groups $G(X)$
is called developable if there is a simply connected simplicial
complex $Y$ with a simplicial group action of $G$ without
inversion, such that $G(X)$ is the associated complex of groups of
$(Y,G)$. In this case, $G=\pi_1(G(X),\ast)$. With Theorem
3.2.10 at hand, it is easy to verify that a complex of groups is
developable if and only if the associated orbihedron has a universal
covering with trivial orbispace structure. Then Theorem 4.1.14 recovers a
theorem of Haefliger on this matter (cf. Theorem 4.1 in \cite{Ha2}). 
In fact, the whole covering theory of complexes of groups developed in 
\cite{Ha2} can be recovered from the orbispace covering theory
developed in this section.

\subsection{Orbispace fibration}

\hspace{5mm}
This section deals with generalizing the notion of fibration and the
associated Serre exact sequence of homotopy groups to the orbispace category.

Recall that a continuous map $p:E\rightarrow B$ between topological spaces
is said to have the homotopy lifting property with respect to a topological
space $T$ if for every map $u:T\rightarrow E$ and homotopy $H:T\times [0,1]
\rightarrow B$ of $p\circ u$ there is a homotopy 
$\widetilde{H}:T\times [0,1]\rightarrow
E$ with $u=\widetilde{H}(\cdot,0)$ and $p\circ \widetilde{H}=H$. 
($\widetilde{H}$ is said to be a lifting of
$H$.) A map $p:E\rightarrow B$ is called a fibration (resp. weak fibration)
if it has the homotopy lifting property with respect to all spaces (resp.
all disks $D^n$, $n\geq 0$). Let $b_0\in B$ and $e_0\in E$ be the base points
satisfying $p(e_0)=b_0$, and let $F=p^{-1}(b_0)$ be the fiber over base point
$b_0$. Associated to a weak fibration is the Serre exact sequence of 
homotopy groups:
$$
\begin{array}{c}
\cdots\stackrel{p_\ast}{\rightarrow}\pi_{k+1}(B,b_0)\stackrel{\partial}{
\rightarrow}\pi_k(F,e_0)\stackrel{i_\ast}{\rightarrow}\pi_k(E,e_0)
\stackrel{p_\ast}{\rightarrow}\pi_k(B,b_0)\stackrel{\partial}{\rightarrow}
\cdots\\
\stackrel{p_\ast}{\rightarrow}\pi_1(B,b_0)\stackrel{\partial}
{\rightarrow}
\pi_0(F,e_0)\stackrel{i_\ast}{\rightarrow}\pi_0(E,e_0)\stackrel{p_\ast}
{\rightarrow}\pi_0(B,b_0).
\end{array} \leqno (4.2.1)
$$ 

In order to motivate the concept of fibration in the orbispace category, let
us examine the following example. We consider the orbispaces defined 
canonically by G-spaces $(\{pt\},G)$, where $G$ acts trivially. We denote
such an orbispace by $X_G$. A morphism between such orbispaces $\tilde{\pi}:
X_G\rightarrow X_\Gamma$ is simply a homomorphism $\rho:G\rightarrow\Gamma$.
Let $H$ be the kernel of $\rho$. Then a natural candidate of the fiber of
$\tilde{\pi}:X_G\rightarrow X_\Gamma$ would be the orbispace $X_H$. According
to Theorem 3.4.1, the homotopy groups $\pi_k(X_G,\ast)$ of such orbispaces
are naturally isomorphic to $\pi_k(BG,\ast)$ for $k\geq 1$. For $k=0$, we
have $\pi_0(X_G,\ast)=\{\ast\}$. Hence the corresponding Serre exact sequence
associated to the morphism $\tilde{\pi}$ would be
$$
\begin{array}{c}
\cdots\stackrel{\rho_\ast}{\rightarrow}\pi_{k+1}(B\Gamma,\ast)\stackrel
{\partial}{\rightarrow}\pi_k(BH,\ast)\stackrel{i_\ast}{\rightarrow}
\pi_k(BG,\ast)\stackrel{\rho_\ast}{\rightarrow}\pi_k(B\Gamma,\ast)
\stackrel{\partial}{\rightarrow}\cdots\\
\stackrel{\rho_\ast}{\rightarrow}\pi_2(B\Gamma,\ast)
\stackrel{\partial}{\rightarrow}\pi_1(BH,\ast)\stackrel{i_\ast}{\rightarrow}
\pi_1(BG,\ast)\stackrel{\rho_\ast}{\rightarrow}\pi_1(B\Gamma,\ast)
\stackrel{\partial}{\rightarrow} 1. 
\end{array} \leqno (4.2.2)
$$ 
In order to ensure such an exact sequence, we need to impose the condition
that the homomorphism $\rho:G\rightarrow\Gamma$ is a {\it surjective weak
fibration}.

\vspace{1.5mm}

Before introducing the notion of orbispace fibration, we need a
digression where we generalize slightly our notion of $G$-structure of
a basic open set. Let $U$ be a basic open set of a given orbispace
$X$. A {\it generalized $G$-structure} of $U$ is a triple
$(\widehat{U},G_U,\pi_U)$ where $\widehat{U}$ is a locally connected
topological space with a continuous action of a topological group
$G_U$ such that (1) $\pi_U:\widehat{U}\rightarrow U$ induces a
homeomorphism between $\widehat{U}/G_U$ and $U$, (2) there is a
connected component $\widehat{U}^0$ of $\widehat{U}$ such that, if we
set $G_U^0:=\{g\in G_U|g\cdot \widehat{U}^0=\widehat{U}^0\}$ and
$\pi_U^0=\pi_U|_{\widehat{U}^0}$, then $(\widehat{U}^0,G_U^0,\pi_U^0)$
is the given $G$-structure of $U$. We will call
$(\widehat{U}^0,G_U^0,\pi_U^0)$ the {\it base component} of the
generalized $G$-structure $(\widehat{U},G_U,\pi_U)$. Let $U_i$,
$i=1,2$, be two basic open sets and
$(\widehat{U}_i,G_{U_i},\pi_{U_i})$, $i=1,2$, be the given generalized
$G$-structures of $U_i$. Then the set of transition maps 
between the two generalized $G$-structures is defined to be 
$$
\{g_2\circ\xi\circ g_1|g_i\in G_{U_i}, i=1,2, \xi\in Tran(U_1,U_2)\}.
\leqno (4.2.3)
$$
When the generalized $G$-structures are already specified, we will
still write $Tran(U_1,U_2)$ for the set $(4.2.3)$ for
simplicity. Morphisms between orbispaces may be represented by {\it
generalized systems} $(\{U_i\},\{U_i^\prime\},\{f_i\},\{\rho_{ji}\})$
where each map $f_i$ is defined between generalized $G$-structures and
each $\rho_{ji}$ is defined between the sets of $(4.2.3)$, satisfying
$(2.2.1a-b)$. In the based version, we require that the base-point 
structures are given by the base component of the generalized 
$G$-structures. End of digression. 

\vspace{2mm}

\noindent{\bf Definition 4.2.1:}
{\em Let $\tilde{\pi}:Y\rightarrow X$ be a morphism between orbispaces, and
$\pi:Y_{top}\rightarrow X_{top}$ be the induced map between underlying 
topological spaces. The morphism $\tilde{\pi}$ is called an orbispace 
fibration if $\tilde{\pi}$ is represented by a generalized system 
$\Pi=(\{V_\alpha\},\{U_\alpha\},\{\pi_\alpha\},\{\rho_{\beta\alpha}\})$,
$\alpha\in\Lambda$, such that
\begin{itemize}
\item [{a)}] each homomorphism $\rho_{\alpha}:G_{V_\alpha}\rightarrow 
G_{U_\alpha}$ is a surjective weak fibration, and each 
$\rho_\alpha$-equivariant map $\pi_\alpha:\widehat{V_\alpha}\rightarrow
\widehat{U_\alpha}$ is a weak fibration, and
\item [{b)}] for any $U\in\{U_\alpha\}$, $\{V_\alpha|\alpha\in\Lambda(U)\}$ 
forms a cover of $\pi^{-1}(U)$, where $\Lambda(U)$ is the subset of 
$\Lambda$ defined by $\Lambda(U):=\{\alpha\in\Lambda|U=U_\alpha\}$.
\end{itemize}
The orbispace $Y$ will be called the total space and $X$ will be
called the base space of the orbispace fibration.
}   

\hfill $\Box$

\noindent{\bf Remark 4.2.2 a:}\hspace{2mm}
The condition that each homomorphism $\rho_{\alpha}:G_{V_\alpha}\rightarrow 
G_{U_\alpha}$ is a surjective weak fibration implies that each
$\rho_{\beta\alpha}:Tran(V_\alpha,V_\beta)\rightarrow
Tran(U_\alpha,U_\beta)$ is also a surjective weak fibration. 

\hfill $\Box$

\noindent{\bf Remark 4.2.2 b:}\hspace{2mm}
We can choose base-point structures of $Y$ and $X$ respectively, such
that $\Pi$ becomes a based generalized system with respect
to the chosen base-point structures.

\hfill $\Box$

\noindent{\bf Remark 4.2.2 c:}\hspace{2mm}
Like the case of orbispace covering, here we may also add-in all
the connected open subsets of $U\in\{U_\alpha\}$ to enlarge the
generalized system $\Pi$ (cf. Remark 4.1.2 a). But it is necessary
that we allow $\Pi$ in Definition 4.2.1 to be a generalized system
rather than a system. Here is an example. Let $Y$ be the trivial
orbispace $S^1$, but equipped with a global $G$-structure $(S^1,\Z_2)$
via a degree two covering. Let $X$ be the global orbispace defined by
the $G$-space $(S^1,\Z_2)$ where $\Z_2$ acts trivially. There is a
morphism $\tilde{\pi}:Y\rightarrow X$ defined by the system
$\Pi:=(\{S^1\},\{S^1\},\{\pi\},\{\rho\})$ where $\pi:S^1\rightarrow
S^1$ is the degree two covering map and $\rho:\Z_2\rightarrow\Z_2$ is the
identity isomorphism. Clearly $\tilde{\pi}$ is an orbispace fibration in the
sense of Definition 4.2.1. But if we consider to add-in any connected
open subset of $S^1$ without allowing generalized systems, 
the condition that each $\rho_\alpha$ is surjective will fail!

Another reason for which we allow $\Pi$ to be a generalized system is
to include {\it orbispace fiber bundles} with disconnected fibers as
examples of orbispace fibrations (cf. Example 4.2.6).

\hfill $\Box$

\noindent{\bf Remark 4.2.2 d:}\hspace{2mm}
We can assume without loss of generality that in the generalized
system $\Pi$ the generalized $G$-structure of each $U_\alpha$ is
actually the $G$-structure. We will assume this modification throughout.

\hfill $\Box$

We shall next construct the fiber of an orbispace fibration. Suppose
we are given a $U\in\{U_\alpha\}$ and a point $\hat{p}\in \widehat{U}$
such that for any $\alpha\in\Lambda(U)$,
$W_\alpha:=\pi_\alpha^{-1}(\hat{p})$ is a locally connected subspace
of $\widehat{V_\alpha}$. Let $H_\alpha$ be the kernel of
$\rho_\alpha:G_{V_\alpha}\rightarrow G_U$. Then $H_\alpha$ has an
induced action on $W_\alpha$ so that we have a locally connected $G$-space
$(W_\alpha,H_\alpha)$ for each $\alpha\in\Lambda(U)$.

\vspace{1.5mm}

\noindent{\bf Lemma 4.2.3:} {\em The $G$-spaces
$\{(W_\alpha,H_\alpha)|\alpha\in\Lambda(U)\}$ can be patched together
to define an orbispace $Z$ with a canonical pseudo-embedding
$\tilde{i}:Z\rightarrow Y$. 
}

\vspace{2mm}

\noindent{\bf Proof:} We set $Z_\alpha$ for the global orbispace
associated to the $G$-space $(W_\alpha,H_\alpha)$ for each
$\alpha\in\Lambda(U)$. We shall define a collection of homeomorphisms
$\{f_{\beta\alpha}|\alpha,\beta\in\Lambda(U)\}$ where each
$f_{\beta\alpha}$ is from an open subset $Dom(f_{\beta\alpha})$ of
$(Z_\alpha)_{top}$ into $(Z_\beta)_{top}$ (here $Dom(f_{\beta\alpha})$
may be empty) satisfying the following conditions:
$$
f_{\alpha\alpha}=Id, \hspace{2mm} f_{\beta\alpha}=f^{-1}_{\alpha\beta}, 
\hspace{2mm} f_{\gamma\beta}\circ
f_{\beta\alpha}|_{Dom(f_{\gamma\alpha})}=f_{\gamma\alpha},
\hspace{2mm} \forall \alpha,\beta,\gamma\in\Lambda(U).\leqno (4.2.4)
$$
The underlying topological space of the orbispace $Z$ will be defined
as the quotient space of $\cup_{\alpha\in\Lambda(U)} (Z_\alpha)_{top}$
under the maps $\{f_{\beta\alpha}\}$.

Let $H_{\beta\alpha}=\{\xi\in
Tran(V_\alpha,V_\beta)|\rho_{\beta\alpha}(\xi)=1_{G_U}\in Tran(U,U)\}$
for any $\alpha,\beta\in\Lambda(U)$. For any $\xi\in H_{\beta\alpha}$,
we have $\pi_\beta\circ\xi=1_{G_U}\circ\pi_\alpha=\pi_\alpha$, hence
$\pi_\beta(\xi(x))=\pi_\alpha(x)=\hat{o}$ for any $x\in
Domain(\xi)\cap W_\alpha\subset\widehat{V_\alpha}$. It is easily seen
that $\xi$ induces a homeomorphism $\phi_\xi:Domain(\xi)\cap
W_\alpha\rightarrow Range(\xi)\cap W_\beta$, with inverse
$\phi_{\xi^{-1}}$ which is the homeomorphism induced by the inverse of
$\xi$. Let $H_{\alpha,\xi}$ be the subgroup of $H_\alpha$ consisting
of elements $g\in H_\alpha$ such that $g\cdot
Domain(\xi)=Domain(\xi)$, and $H_{\beta,\xi}$ be the subgroup of
$H_\beta$ consisting of elements $h\in H_\beta$ such that $h\cdot
Range(\xi)=Range(\xi)$. We claim that $\xi$ induces an isomorphism
$\lambda_\xi:H_{\alpha,\xi}\rightarrow H_{\beta,\xi}$. This could be
seen as follows. For any $u\in H_{\alpha,\xi}$, let $u^\prime$ be the
image of $u$ under $\xi$. It suffices to show that $u^\prime\in
H_{\beta,\xi}$. Observe that $\xi\circ u=u^\prime\circ\xi$ as
transition maps, hence 
$$
\rho_\beta(u^\prime)\circ\rho_{\beta\alpha}(\xi)=\rho_{\beta\alpha}(u^\prime\circ\xi)=\rho_{\beta\alpha}(\xi\circ
u)=\rho_{\beta\alpha}(\xi)\circ\rho_\alpha(u)=\rho_{\beta\alpha}(\xi),
\leqno (4.2.5)
$$ 
which implies that $\rho_\beta(u^\prime)=1_{G_U}$ and therefore
$u^\prime$ is in $H_\beta$. 
Hence $u^\prime\in H_{\beta,\xi}$ since it satisfies $u^\prime\cdot
Range(\xi)=Range(\xi)$.

Thus we have obtained for each $\xi\in H_{\beta\alpha}$ a pair
$(\phi_\xi,\lambda_\xi)$ where $\phi_\xi:Domain(\xi)\cap
W_\alpha\rightarrow Range(\xi)\cap W_\beta$ is a
$\lambda_\xi$-equivariant homeomorphism. We observe the following
facts: (1) if $\xi^\prime\in H_{\beta\alpha}$ such that
$Domain(\xi^\prime)=Domain(\xi)$ and $Range(\xi^\prime)=Range(\xi)$,
then there exists a $u\in H_{\alpha,\xi}$ such that
$\xi^\prime=\xi\circ u=\lambda_\xi(u)\circ\xi$, and we have
$\phi_{\xi^\prime}=\phi_\xi\circ u$, and
$H_{\alpha,\xi^\prime}=H_{\alpha,\xi}$,
$H_{\beta,\xi^\prime}=H_{\beta,\xi}$,
$\lambda_{\xi^\prime}(\cdot)=\lambda_\xi(u)\lambda_\xi(\cdot)
\lambda_\xi(u)^{-1}$,
(2) if for $\xi^\prime\in H_{\beta\alpha}$, there are $g\in H_\alpha$,
$h\in H_\beta$ such that $g\cdot Domain(\xi)=Domain(\xi^\prime)$,
$h\cdot Range(\xi)=Range(\xi^\prime)$, then
$h^{-1}\circ\xi^\prime\circ g\in H_{\beta\alpha}$ and
$h^{-1}\circ\xi^\prime\circ g$ has the same domain and range as $\xi$.
The facts (1) and (2) imply that the collection of equivariant maps
$\{(\phi_\xi,\lambda_\xi)|\xi\in H_{\beta\alpha}\}$ induces a
homeomorphism from an open subset of $(Z_\alpha)_{top}$ into
$(Z_\beta)_{top}$. This homeomorphism is defined to be the map
$f_{\beta\alpha}$, which clearly satisfies $(4.2.4)$ by the nature of
construction. 

Thus we have obtained the underlying topological space of the
orbispace $Z$ as the quotient space of $\cup_{\alpha\in\Lambda(U)}
(Z_\alpha)_{top}$ under the maps $\{f_{\beta\alpha}\}$. 
For each connected component $(Z_\alpha)_i$ of $(Z_\alpha)_{top}$, 
we fix a choice of a connected component $W_{\alpha,i}$ of $W_\alpha$
which lies in the base component of the generalized $G$-structure of
$V_\alpha$, and define $(W_{\alpha,i},H_{\alpha,i})$ 
to be the $G$-structure of $(Z_\alpha)_i$ where $H_{\alpha,i}$ is the
subgroup of $H_\alpha$ fixing the component $W_{\alpha,i}$. Then
the orbispace structure on each $Z_\alpha$ together with the set of 
equivariant maps $\{(\phi_\xi,\lambda_\xi)|\xi\in H_{\beta\alpha}\}$ 
define the orbispace structure of $Z$. There is a nature morphism
$\tilde{i}:Z\rightarrow Y$ induced by the inclusions
$$
(W_\alpha,H_\alpha)\hookrightarrow (\widehat{V_\alpha},G_{V_\alpha}),
\hspace{2mm} \alpha\in\Lambda(U), \leqno (4.2.6)
$$ 
which also gives rise to a generalized 
system $\sigma$ representing the morphism $\tilde{i}:Z\rightarrow Y$. 
Finally, we observe that each projection $G_{V_\alpha}\rightarrow
G_{V_\alpha}/H_{\alpha,i}$ is a weak fibration since
$G_{V_\alpha}\rightarrow G_{V_\alpha}/H_{\alpha}=G_U$ is a weak
fibration by assumption and $H_\alpha/H_{\alpha,i}$ has a discrete
topology, which is sufficient to show that the morphism $\tilde{i}$ is
a pseudo-embedding. 

\hfill $\Box$

The orbispace $Z$ together with the pseudo-embedding
$\tilde{i}:Z\rightarrow Y$ is called the {\it fiber} of the fibration
$\tilde{\pi}$ over the base-point structure
$\underline{p}=(p,U_o,\hat{p})$ where $U_o=U$. It is clear from the
construction that for any base-point structure
$\underline{q}=(q,V_o,\hat{q})$ of $Y$ such that $\tilde{\pi}$ is a
based morphism with respect to $\underline{p}$ and $\underline{q}$,
then there is a natural base-point structure
$\underline{r}=(r,W_o,\hat{r})$ of $Z$ such that the pseudo-embedding
$\tilde{i}:Z\rightarrow Y$ is a based morphism with respect to 
$\underline{r}$ and $\underline{q}$. Lemma 4.2.3 ensures the existence
of the fiber of a given orbispace fibration over a base-point structure
$\underline{p}$ as long as each $W_\alpha=\pi_\alpha^{-1}(\hat{p})$ is
locally connected. In this case we say that the fiber over
$\underline{p}$ is well-defined. 

\vspace{2mm}

\noindent{\bf Remark 4.2.4:}\hspace{2mm}
The underlying topological space of a fiber $Z$ is not necessarily a
subspace of the underlying topological space of $Y$. For example, we
consider the fibration $\tilde{\pi}:Y\rightarrow X$ in Remark 4.2.2
c. The fiber $Z$ over any base-point structure is an orbispace of two
points, with a discrete topology and a trivial orbispace structure. On the
other hand, the map $\pi:Y_{top}\rightarrow X_{top}$ induced by the
fibration $\tilde{\pi}:Y\rightarrow X$ is clearly a homeomorphism,
hence $i:Z_{top}\rightarrow Y_{top}$ is not an embedding.
 
\hfill $\Box$

\noindent{\bf Theorem 4.2.5:}
{\em Let $\tilde{\pi}:(Y,\underline{q})\rightarrow (X,\underline{p})$
be an orbispace fibration, and $\tilde{i}:(Z,\underline{r})\rightarrow
(Y,\underline{q})$ be the fiber of $\tilde{\pi}$ over the base-point
structure $\underline{p}$. Then there is a long exact sequence 
$$
\begin{array}{c}
\cdots\stackrel{\pi_{\#}}{\rightarrow}\pi_{k+1}(X,\underline{p})
\stackrel{\partial}{\rightarrow}\pi_k(Z,\underline{r})\stackrel{i_{\#}}
{\rightarrow}\pi_k(Y,\underline{q})\stackrel{\pi_{\#}}{\rightarrow}
\pi_k(X,\underline{p})\stackrel{\partial}{\rightarrow}\cdots\\
\stackrel{\pi_{\#}}{\rightarrow}\pi_1(X,\underline{p})\stackrel{\partial}
{\rightarrow}\pi_0(Z,r)\stackrel{i_{\#}}{\rightarrow}\pi_0(Y,q)
\stackrel{\pi_{\#}}{\rightarrow}\pi_0(X,p).
\end{array}\leqno (4.2.7)
$$   
}

\vspace{2mm}

\noindent{\bf Proof:}
The exact sequence $(4.2.7)$ will be derived from the long exact
sequence $(3.3.4)$ associated to the pseudo-embedding
$\tilde{i}:(Z,\underline{r})\rightarrow (Y,\underline{q})$. 

We first show that the sequence
$\pi_0(Z,r)\stackrel{i_{\#}}{\rightarrow}\pi_0(Y,q)
\stackrel{\pi_{\#}}{\rightarrow}\pi_0(X,p)$ is exact at
$\pi_0(Y,q)$. It suffices to show that for any $y\in Y$, if $\pi(y)\in
X$ is path-connected to the base point $p\in X$ in $X$, then there is
a $z\in Z$ in the fiber such that $y\in Y$ is path-connected to
$i(z)\in Y$ in $Y$. Let
$\tau=(\{I_i\},\{U_i\},\{\gamma_i\},\{\xi_{ji}\})$, $i=0,1,\cdots,n$,
be a representing system of a path connecting $\pi(y)$ to $p$ in
$X$. By $b)$ of Definition 4.2.1, there is a $\alpha_0\in
\Lambda(U_0)$ such that $y\in V_{\alpha_0}$. Let $\hat{y}\in
\widehat{V_{\alpha_0}}$ such that
$\pi_{V_{\alpha_0}}(\hat{y})=y$. Since
$\pi_{\alpha_0}:\widehat{V_{\alpha_0}}\rightarrow \widehat{U_0}$ is a
weak fibration, there is a continuous map
$\gamma_0^\prime:I_0\rightarrow \widehat{V_{\alpha_0}}$ satisfying
$\gamma_0^\prime(0)=\hat{y}$ and $\pi_{\alpha_0}\circ
\gamma_0^\prime=\gamma_0$. Pick a point $t_1\in I_0$ such that
$\pi_{U_0}(\gamma_0(t_1))$ lies in $U_0\cap U_1$. Then by $b)$ of
Definition 4.2.1 again, there is a $\alpha_1\in\Lambda(U_1)$ such that
$\pi_{V_{\alpha_0}}(\gamma_0^\prime(t_1))$ lies in $V_{\alpha_0}\cap
V_{\alpha_1}$. By the assumption that
$\rho_{\alpha_1\alpha_0}:Tran(V_{\alpha_0},V_{\alpha_1})\rightarrow
Tran(U_0,U_1)$ is surjective, there is a transition map
$\xi^\prime_{10}$ such that
$\rho_{\alpha_1\alpha_0}(\xi^\prime_{10})=\xi_{10}$. Now using the
assumption that $\pi_{\alpha_1}:\widehat{V_{\alpha_1}}\rightarrow
\widehat{U_1}$ is a weak fibration, we find a continuous map
$\gamma_1^\prime: I_1\rightarrow \widehat{V_{\alpha_1}}$ such that
$\gamma^\prime_1=\xi_{10}^\prime\circ \gamma_0^\prime$ on $I_0\cap
I_1$ and $\pi_{\alpha_1}\circ\gamma_1^\prime=\gamma_1$. We continue
with this process and construct a generalized system
$\tau^\prime=(\{I_i\},\{V_{\alpha_i}\},\{\gamma_i^\prime\},
\{\xi^\prime_{ji}\})$, $i=0,1,\cdots,n$, satisfying
$\Pi\circ\tau^\prime=\tau$. Clearly $\gamma_n^\prime(1)$ lies in the 
image of the generalized system $\sigma$ representing the pseudo-embedding
$\tilde{i}$ (cf. (4.2.6)). Let $z$ be a point in the fiber
$Z$ such that $i(z)=\pi_{V_{\alpha_n}}(\gamma_n^\prime(1))$. Then it is
clear that $y\in Y$ is path-connected to $i(z)$ in $Y$ through path
$[\tau^\prime]$. This concludes the proof that the sequence
$\pi_0(Z,r)\stackrel{i_{\#}}{\rightarrow}\pi_0(Y,q)
\stackrel{\pi_{\#}}{\rightarrow}\pi_0(X,p)$ is exact at
$\pi_0(Y,q)$.

Now we take care of the rest of $(4.2.7)$. By Theorem 3.3.3, there is
a long exact sequence 
$$
\begin{array}{c}
\cdots\stackrel{j_{\#}}{\rightarrow}\pi_{k+1}(Y,Z,\tilde{i})
\stackrel{\partial^\prime}{\rightarrow}\pi_k(Z,\underline{r})\stackrel{i_{\#}}
{\rightarrow}\pi_k(Y,\underline{q})\stackrel{j_{\#}}{\rightarrow}
\pi_k(Y,Z,\tilde{i})\stackrel{\partial^\prime}{\rightarrow}\cdots\\
\stackrel{j_{\#}}{\rightarrow}\pi_1(Y,Z,\tilde{i})\stackrel{\partial^\prime}
{\rightarrow}\pi_0(Z,r)\stackrel{i_{\#}}{\rightarrow}\pi_0(Y,q)
\end{array}\leqno (4.2.8)
$$   
associated to the pseudo-embedding
$\tilde{i}:(Z,\underline{r})\rightarrow (Y,\underline{q})$. On the
other hand, the fibration $\tilde{\pi}$ induces a continuous map
$\pi_!:\Omega(Y,Z,\tilde{i})\rightarrow \Omega(X,\underline{p})$ given
by $\tilde{\gamma}\mapsto \tilde{\pi}\circ\tilde{\gamma}$, which
satisfies $\Omega(\tilde{\pi})=\pi_!\circ j$ where $j$ is the
canonical map from $\Omega(Y,\underline{q})$ to
$\Omega(Y,Z,\tilde{i})$. It is easily seen that $(4.2.7)$ is a
consequence of $(4.2.8)$ by setting $\partial=\partial^\prime\circ
(\pi_!)_{\#}^{-1}$ for the connecting homomorphisms
$\partial:\pi_k(X,\underline{p})\rightarrow
\pi_{k-1}(Z,\underline{r})$ for all $k\geq 1$, if we show that $\pi_!$
induces isomorphisms on homotopy groups
$(\pi_!)_{\#}:\pi_k(Y,Z,\tilde{i})\cong \pi_k(X,\underline{p})$ for
all $k\geq 1$.

\vspace{1.5mm}

\noindent{(1)} Surjectivity of
$(\pi_!)_{\#}:\pi_k(Y,Z,\tilde{i})\rightarrow\pi_k(X,\underline{p})$. It
suffices to show that given any continuous map
$u:(S^k,\ast)\rightarrow (\Omega(X,\underline{p}),\tilde{p})$ where
$k\geq 0$, there is a continuous map $\ell(u):(S^k,\ast)\rightarrow
(\Omega(Y,Z,\tilde{i}),\tilde{q})$ satisfying $\pi_!\circ
\ell(u)=u$. We take a simplicial decomposition of $S^k$,
$S^k=\cup_{a\in A} K_a$, such that the base point $\ast\in S^k$ is a
vertex of some simplex $K_o$. We may assume that the simplicial
decomposition is sufficiently fine so that the restriction of map $u$
to each simplex $K_a$ lies in a canonical neighborhood of a based loop
in $(X,\underline{p})$ (cf. Lemma 3.1.2). Hence there are systems 
$$
\tau_a=(\{K_a\times I_i\},\{U_i^a\},\{\gamma_i^a\},\{\xi_{ji}^a\}),
\hspace{2mm} i=0,1,\cdots,n, \leqno (4.2.9)
$$
representing the map $u|_{K_a}$. These systems satisfy the following
conditions. Each $U_i^a$ is an element of the set $\{U_\alpha\}$ in
the generalized system $\Pi$ representing the orbispace 
fibration $\tilde{\pi}$. Each
$\gamma_i^a:K_a\times I_i\rightarrow \widehat{U_i^a}$ is continuous in
both variables. Each $\xi_{ji}^a:K_a\rightarrow Tran(U_i^a,U_j^a)$ is
a constant map for $j=1,2,\cdots,n$, and a continuous map for $j=0$
and $i=n$. Moreover, $\gamma_0^a(x,0)=\hat{p}$ ($U_0^a=U_o$ for each
$a\in A$) for all $x\in K_a$, and $\gamma_j^a=\xi_{ji}^a\circ
\gamma_i^a$ on $K_a\times (I_i\cap I_j)$. Let $K_o$ be the simplex
containing the base point $\ast\in S^k$. We have
$\gamma_i^o(\ast,t)=\hat{p}$ for all $t\in I_i$, and
$\xi_{ji}^o(\ast)=1_{G_{U_o}}$ for all $j=0,\cdots,n$.  Let $K_{ab}$
be the common face of $K_a$ and $K_b$. Since the restriction of
$\tau_a$ and $\tau_b$ on $K_{ab}$ defines the same family of based
loops, there are transition maps $\eta_{ba}^i(x)\in Tran(U_i^a,U_i^b)$
for $x\in K_{ab}$, such that 
$$
\eta_{ba}^0(x)=1_{G_{U_o}}, \hspace{2mm}
\gamma_i^b(x,\cdot)=\eta_{ba}^i(x) \circ\gamma_i^a(x,\cdot), \hspace{2mm}
\xi_{ji}^b(x)=\eta_{ba}^j(x)\circ\xi_{ji}^a(x)\circ
(\eta_{ba}^i(x))^{-1} \leqno (4.2.10)
$$
for all $x\in K_{ab}$. It is easily seen from
these equations that $\eta_{ba}^i(x)$ are uniquely determined and
constant in $x$ for all $i=0,1,\cdots,n$. Hence on triple
intersections $K_a\cap K_b\cap K_c$, we have
$$
\eta_{ca}^i(x)=\eta_{cb}^i(x)\circ \eta_{ba}^i(x), \hspace{2mm}
i=0,1,\cdots, n. \leqno (4.2.11)
$$

We split $I_0$ into a union $I_{0,-}\cup I_{0,+}$ along $0\in I_0$ with
the convention that $I_{0,+}\cap I_1\neq \emptyset$. We set
$J_0=I_{0,+}$, $J_1=I_1, \cdots, J_n=I_n$ and $J_{n+1}=I_{0,-}$. Then
the systems $\tau_a$ can be regarded as
$$
\tau_a=(\{K_a\times J_k\},\{U_k^a\},\{\gamma_k^a\},\{\xi_{lk}^a\}),
\hspace{2mm} k=0,1,\cdots,n+1, \leqno (4.2.12)
$$
satisfying similar conditions, and we correspondingly have
$\eta_{ba}^k\in Tran(U_k^a,U_k^b)$ with $\eta_{ba}^{n+1}=1_{G_{U_o}}$.
Since $\gamma_0^a(x,0)=\hat{p}$ for all $x\in K_a$ and any $a\in A$,
and $\pi_o:\widehat{V_o}\rightarrow \widehat{U_o}$ is a weak
fibration, we can lift $\gamma_0^a: K_a\times
J_0\rightarrow\widehat{U_o}$ to a continuous map
$\ell(\gamma_0^a):K_a\times J_0\rightarrow \widehat{V_o}$ satisfying
$\ell(\gamma_0^a(x,0))=\hat{q}$ and $\pi_o\circ
\ell(\gamma_0^a)=\gamma_0^a$. We take the lifting
$\ell(\eta_{ba}^0)=1_{G_{V_o}}$ of $\eta_{ba}^0$, and require that
$\ell(\gamma_0^b)=\ell(\eta_{ba}^0)\circ \ell(\gamma_0^a)$. Now by
$b)$ of Definition 4.2.1, for each $U_1^a$, there is an index
$\alpha(1,a)\in \Lambda(U_1^a)$ such that
$\pi_{V_o}(\ell(\gamma_0^a)(J_0\cap J_1))$ lies in
$V_{\alpha(1,a)}$. These basic open sets have the property that
$V_{\alpha(1,a)}\cap V_{\alpha(1,b)}\neq \emptyset$ if $K_a\cap
K_b\neq \emptyset$.

Recall that the base point $\ast\in S^k$ is contained in $K_o$ and
$\gamma_i^o(\ast,t)=\hat{p}$ for all $t\in I_i$, and
$\xi_{ji}^o(\ast)=1_{G_{U_o}}$ for all $j=0,\cdots,n$. We take
arbitrary liftings $\ell(\xi_{10}^a)$ of $\xi_{10}^a$, only requiring 
$\ell(\xi_{10}^o)=1_{G_{V_o}}$ (the existence of $\ell(\xi_{10}^a)$ is
ensured by the assumption that
$\rho_{\alpha(1,a)o}:Tran(V_o,V_{\alpha(1,a)}) \rightarrow
Tran(U_o,U_1^a)$ is surjective). We define $\ell(\eta_{ba}^1)$ by
$$
\ell(\eta_{ba}^1):=\ell(\xi_{10}^b)\circ\ell(\eta_{ba}^0)\circ 
(\ell(\xi_{10}^a))^{-1}, \leqno (4.2.13)
$$
which satisfy the following compatibility conditions
$$
\ell(\eta_{ca}^1)=\ell(\eta_{cb}^1)\circ\ell(\eta_{ba}^1) \leqno (4.2.14)
$$
on triple intersections $K_a\cap K_b \cap K_c$, since $(4.2.14)$ holds
for $\ell(\eta_{ba}^0)$ trivially. Now we define
$\ell(\gamma_1^a):=\ell(\xi_{10}^1)\circ\ell(\gamma_0^a)$ as a map
$K_a\times (J_0\cap J_1)\rightarrow \widehat{V_{\alpha(1,a)}}$, and
want to extend it to a map $\ell(\gamma_1^a): K_a\times J_1\rightarrow
\widehat{V_{\alpha(1,a)}}$ satisfying $\pi_{\alpha(1,a)}\circ
\ell(\gamma_1^a)=\gamma_1^a$. The extension is carried out step by
step as follows. For any vertex $x$, choose a simplex $K_a$ containing
$x$. Since $\pi_{\alpha(1,a)}:\widehat{V_{\alpha(1,a)}}\rightarrow
\widehat{U_1^a}$ is a weak fibration, $\ell(\gamma_1^a)(x,t)$, $t\in
J_1$, can be defined as a lifting of $\gamma_1^a(x,t)$. If $x$ is also
contained in another simplex $K_b$, we define $\ell(\gamma_1^b)(x,t):=
\ell(\eta_{ba}^1)(\ell(\gamma_1^a(x,t)))$ for $t\in J_1$. 
The compatibility conditions in $(4.2.14)$ ensure that this is
well-defined. Now we appeal to a general fact that for any simplex 
$\sigma$ and the interval $I=[0,1]$, there is a self-diffeomorphism of
$\sigma\times I$ such that $\sigma\times\{0\}\cup
(\partial\sigma\times I)$ is mapped to $\sigma\times\{0\}$. 
It is easily seen that $\ell(\gamma_1^a)$ can be defined over 
edges cross $J_1$, etc., and finally over $K_a\times J_1$.

We can repeat this procedure and define $V_{\alpha(k+1,a)}$,
$\ell(\xi_{(k+1)k}^a)$, $\ell(\eta_{ba}^{k+1})$, and
$\ell(\gamma_{k+1}^a)$ whenever they are already defined for $k$. At
the last step $k=n$, we have to use the assumption that 
$\rho_{\alpha(n+1,a)\alpha(n,a)}:Tran(V_{\alpha(n,a)},V_{\alpha(n+1,a)})
\rightarrow Tran(U_n^a,U_o)$ is a surjective weak fibration, instead
of merely a surjective map as required in the previous steps. We also
observe that the liftings $\ell(\eta_{ba}^{n+1})$ defined by
$$
\ell(\eta_{ba}^{n+1}):=\ell(\xi_{(n+1)n}^b)\circ\ell(\eta_{ba}^n)\circ 
(\ell(\xi_{(n+1)n}^a))^{-1}\leqno (4.2.15)
$$
are continuous maps from $K_a\cap K_b$ to 
$Tran(V_{\alpha(n+1,a)},V_{\alpha(n+1,b)})$, which lie in the image of
the generalized system $\sigma$ representing the pseudo-embedding 
$\tilde{i}$ (cf. $(4.2.6)$).  The points $\ell(\gamma_{n+1}^a(x,1))$ also 
lie in the image of $\sigma$ for all $x\in K_a$. Thus we have
obtained a family of generalized systems
$$
\ell(\tau_a)=(\{K_a\times J_k\},\{V_{\alpha(k,a)}\},
\{\ell(\gamma_k^a)\},\{\ell(\xi_{lk}^a)\}), \hspace{2mm} 
k=0,1,\cdots,n+1, \leqno (4.2.16)
$$
which are patched together by continuous families of transition maps
$\ell(\eta_{ba}^k)(x)$, $x\in K_a\cap K_b$, to define a continuous map 
$\ell(u):(S^k,\ast)\rightarrow (\Omega(Y,Z,\tilde{i}),\tilde{q})$ 
satisfying $\pi_!\circ\ell(u)=u$. This concludes the proof of (1).

\vspace{1.5mm}

\noindent{(2)} Injectivity of $(\pi_!)_{\#}$. We need to show that
given a continuous map $u:(S^k,\ast)\rightarrow  
(\Omega(Y,Z,\tilde{i}),\tilde{q})$, if there is a continuous map 
$H:(CS^k,\ast)\rightarrow (\Omega(X,\underline{p}),\tilde{p})$ such
that the restriction of $H$ to $S^k\subset CS^k$ equals $\pi_!\circ
u$, then there is a continuous map $\ell(H):(CS^k,\ast)\rightarrow
(\Omega(Y,Z,\tilde{i}),\tilde{q})$ satisfying $\pi_!\circ \ell(H)=H$
and $\ell(H)|_{S^k}=u$. We can represent $H$ by a collection of
systems as in $(4.2.9)$, and repeat the construction in (1) for $H$. 
The only difference occured here is that the lifting $\ell(H)$ is already
determined on the subset $S^k$ of $CS^k$. But this does not effect the
argument. 

The proof of Theorem 4.2.5 is completed.

\hfill $\Box$

We end the discussion of this section with some examples of orbispace
fibration. 

\vspace{1.5mm}

\noindent{\bf Example 4.2.6:} (orbispace fiber bundles)

Let $\tilde{\pi}:Y\rightarrow X$ be an orbispace fibration,
$\Pi=(\{V_\alpha\},\{U_\alpha\},\{\pi_\alpha\},\{\rho_{\beta\alpha}\})$,
$\alpha\in \Lambda$, be a generalized system representing
$\tilde{\pi}$. We will call $\tilde{\pi}:Y\rightarrow X$
an {\it orbispace fiber bundle} over $X$ if for each $\alpha$ 
there is a locally connected
topological space $F_\alpha$ such that $\widehat{V_\alpha}
=\widehat{U_\alpha}\times F_\alpha$ with $\pi_\alpha$ 
being the natural projection, $G_{V_\alpha}=G_{U_\alpha}$ with 
$\rho_\alpha: G_{V_\alpha}\rightarrow G_{U_\alpha}$ being the identity
homomorphism, and for any $U\in\{U_\alpha\}$, $V_\alpha\cap
V_\beta=\emptyset$ if $\alpha,\beta\in\Lambda(U)$ are distinct. The
fiber of $\tilde{\pi}$ over a base-point structure $\underline{o}=
(o,U,\hat{o})$ for any $U\in\{U_\alpha\}$ is the topological space
$\sqcup_{\alpha\in\Lambda(U)} F_\alpha$, whose homeomorphism class
depends only on the connected component of $X$ that contains $U$.

For example, let $X$ be an orbifold and $E$ be an orbifold
vector bundle of rank $n$ over $X$. Then $E$ is an orbispace fiber
bundle over $X$ with fiber $\R^n$. Moreover, one can derive many
orbispace fiber bundles from $E$. For instance, the total space $Y$ of
orthonormal frames of $E$ (assuming some metric on $E$) is naturally 
an orbifold, and is also an orbispace fiber bundle over $X$ under the 
natural projection. The fiber of $Y$ is $O(n)$. 

Let $X$ be an orbifold of dimension $n$ and $E=TX$ be the
tangent bundle of $X$. Then the fiber bundle $Y$ of orthonormal frames
of $E$ is actually a smooth manifold since over
each uniformizing system $(V_p,G_p,\pi_p)$ of $X$, $p\in X$, the
action of $G_p$ on $V_p$ is effective. Furthermore, there is a natural
right action of $O(n)$ on $Y$, $(y,A)\mapsto yA$, such that 
the orbit space $Y/O(n)$ is identified with $X_{top}$ under the
natural projection. The homotopy groups of the orbifold $X$ defined by
Haefliger in \cite{Ha1} are isomorphic to the homotopy groups of the
Borel space of the $G$-space $(Y,O(n))$ (cf. \cite{Boy}). As a natural
application of the exact sequence $(4.2.7)$, we shall
prove next that the definition of homotopy groups of orbifolds in this
paper is equivalent to Haefliger's definition in \cite{Ha1}.

\hfill $\Box$

\noindent{\bf Theorem 4.2.7:} 
{\em Let $X$ be a connected orbifold of dimension $n$, $Y$ be the bundle of
orthonormal frames on $X$ with the standard right action of
$O(n)$. Then for any $k\geq 1$, $\pi_k(X)$ is isomorphic to the $k$-th
homotopy group of the Borel space of $(Y,O(n))$.
}

\vspace{2mm}

\noindent{\bf Proof:}
We regard the right action of $O(n)$ on $Y$ as a left action:
$O(n)\times Y\rightarrow Y$ by $(A,y)\mapsto yA^{-1}$. Since
$Y/O(n)=X_{top}$, we can put another orbispace structure on $X_{top}$
which is canonically defined from the $G$-space $(Y,O(n))$ (with the left
action). We denote the resulting orbispace by $X^\prime$. By Theorem
3.4.1, it suffices to show that $\pi_k(X)$ is isomorphic to
$\pi_k(X^\prime)$ for all $k\geq 1$. 

We shall first choose a generalized system $\Pi$ to represent the
orbispace fibration $\tilde{\pi}:Y\rightarrow X$, which induces the
natural projection $\pi:Y\rightarrow X$. Let $\{U_i\}$ be a
collection of geodesic neighborhoods which covers $X$, and for each
$i$, $(V_i,G_i,\pi_i)$ be the uniformizing system of $U_i$ where $V_i$
is a geodesic ball. Then each $Y_i:=\pi^{-1}(U_i)\subset Y$ is an open
smooth submanifold and there is a finite covering map
$V_i\times O(n)\rightarrow Y_i$ given by the natural
action of $G_i$ on $V_i\times O(n)$. Each $Y_i$ is either connected or
of two connected components, depending on whether $G_i$ acts on $V_i$ 
orientation reversingly or orientation preservingly. 
Let $\{Y_i^\nu\}$ be the set of
components of $Y_i$. We set $\widehat{Y}_i^\nu=V_i\times
O(n)$ if $Y_i$ is connected, and set $\widehat{Y}_i^\nu=V_i\times
SO(n)$ if $Y_i$ is disconnected. Then we define 
$(\widehat{Y}_i^\nu,G_i)$ to be the generalized $G$-structure of
$Y_i^\nu$, and define $\pi^{i,\nu}:\widehat{Y}_i^\nu\rightarrow V_i$ to be 
the projection onto the first factor. Observe
that each transition map $\xi\in Tran(U_i,U_j)$ induces a transition
map $\xi^\prime\in Tran(Y_i^{\nu_i},Y_j^{\nu_j})$ between the generalized
$G$-structures $(\widehat{Y}_i^{\nu_i},G_i)$ and 
$(\widehat{Y}_j^{\nu_j},G_j)$ whenever $Y_i^{\nu_i}\cap Y_j^{\nu_j}
\neq \emptyset$, and we define $\rho_{ji,\nu}:\xi^\prime\mapsto \xi$ 
to be the corresponding map $Tran(Y_i^{\nu_i},Y_j^{\nu_j})
\rightarrow Tran(U_i,U_j)$ which is one to one and onto,
and restricts to the identity isomorphism $\rho_i:G_i\rightarrow
G_i$. We set $\Pi:=(\{Y_i^\nu\},\{U_i\},\{\pi^{i,\nu}\},\{\rho_{ji,\nu}\})$.

Fix a base point $o\in X$ which lies in some $U_o\in\{U_i\}$ for which
$G_o$ is trivial. Let $\underline{o}=(o,U_o,o)$ be the associated
base-point structure of $X$. Fix an inverse image $o^\prime\in
\pi^{-1}(o)$ in $Y$, and set $\underline{o}^\prime=(o,Y_o,o^\prime)$ for
the corresponding base-point structure of $X^\prime$, where $Y_o$ is
the connected component of $Y$ containing $o^\prime$. We shall
construct a continuous map
$\phi:\Omega(X^\prime,\underline{o}^\prime)\rightarrow
\Omega(X,\underline{o})$ as follows. Let $u=(\gamma,A)$ be an element
of $\Omega(X^\prime,\underline{o}^\prime)$, where $A\in O(n)$ and
$\gamma: [0,1]\rightarrow Y$ such that $\gamma(0)=o^\prime$ and
$\gamma(1)=o^\prime A^{-1}$. We set $I_{\nu_i}=\gamma^{-1}(Y_i^{\nu_i})$, 
and for each $\nu_i$ pick a lifting $\gamma_{\nu_i}:
I_{\nu_i}\rightarrow \widehat{Y}_i^{\nu_i}$ of
$\gamma|_{I_{\nu_i}}$. Then there is a unique $\xi_{\nu_j\nu_i}\in 
Tran(Y_i^{\nu_i},Y_j^{\nu_j})$
such that $\gamma_{\nu_j}=\xi_{\nu_j\nu_i}\circ\gamma_{\nu_i}$. 
The generalized system 
$\sigma=(\{I_{\nu_i}\},\{Y_i^{\nu_i}\},\{\gamma_{\nu_i}\},
\{\xi_{\nu_j\nu_i}\})$ defines a based
path in $Y$, which is independent of the choices of the liftings
$\gamma_{\nu_i}$. We define 
$\phi:\Omega(X^\prime,\underline{o}^\prime)\rightarrow
\Omega(X,\underline{o})$ by setting $\phi(u)=\tilde{\pi}\circ
[\sigma]$. 

The exact sequences $(3.4.3)$ and $(4.2.7)$ can be put together to
form the following commutative diagram 
$$
\begin{array}{ccccccccc}
\cdots \rightarrow & \pi_k(O(n),I_n) &
\stackrel{\partial}{\rightarrow} & \pi_k(Y,o^\prime) &
\stackrel{i_\#}{\rightarrow} & \pi_k(X^\prime,\underline{o}^\prime) & 
\stackrel{\pi_\#}{\rightarrow} & \pi_{k-1}(O(n),I_n) & \rightarrow
\cdots\\
         &\downarrow\lambda_\ast &     & \parallel  &
& \downarrow \phi_\# &  & \downarrow\lambda_\ast &   \\
\cdots \rightarrow & \pi_k(O(n),I_n) &
\stackrel{i_\#}{\rightarrow} & \pi_k(Y,o^\prime) &
\stackrel{\pi_\#}{\rightarrow} & \pi_k(X,\underline{o}) & 
\stackrel{\partial}{\rightarrow} & \pi_{k-1}(O(n),I_n) & \rightarrow
\cdots
\end{array} \leqno (4.2.17)
$$
where $\lambda_\ast:\pi_k(O(n),I_n)\rightarrow \pi_k(O(n),I_n)$ is
the induced isomorphism of the map $\lambda:O(n)\rightarrow O(n)$
defined by $A\mapsto A^{-1}$. The commutative diagram $(4.2.17)$
implies that $\phi_\#:\pi_k(X^\prime,\underline{o}^\prime)\rightarrow
\pi_k(X,\underline{o})$ is an isomorphism for all $k\geq 1$. This
completes the proof.

\hfill $\Box$

\noindent{\bf Example 4.2.8:} (Seifert fibered spaces)

Let $Y$ be a 3-manifold which is a Seifert fibered space. Then $Y$ is an
orbispace fiber bundle over an orbifold Riemann surface with fiber
$S^1$. In this example we examine the associated exact sequence
$(4.2.7)$ for the case when $Y$ is the lens space $L(p,q)$ where $p,q$
are relatively prime. 

The lens space $Y=L(p,q)$ can be exhibited as the orbit space of a
free $\Z_p$ action on $S^3\subset \C^2$. The $\Z_p$ action is defined by
$e^{2\pi i/p}\cdot (z_1,z_2)=(e^{2\pi i/p}z_1,e^{2\pi qi/p}z_2)$. The
Hopf fibration $S^3\rightarrow S^2$ is preserved by the $\Z_p$ action,
hence induces a projection $\pi:Y\rightarrow X$, where $X$ is
$S^2$. In fact $\pi:Y\rightarrow X$ is a Seifert fibration. Let us
write $\frac{q-1}{p}=\frac{n}{m}$ where $n,m$ are relatively prime,
and let $l=gcd(p,q-1)$ hence $p=lm, q-1=ln$. The manifold $Y$ can be
written as $(Y_1\cup_\phi Y_2)/\Z_m$ where each $Y_j$, $j=1,2$, is
$S^1\times D^2$, $\Z_m$ acts on $Y_1$ by $e^{2\pi i/m}\cdot
(w,z)=(e^{2\pi i/m}w,e^{2\pi ni/m}z)$ and on $Y_2$ by $e^{2\pi
i/m}\cdot (w,z)=(e^{2\pi qi/m}w,e^{-2\pi ni/m}z)$, and $\phi:S^1\times
\partial D^2\rightarrow S^1\times \partial D^2$ is the $\Z_m$-equivariant
diffeomorphism defined by $(w,z)\mapsto (wz^l,z^{-1})$. Hence $X$ has
an orbifold structure defined by the $\Z_m$ action on $S^2$: $e^{2\pi
i/m}\cdot [z_1,z_2]=[e^{2\pi i/m}z_1,z_2]$. The exact sequence
$(4.2.7)$ associated to $\pi:Y\rightarrow X$ is 
$$
0\rightarrow \pi_2(X)\stackrel{\partial}{\rightarrow}
\pi_1(S^1)\rightarrow \pi_1(Y)\rightarrow \pi_1(X)\rightarrow 0
\leqno (4.2.18)
$$
where $\pi_2(X)=\Z$, $\pi_1(S^1)=\Z$, $\pi_1(Y)=\Z_p$ and
$\pi_1(X)=\Z_m$. The connecting homomorphism $\partial:\Z\rightarrow
\Z$, which must be given by $x\mapsto lx$ for any $x\in\Z$, has the
following geometric interpretation. A generator of $\pi_2(X)$ can be
represented by a $S^1$-family of based loops $\gamma_s$ in $S^2$,
$s\in S^1$, where if we identify $S^2$ with the suspension of $S^1$,
each $\gamma_s$ goes around $S^1$ $m$ times. One lifts $\gamma_s$ to
$Y$ to get a $S^1$-family of paths $\gamma_s^\prime$ such that the
terminal point of $\gamma_s^\prime$ lies in the fiber over the base
point in $X$, which is $S^1$. The connecting homomorphism $\partial$
is defined by the degree of the map $s\mapsto 
\gamma_s^\prime(1)$. The degree of this map is $l$, as seen in the
gluing map $\phi: (w,z)\rightarrow (wz^l,z^{-1})$.

\hfill $\Box$

\noindent{\bf Example 4.2.9:} (normal orbispaces)

Recall the definition of normal orbispaces in Remark 2.1.4 e. Let $X$
be an orbispace. For any basic open set $U$, let $K_U=\{g\in
G_U|g\cdot x=x,\hspace{1.5mm}\forall x\in\widehat{U}\}$. The orbispace
$X$ is called normal if for any $W\subset U$, $K_W$ is mapped
isomorphically to $K_U$ under any injection in $Tran(W,U)$. 
The canonical reduction
$X^{red}$ of a normal orbispace $X$ is defined as follows. The
underlying topological space $(X^{red})_{top}$ is still $X_{top}$. The
$G$-structure on $U$, however, is changed to
$(\widehat{U},G_U/K_U,\pi_U)$. There is a canonical morphism
$\tilde{\pi}:X\rightarrow X^{red}$ defined by the system
$\Pi=(\{U_\alpha\},\{U_\alpha\},\{\pi_\alpha\},\{\rho_{\beta\alpha}\})$
where $\{U_\alpha\}$ is a cover of $X$, $\pi_\alpha:\widehat{U_\alpha}
\rightarrow \widehat{U_\alpha}$ is the identity map, and
$\rho_{\beta\alpha}$ sends each transition map $(\phi,\lambda)$ in
$Tran(U_\alpha,U_\beta)$ to the transition map $(\phi,[\lambda])$
where $[\lambda]$ is the isomorphism on the quotient groups induced by
$\lambda$. It is easily seen that $\tilde{\pi}:X\rightarrow X^{red}$
is an orbispace fibration if and only if each projection 
$G_{U_\alpha}\rightarrow G_{U_\alpha}/K_{U_\alpha}$ is a weak
fibration. For a connected normal
orbispace $X$, the abstract group which is isomorphic to all $K_U$ is
denoted by $K_X$ and called the kernel of $X$.
Without loss of generality, we assume that the orbispace $X$ is
connected. Then the fiber over any base-point structure $\underline{p}$ of 
$X^{red}$ is the global orbispace defined by the $G$-space 
$(\{pt\},K_X)$. In this case, the associated exact sequence $(4.2.7)$ 
becomes
$$
\begin{array}{c}
\cdots\stackrel{\pi_{\#}}{\rightarrow}\pi_{k+1}(X^{red},\underline{p})
\stackrel{\partial}{\rightarrow}\pi_{k-1}(K_X,1)\stackrel{i_{\#}}
{\rightarrow}\pi_k(X,\underline{q})\stackrel{\pi_{\#}}{\rightarrow}\\
\pi_k(X^{red},\underline{p})\stackrel{\partial}{\rightarrow} \cdots
\stackrel{\pi_{\#}}{\rightarrow}\pi_1(X^{red},\underline{p})
\rightarrow 1. 
\end{array} \leqno (4.2.19)
$$   
  
Let us examine the example in Remark 2.1.4 e, where
$X=(S^1,\U_\tau)$. The canonical reduction $X^{red}$ is the trivial
orbispace $S^1$, and the kernel $K_X=G$. The exact sequence $(4.2.19)$
becomes 
$$
1\rightarrow \pi_0(G)\rightarrow \pi_1(X)\rightarrow \Z\rightarrow
1. \leqno (4.2.20)
$$
The monomorphism $\pi_0(G)\rightarrow \pi_1(X)$ in $(4.2.20)$ can be
realized as follows. Each $g\in G$ canonically defines a based loop
$\tilde{\gamma}_g$ in $X$, where $\tilde{\gamma}_g$ is defined by
$(\{I_0,I_1\},\{U_1,U_1\},\{\gamma_0,\gamma_1\},\{\xi_{10}\})$ in
which $\gamma_0,\gamma_1$ are constant maps into the base point and
$\xi_{10}=g^{-1}$. The monomorphism $\pi_0(G)\rightarrow \pi_1(X)$ is given
by $[g]\mapsto [\tilde{\gamma}_g]$. On the other hand, there is a
based loop $\tilde{f}$ in $X$ which is defined by the system
$(\{I_0,I_1\},\{U_2,U_1\},\{f_0,f_1\},\{\eta_{10}\})$, where
$f_0$ and $f_1$ are homeomorphisms, and $\eta_{10}$ is
defined by the isomorphism $\tau:G\rightarrow G$. The class
$[\tilde{f}]$ is mapped to a generator of $\Z$ under the
epimorphism $\pi_1(X)\rightarrow \Z$ in $(4.2.20)$. Now it is easy to
see that $\nu(\tilde{f})\#\tilde{\gamma}_g\#\tilde{f}$
is homotopic to the based loop $\tilde{\gamma}_{\tau(g)}$. Hence the
exact sequence $(4.2.20)$ identifies $\pi_1(X)$ as the semi-direct
product of $\pi_0(G)$ by $\Z$ with respect to the homomorphism
$\Z\rightarrow Aut(\pi_0(G))$ given by $1\mapsto \tau_\ast$. 

\subsection{Seifert-Van Kampen Theorem}

\hspace{5mm}
In this section, we generalize to the case of orbispaces the classical 
Seifert-Van Kampen Theorem expressing the fundamental group of the union of 
two open sets in terms of the fundamental groups of each open set and 
their intersection\footnote{in the case of orbifolds it has been addressed
in \cite{Pr} via a different setup}. The proof of our theorem is a 
suitable modification of 
the one given in \cite{Ma} for the case of topological spaces.

\vspace{2mm}

\noindent{\bf Theorem 4.3.1:}
{\em Suppose $X_1$ and $X_2$ are two open sub-orbispaces of an orbispace 
$X$ such that $X=X_1\cup X_2$ and $X_1\cap X_2\neq\emptyset$, and 
$X_1$, $X_2$ and $X_1\cap X_2$ are path-connected. Let $\underline{o}$ be a 
common base-point structure of $X_1$, $X_2$ and $X_1\cap X_2$. For any group 
$H$ which satisfies the following commutative diagram 
$$
\begin{array}{ccccc}
   &   & \pi_1(X_1,\underline{o}) &    &   \\
   &\stackrel{j_1}{\nearrow} &        &\stackrel{\rho_1}{\searrow} & \\
\pi_1(X_1\cap X_2,\underline{o}) &  &\stackrel{\rho_{12}}{\longrightarrow} &  & H \\
   &\stackrel{j_2}{\searrow} &        &\stackrel{\rho_2}{\nearrow} &  \\
   &   & \pi_1(X_2,\underline{o}) &    &   
\end{array} \leqno (4.3.1)
$$
for some homomorphisms $\rho_1,\rho_2,\rho_{12}$, where 
$j_1,j_2$ are induced by the corresponding inclusions between
orbispaces, then there exists a homomorphism 
$\sigma:\pi_1(X,\underline{o})\rightarrow H$ such that 
$$
\begin{array}{ccc}
   &   & \pi_1(X,\underline{o}) \\
   &\stackrel{i_\alpha}{\nearrow} &  \\
\pi_1(X_\alpha,\underline{o}) &   &\downarrow \sigma \\
   &\stackrel{\rho_\alpha}{\searrow} &  \\
   &   &   H 
\end{array} \leqno (4.3.2)
$$
is commutative, where $i_\alpha$ is induced by inclusion 
$X_\alpha\hookrightarrow X$, for $\alpha=1,2,12$. Here $X_{12}$ stands 
for $X_1\cap X_2$.
}

\vspace{2mm}

\noindent{\bf Proof:}
First of all, since $X_1\cap X_2$ is path-connected, we can prove as in the 
classical case that $\pi_1(X,\underline{o})$ is generated by the subgroups 
$i_\alpha (\pi_1(X_\alpha,\underline{o}))$ where $\alpha=1,2$. 
For any $\beta\in \pi_1(X,\underline{o})$, we write 
$$
\beta=i_{\alpha_1}(\beta_1)\cdot i_{\alpha_2}(\beta_2)
\cdots i_{\alpha_k}(\beta_k) \leqno (4.3.3)
$$ 
where $\alpha_i$ is either $1$ or $2$ while $\beta_i$ is correspondingly in
$\pi_1(X_1,\underline{o})$ or $\pi_1(X_2,\underline{o})$. We define 
the homomorphism $\sigma:\pi_1(X,\underline{o})\rightarrow H$ by
$$
\sigma(\beta):=\rho_{\alpha_1}(\beta_1)\cdot\rho_{\alpha_2}(\beta_2)
\cdots\rho_{\alpha_k}(\beta_k). \leqno (4.3.4)
$$
The diagram $(4.3.2)$ is automatically satisfied once we verify that the 
so-defined $\sigma$ is a well-defined homomorphism. This boils down to show 
that if $i_{\alpha_1}(\beta_1)\cdot i_{\alpha_2}(\beta_2)\cdots 
i_{\alpha_k}(\beta_k)=1$ in $\pi_1(X,\underline{o})$, then 
$\rho_{\alpha_1}(\beta_1)\cdot\rho_{\alpha_2}(\beta_2)\cdots\rho_{\alpha_k}
(\beta_k)=1$ in $H$.

Represent each $\beta_j$ by an element $\tilde{v}_j$ in 
$\Omega(X_{\alpha_j},\underline{o})$. Then 
$i_{\alpha_1}(\beta_1) \cdot i_{\alpha_2}(\beta_2)\cdots 
i_{\alpha_k}(\beta_k)$ is represented by 
$\tilde{\gamma}_0=\tilde{v}_1\# \tilde{v}_2\# \cdots\tilde{v}_k$ in 
$\Omega(X,\underline{o})$. The assumption that $i_{\alpha_1}(\beta_1)
\cdot i_{\alpha_2}(\beta_2)\cdots i_{\alpha_k}(\beta_k)=1$ in 
$\pi_1(X,\underline{o})$ means that there is a continuous map $f:[0,1]
\rightarrow \Omega(X,\underline{o})$ such that $f(0)=\tilde{o}$ and 
$f(1)=\tilde{\gamma}_0$. We can subdivide $[0,1]$ into $\cup_{k=0}^{m-1}
[s_k,s_{k+1}]$ and assume that each $f([s_k,s_{k+1}])$ is in a canonical
neighborhood of a based loop (cf. Lemma 3.1.2) so that it is represented by
a system
$$
\tau_k=(\{[s_k,s_{k+1}]\times I_i\},\{U_i^k\},\{f_i^k\},\{\xi_{ji}^k\}),
\hspace{2mm} i=0,1,\cdots,n,\leqno (4.3.5)
$$
where each $U_i^k$ is either contained in $X_1$ or in $X_2$, and the $\{I_i\}$
can be chosen to be independent of $k$ by passing to refinements. Each $f_i^k:
[s_k,s_{k+1}]\times I_i\rightarrow \widehat{U_i^k}$ is continuous, satisfying
$f_0^k(s,0)=f_i^0(0,t)=\hat{o}$ for any $k$, $i$, and $s\in [s_k,s_{k+1}]$
and $t\in I_i$. Each $\xi_{ji}^k:[s_k,s_{k+1}]\rightarrow Tran(U_i^k,U_j^k)$
is continuous and constant in $s\in [s_k,s_{k+1}]$ when $j=1,\cdots,n$. 
Moreover, $f_j^k(s,t)=\xi_{ji}^k(s)(f_i^k(s,t))$ for any $(s,t)\in
[s_k,s_{k+1}]\times (I_i\cap I_j)$. Restricted to each $s_k$, 
$k=1,\cdots,m-1$, since 
$\tau_{k-1}$ and $\tau_k$ define the same based loop $f(s_k)$, there are
transition maps $\eta_{k(k-1)}^i\in Tran(U_i^{(k-1)},U_i^k)$ such that
$$
\begin{array}{c}
\eta_{k(k-1)}^0=1_{G_{U_o}}, \hspace{3mm}
f^k_i(s_k,\cdot)=\eta_{k(k-1)}^i\circ f^{(k-1)}_i(s_k,\cdot),\\
\xi_{ji}^k(s_k)=\eta_{k(k-1)}^j\circ\xi_{ji}^{(k-1)}(s_k)\circ
(\eta_{k(k-1)}^i)^{-1}.
\end{array} \leqno (4.3.6)
$$
We choose $t_0=0$, $t_1\in I_0\cap I_1,\cdots,t_n\in I_{n-1}\cap I_n$, and 
$t_{n+1}=I_n\cap I_0$, $t_{n+2}=0\in I_0$. Without loss of generality, 
we may assume that the end points of the domain of each morphism 
$\tilde{v}_j$ are included in the set $\{t_0,\cdots,t_{n+2}\}$. 

Before proceeding further with the proof, we need to introduce some notations.
Set $v_{k,i}^{-}=f_{i-1}^{(k-1)}(s_k,t_i)\in \widehat{U}_{i-1}^{(k-1)}$ and 
$v_{k,i}^{+}=f_i^{(k-1)}(s_k,t_i)\in \widehat{U}_{i}^{(k-1)}$ for $k\geq 0$
and $i\geq 0$ with the understanding that $i-1=0$ when $i=0$ and $k-1=0$ 
when $k=0$. Set $I_{k,i}=\{s_k\}\times [t_{i-1},t_i]$ for $k\geq 0$, $i\geq 1$,
and $J_{k,i}=[s_{k-1},s_k]\times\{t_i\}$ for $k\geq 1$ and $i\geq 0$. Since
each $X_1$, $X_2$ and $X_1\cap X_2$ is path-connected, for each of 
$v_{k,i}^{-}$, we can choose a system $g_{k,i}^{-}$ from $[0,1]$ into 
$X_\alpha$, where $\alpha=1,2,12$, whose initial point is $\hat{o}$ and 
ending point is $v_{k,i}^{-}$. We require that $\alpha=12$ if the image of 
$v_{k,i}^{-}$ in $X$ is contained in $X_1\cap X_2$. Observe that $v_{k,i}^{-}=
v_{k,i}^{+}$ except for $i=n+1$, since $\xi_{0n}^{(k-1)}(s)$ may not be 
constant in $s$. Although each $\xi_{0n}^{(k-1)}(s_{k-1})\circ f_{n}^{(k-1)}$ 
may not coincide with $f_0^{(k-1)}$ when restricted on $J_{k,n+1}$, 
they are canonically homotopic by Sublemma 3.3.5.
We let $g_{k,i}^{+}$ be the system from $[0,1]$ into $X_\alpha$, whose 
initial point is $\hat{o}$ and ending point is $v_{k,i}^{+}$, such that
$g_{k,i}^{+}=g_{k,i}^{-}$ for $i=0,1,\cdots,n$, and $g_{k,n+1}^{+}$ equals
the composition of $g_{k,n+1}^{-}$ with the map traced out by the canonical 
homotopy between $\xi_{0n}^{(k-1)}(s_{k-1})\circ f_{n}^{(k-1)}$ and 
$f_0^{(k-1)}$ at the end point $(s_k,t_{n+1})$ of $J_{k,n+1}$, 
and introduce $g_{k,n+2}^{-}=g_{k,n+2}^{+}=\tilde{o}$ for convenience.
We let $a_{k,i}$ be the system obtained by composing $g_{k,i-1}^{+}$ with 
$f_{i-1}^{(k-1)}|_{I_{k,i}}$ followed by $(g_{k,i}^{-})^{-1}$, and 
$b_{k,i}^{-}$ be the system obtained by composing $g_{k-1,i}^{-}$ with 
$f_{i-1}^{(k-1)}|_{J_{k,i}}$ followed by $(g_{k,i}^{-})^{-1}$, and 
$b_{k,i}^{+}$ be the system obtained by composing $g_{k-1,i}^{+}$ with 
$f_{i}^{(k-1)}|_{J_{k,i}}$ followed by $(g_{k,i}^{+})^{-1}$. 

We have the following relations amongst the systems introduced in the
previous paragraph:
$$
a_{k-1,i}\# b_{k,i}^{-}\sim_\alpha b_{k,i-1}^{+}\# a_{k,i}, \hspace{1cm}
b_{k,i}^{-}\sim_\alpha b_{k,i}^{+}, \leqno (4.3.7)
$$
where $\sim_\alpha$ stands for ``homotopy equivalent to'' in 
$\Omega(X_\alpha,\underline{o})$, provided that $U_{i-1}^{(k-1)}$ is contained
in $X_\alpha$ for some $\alpha=1,2,12$.

Now by $(4.3.1)$, if a based loop  
$\tilde{\gamma}\in\Omega(X_1\cap X_2,\underline{o})$ represents $x_i$ in 
$\pi_1(X_i,\underline{o})$ for $i=1,2$, then $\rho_1(x_1)=\rho_2(x_2)$ in 
$H$. This enables us to assign each homotopy class of based loops in 
$\Omega(X_\alpha,\underline{o})$, $\alpha=1,2, 12$, with an element in $H$, 
independent of which $\Omega(X_\alpha,\underline{o})$ it is considered in. 
We denote this assignment by $\rho$. To finish the proof,
we simply observe that by $(4.3.7)$, we have 
\begin{eqnarray*}
\rho_{\alpha_1}(\beta_1)
\cdot\rho_{\alpha_2}(\beta_2)\cdots\rho_{\alpha_k}(\beta_k)
& = & \rho(a_{m,1})\cdot \rho(a_{m,2})\cdots\rho(a_{m,n+2})\\
& = & \cdots \\
& = & \rho(a_{0,1})\cdot \rho(a_{0,2})\cdots\rho(a_{0,n+2})\\
& = & \rho([\tilde{o}])\cdot\rho([\tilde{o}])\cdots\rho([\tilde{o}])\\
& = & 1
\end{eqnarray*}
\hfill $\Box$

\end{document}